\documentclass[12pt,a4paper]{article}
\usepackage{amssymb}
\usepackage{amsmath}
\usepackage{a4wide}
\newtheorem{theorem}{Theorem}
\newtheorem{corollary}{Corollary}
\newtheorem{lemma}{Lemma}
\newtheorem{proposition}{Proposition}
\newtheorem{definition}{Definition}

\newtheorem{conjecture}{Conjecture}

\newtheorem{thkh}{Khintchine's Theorem (1924).}
\newtheorem{thjar}{Jarn\'{\i}k's Theorem (1931).}
\newtheorem{thjarg}{Jarn\'{\i}k's General Theorem (1931).}

\newtheorem{thdel}{Theorem A }

\newcommand{\bp}{\mathbf{p}}
\newcommand{\f}{\mathbf{f}}

\newcommand{\cRp}{\cR_{\cC}(\Phi)}
\newcommand{\cJp}{\cJ_{\cC}(\Phi)}
\newcommand{\rao}{R_{\alpha,1}}
\newcommand{\rak}{R_{\alpha,k}}
\newcommand{\waq}{A_Q^*(I)}
\newcommand{\lp}{\lambda_\psi}

\newcommand{\gb}{G_{\! \mbox{\tiny \em B}}^u}
\newcommand{\kt}{u_{t_1}}
\newcommand{\ktk}{u_{t_1+i}}
\renewcommand{\Bbb}[1]{\mathbb{#1}}
\newcommand{\N}{{\Bbb N}}         
\newcommand{\Q}{{\Bbb Q}}         
\newcommand{\R}{{\Bbb R}}         
\newcommand{\Rp}{{\Bbb R}^{+}}    
\newcommand{\Z}{{\Bbb Z}}         
\newcommand{\al}{\alpha}
\newcommand{\ba}{\beta_\alpha}
\newcommand{\ka}{\kappa}
\renewcommand{\r}{\rho}
\newcommand{\p}{\psi}
\newcommand{\bu}{{\bf u}}

\newcommand{\cC}{{\cal C}}

\newcommand{\cE}{{\cal E}}

\newcommand{\cH}{{\cal H}}

\newcommand{\cJ}{{\cal J}}
\newcommand{\cK}{{\cal K}}

\newcommand{\cM}{{\cal M}}

\newcommand{\cQ}{{\cal Q}}
\newcommand{\cR}{{\cal R}}
\newcommand{\cS}{{\cal S}}

\newcommand{\cV}{{\cal V}}

\def\hs{{\cal H}^{s}}

\newcommand{\ve}{\varepsilon}
\newcommand{\tpsi}{\tilde\psi}

\newcommand{\set}[1]{\left\{#1\right\}}

\newcommand{\vv}[1]{{\mathbf{#1}}}
\newcommand{\Veronese}{\cV}
\renewcommand{\le}{\leqslant}
\renewcommand{\ge}{\geqslant}
\newcommand{\ra}{R_{\alpha}}

\newcommand{\mysection}[1]{\refstepcounter{section}%
\section*{{\bf  \thesection.~~#1}}
\addcontentsline{toc}{section}{{\bf  \thesection.~~#1}}}

\newcommand{\mysubsubsection}[1]{\refstepcounter{subsubsection}%
\subsubsection*{\bf \thesubsubsection~~#1}
\addcontentsline{toc}{subsubsection}{{\hspace*{0ex}
\thesubsubsection~~#1}}}

\newcommand{\myappsection}[1]{\refstepcounter{section}%
\section*{{\bf \large  \thesection.~~#1}}
}

\newcommand{\myappsubsection}[1]{\refstepcounter{subsection}%
\subsection*{{\bf \large  \thesubsection~~#1}}
}

\newcommand{\myappsubsubsection}[1]{\refstepcounter{subsubsection}%
\subsubsection*{{\bf \large  \thesubsubsection~~#1}}
}

\begin{document}

\title{Diophantine approximation on planar curves and
the~distribution of rational points}

\author {~~~~~~~~~~~~  Victor Beresnevich \footnote{This work has been partially
supported by INTAS Project 00-429 and by EPSRC grant GR/R90727/01}
\\ ~~~~~~~~~~~~ {\small M{\scriptsize INSK}}
\and ~ \and  Detta Dickinson ~~~~~~~~~~ \\
$\!\!\!\!\!\!\!\!\!\!\!\!\!\!\!\!\!\!\!\!\!\!${\small
M{\scriptsize AYNOOTH}}    \and
\\[0mm] Sanju Velani\footnote{Royal Society University Research
Fellow} \\ {\small Y{\scriptsize ORK}} }

\date{With an Appendix \\ ~ \\  {\Large Sums of two squares near perfect
squares }  \\  
~ \\ R. C. Vaughan \footnote{Research supported by NSA grant
MDA904-03-1-0082.} \\ {\small P{\scriptsize ENN} S{\scriptsize
TATE} }   \\ ~ \\ I{\small N  MEMORY  OF }  P{\small RITISH }
L{\small IMANI} (1983--2003)
}

\maketitle

\vspace{6mm}

\abstract{Let $\cal C$ be a non--degenerate planar curve and for a
real, positive decreasing function $\psi$ let $\cal C(\psi)$
denote the set of simultaneously $\psi$--approximable points lying
on $\cal C$. We show that $\cal C$ is of Khintchine type for
divergence; i.e. if a certain sum diverges then the
one-dimensional Lebesgue measure on $\cal C$ of $\cal C(\psi)$ is
full. We also obtain the Hausdorff measure analogue of the
divergent Khintchine type result.  In the case that $\cal C$ is a
rational quadric  the convergence counterparts of the divergent
results are also obtained. Furthermore, for functions $\psi$ with
lower order in a critical range we determine a general, exact
formula for the Hausdorff dimension of $\cal C(\psi)$.  These
results constitute the first precise and general results in the
theory of simultaneous Diophantine approximation on manifolds.}

\bigskip

\bigskip


\newpage
\parindent=4ex
\parskip=0.6ex

\newpage

\tableofcontents



\newpage

\mysection{Introduction }

In $n$--dimensional Euclidean space  there are two main types of
Diophantine approximation which can be considered, namely
simultaneous and dual. Briefly, the simultaneous case involves
approximating  points $\vv y  =(y_1,\dots,y_n)$ in $\R^n$  by
rational points $\{{\bf p}/q : ({\bf p},q) \in \Z^n \times \Z \} $.
On the other hand, the dual case involves approximating points $\vv
y$  by rational hyperplanes $\{{\bf q.x}= p : (p,{\bf q}) \in \Z
\times \Z^n  \}$ where ${\bf x.y} = x_1y_1 + \dots +x_ny_n$ is the
standard scalar product of two vectors ${\bf x}, {\bf y} \in \R^n$.
In both cases the `rate' of approximation is governed by some given
approximating function. In this paper we consider the general
problem of simultaneous Diophantine approximation on manifolds.
Thus, the points in $\R^n$ of interest are restricted to some
manifold ${\cal M}$  embedded in $\R^n$. Over the past ten years or
so, major advances have been made towards  developing a complete
`metric' theory for the dual form of approximation. However, no such
theory exists for the simultaneous case. To some extent this work is
an attempt to address this in balance.

\subsection{Background and the general problems \label{intro}}

\noindent{\em Simultaneous approximation in $\R^n$. } In order to
set the scene we recall two fundamental results in the theory of
simultaneous Diophantine approximation in $n$--dimensional
Euclidean space. Throughout $\p:\R^+\to\R^+$ will denote a real,
positive decreasing function and will be referred to as an
\emph{approximating function}. Given an approximating function
$\psi$, a point $\vv y=(y_1,\dots,y_n)\in\R^n$ is called {\it
simultaneously $\psi$--approximable} if there are infinitely many
$q\in\N$ such that $$ \max_{1\le i\le n}\|q y_i\|<\psi(q) $$ where
 $\|x\|=\min\{|x-m|:m\in\Z\}$.
In the case $\p$ is $\p_v:h\to h^{-v}$ with $v>0$ the point $\vv y$
is said to be \emph{simultaneously $v$--approximable}. The set of
simultaneously $\p$--approximable points will be denoted by
$\cS_n(\psi)$ and similarly $\cS_n(v)$ will denote the set of
simultaneously $v$--approximable points in $\R^n$. Note that in view
of Dirichlet's theorem ($n$-dimensional simultaneous version),
$\cS_n(v)= \R^n $ for any $v \leq 1/n$.

 The following  fundamental result provides a beautiful and
simple criteria for the  `size' of the set $\cS_n(\psi)$ expressed
in terms of  $n$--dimensional Lebesgue measure $|\  \ |_{\R^n}$.

\begin{thkh} Let   $\p$ be an approximating function. Then
$$|\cS_n(\psi)|_{\R^n} =\left\{\begin{array}{ll} \mbox{\rm
Z{\scriptsize ERO}} & {\rm if} \;\;\; \sum \;  \p(h)^n \;\;
<\infty\\ &
\\ \mbox{\rm F{\scriptsize ULL}} & {\rm if} \;\;\; \sum  \;  \p(h)^n \;\;
 =\infty \; \;
\end{array}\right..$$
\end{thkh}

\noindent Here `full' simply means that the complement of the set
under consideration is of zero measure. Thus the $n$--dimensional
Lebesgue measure of the set of simultaneously $\p$--approximable
points in $\R^n$ satisfies a `zero-full' law. The divergence part of
the above statement  constitutes the main substance of the theorem.
The convergence part is a simple consequence of the Borel-Cantelli
lemma from probability theory. Note that $|\cS_n(v)|_{\R^n} = 0$ for
$v > 1/n$ and so $\R^n$ is extremal -- see below.

 The next fundamental result is a Hausdorff measure version
of the above theorem and shows that the  $s$--dimensional
Hausdorff measure $\hs(\cS_n(\psi))$  of the set $\cS_n(\psi)$
satisfies an elegant  `zero-infinity' law.

\begin{thjar}
 Let $s \in (0,n)$ and  $\p$ be  an
approximating function. Then $$
\hs\left(\cS_n(\psi)\right)=\left\{\begin{array}{ll} 0 & {\rm if}
\;\;\; \sum \; h^{n-s} \,  \p(h)^s \;\;
 <\infty\\ &
\\ \infty & {\rm if} \;\;\; \sum \;
 h^{n-s} \,  \p(h)^s  \;\;  =\infty
\end{array}\right..$$
Furthermore $$ \dim \cS_n(\psi) \ = \ \inf \{ s : \mbox{$\sum $}
\; h^{n-s} \,  \p(h)^s  < \infty \} \; . $$ \label{main}
\end{thjar}

\noindent The dimension part of the statement follows directly from
the definition of Hausdorff dimension -- see  \S\ref{HM}.
In
Jarn\'{\i}k's original statement the additional hypotheses that $r
\p(r)^n \to 0 $ as $ r \to \infty $,   $r \p(r)^n$ is decreasing and
that $r^{1+n-s} \p(r)^s $ is decreasing were assumed. However, these
are not necessary -- see \cite[\S1.1 and \S12.1]{BDV03}. Also,
Jarn\'{\i}k obtained his theorem for general Hausdorff measures $
{\cal H}^h $ where $h$ is a dimension function -- see \S\ref{GHMTHM}
and  \cite[\S1.1 and \S12.1]{BDV03}. However, for the sake of
clarity and ease of discussion we have specialized to
$s$-dimensional Hausdorff measure. Note that the above theorem
implies that for $v > 1/n$ $$ {\cal H}^d \left(\cS_n(v)\right) \ = \
\infty \hspace{10mm} {\rm where } \hspace{10mm} d :=  \dim \cS_n(v)
= \frac{1+n}{v+1} \ .
$$

\noindent The two fundamental theorems stated above  provide a
complete measure theoretic description of  $\cS_n(\psi)$. For a
more detailed discussion and various generalizations of these
theorems  see  \cite{BDV03}.

\medskip
\noindent{\em Simultaneous approximation restricted to manifolds.}
 Let ${\cal M}$ be a manifold of dimension $m$ embedded in
$\R^n$. Given an approximating function $\psi$ consider the set $$
{\cal M} \cap \cS_n(\psi) \  $$ consisting  of points $\vv y $ on
${\cal M}$ which are simultaneously $\psi$--approximable. Two
natural problems now arise.

\medskip

\noindent{\bf Problem 1.} To develop a Khintchine type theory for
$ {\cal M} \cap \cS_n(\psi) $.

\medskip

\noindent{\bf Problem 2.} To develop a Hausdorff measure/dimension
theory for $ {\cal M} \cap \cS_n(\psi) $.

\medskip

\noindent In short, the aim is to establish analogues of the two
fundamental theorems described above and thereby  provide a complete
measure theoretic description of the sets $ {\cal M} \cap
\cS_n(\psi) $. The fact that the points $\vv y$ of interest are of
dependent variables, reflecting the fact that $\vv y \in {\cal M}$
introduces major difficulties in attempting to describe the measure
theoretic structure of $ {\cal M} \cap \cS_n(\psi) $. This is true
even in the specific case that ${\cal M}$ is a planar curve. More to
the point, even for seemingly simple curves such as the unit circle
or the parabola the problem is fraught with difficulties.

\medskip

\noindent{\em Non-degenerate manifolds. } In order to make any
reasonable progress with the above problems it is not unreasonable
to assume that the manifolds $ {\cal M}$ under consideration are
{\bf non-degenerate} \cite{KM98}. Essentially, these are smooth
sub-manifolds of $\R^n$ which are sufficiently curved so as to
deviate from any hyperplane. Formally, a  manifold $\cM$ of
dimension $m$ embedded in $\R^n$ is said to be non-degenerate if it
arises from a non--degenerate map $\f:U\to \R^n$ where $U$ is an
open subset of $\R^m$ and $\cM:=\f(U)$. The map $\f:U\to
\R^n:\bu\mapsto \f(\bu)=(f_1(\bu),\dots,f_n(\bu))$ is said to be
\emph{non--degenerate at} $\bu\in U$ if there exists some $l\in\N$
such that $\f$  is $l$ times continuously differentiable on some
sufficiently small ball centred at $\bu$ and the partial derivatives
of $\f$ at $\bu$ of orders up to $l$ span $\R^n$. The map $\f$ is
\emph{non--degenerate} if it is non--degenerate at almost every (in
terms of $m$--dimensional Lebesgue measure) point in $U$; in turn
the manifold $\cM=\f(U)$ is also said to be non--degenerate.  Any
real, connected analytic manifold not contained in any hyperplane of
$\R^n$  is non--degenerate.

 Note that in the case the manifold $\cM$ is a planar curve
${\cal C}$,  a point on ${\cal C}$ is non-degenerate if  the
curvature at that point is non-zero. Thus, ${\cal C}$ is a
non-degenerate planar curve if   the set of points on ${\cal C}$
at which the  curvature vanishes is a set of one--dimensional
Lebesgue measure zero. Moreover, it is not difficult to show that
the set  of points on a planar curve at which the curvature
vanishes but the curve  is non-degenerate is at most countable. In
view of this, the curvature completely describes the
non-degeneracy of planar curves. Clearly, a straight line is
degenerate everywhere.

\subsection{The Khintchine type theory \label{Ktheory} }

The aim is to obtain an analogue of Khintchine's theorem for the set
$ {\cal M} \cap \cS_n(\psi) $ of simultaneously $\psi$--approximable
points lying on ${\cal M}$. First of all notice that if the
dimension $m$ of the manifold ${\cal M}$ is strictly less than $n$
then $ |{\cal M} \cap \cS_n(\psi)|_{\R^n} = 0  $ irrespective of the
approximating function $\psi$.  Thus, when referring to the Lebesgue
measure of the set $ {\cal M} \cap \cS_n(\psi) $ it is always  with
reference to the induced Lebesgue measure on ${\cal M}$. More
generally, given a subset $S$ of ${\cal M}$ we shall write
$|S|_{\cal M} $ for the  measure of $S$ with respect to the induced
Lebesgue measure on ${\cal M}$. Notice that for $v \leq 1/n $,  we
have that $|{\cal M} \cap \cS_n(v)|_{{\cal M}} = |{\cal M}|_{{\cal
M}} := \mbox{F{\scriptsize ULL}} $  as it should be since $\cS_n(v)=
\R^n $.

 To develop the Khintchine theory it is natural to
consider the convergence and divergence cases separately and the
following terminology is most useful.

\begin{definition}\sl
Let $\cM\subset\R^n$ be a manifold. Then we say that
\begin{enumerate}
\item $\cM$ is  of {\it Khintchine type for convergence}\/ if
$|\cM\cap\cS_n(\psi)|_\cM= \mbox{{\rm Z{\scriptsize ERO}}}$ for
any approximating function $\psi$ with
$\sum_{h=1}^\infty\psi(h)^n<\infty$.
\item $\cM$ is  of {\it Khintchine type for divergence}\/ if
$|\cM\cap\cS_n(\psi)|_\cM= \mbox{{\rm F{\scriptsize ULL}}}$ for
any approximating function $\psi$ with
$\sum_{h=1}^\infty\psi(h)^n=\infty$.
\end{enumerate}
\end{definition}

\noindent The set of manifolds which are of  Khintchine type for
convergence will be denoted by $\cK_{<\infty}$. Similarly, the set
of manifolds which are of Khintchine type for divergence will be
denoted by $\cK_{=\infty}$. Also, we define
$\cK:=\cK_{<\infty}\cap\cK_{=\infty}$. By definition, if $\cM \in
\cK$ then an analogue of Khintchine's theorem exists for
$\cM\cap\cS_n(\psi)$ and $\cM$ is simply said to be of Khintchine
type. Thus Problem 1 mentioned above, is equivalent to describing
the set of Khintchine type manifolds. Ideally, one would like to
prove that any non-degenerate manifold is of Khintchine type.
Similar terminology exists for the dual form of approximation in
which `Khintchine type' is replaced by `Groshev type'; for further
details see \cite[pp.\,29--30]{BD99}.

 A weaker notion than `Khintchine type for convergence'
is that of extremality. A manifold $\cM$ is said to be {\it
extremal}\/ if $|\cM\cap\cS_n(v)|_\cM=0$ for any $v>1/n$. The set of
extremal manifolds of $\R^n$ will be denoted by $\cE$ and it is
readily verified that  $\cK_{<\infty}\subset\cE$. In 1932, Mahler
made the conjecture that for any $n \in \N$ the Veronese curve
$\Veronese_n=\{(x,x^2,\dots,x^n):x\in\R\}$ is extremal. The
conjecture  was eventually  settled   in 1964 by Sprindzuk
\cite{Spr64}  -- the special cases $n=2$ and $3$ had been done
earlier. Essentially, it is this  conjecture and its investigations
which gave rise  to the now flourishing  area of `Diophantine
approximation on manifolds' within metric number
theory. 
Up to 1998, manifolds satisfying a variety of analytic, arithmetic
and geometric constraints had been shown to be extremal. For
example, Schmidt in 1964  proved that  any $C^3$ planar curve with
non-zero curvature almost everywhere is extremal. However,
Sprindzuk in the 1980's, had conjectured that any analytic
manifold  satisfying a necessary non--degeneracy condition is
extremal.  In 1998, Kleinbock and Margulis \cite{KM98} showed that
any non-degenerate manifold is extremal and thereby settled the
conjecture of Sprindzuk.

Regarding the `Khintchine theory' very little is known. The
situation for the dual form of approximation  is very different. For
the dual case, it has recently been shown that any non--degenerate
manifold is of Groshev type -- the analogue of Khintchine type in
the dual case (see  \cite{BBKM02}, \cite{BKM01} and
\cite[\S12.7]{BDV03}). For the simultaneous case, the current state
of the Khintchine theory is somewhat ad-hoc. Either a specific
manifold or a special class of manifolds satisfying various
constraints is studied.
For example it has been shown that (i) manifolds which are a
topological product of at least four non--degenerate planar curves
are in $\cK$ \cite{Ber73};
(ii) the parabola $\Veronese_2$ is in $\cK_{<\infty}$
\cite{Ber79}; (iii) the so called 2--convex manifolds of dimension
$m \geq 2$
are in $\cK_{<\infty}$ \cite{DRV91}  
and (iv) straight lines through the origin satisfying a natural
Diophantine condition are in $\cK_{<\infty}$
\cite{Kov00}. 
Thus, even in the simplest geometric and arithmetic situation in
which the manifold is a genuine curve in $\R^2$ the only known
result to date is that of  the parabola $\Veronese_2$. To our
knowledge, no curve has ever been shown to be in $\cK_{=\infty}$.

\medskip

In this paper we address the fundamental problems of \S\ref{intro}
in the case that the manifold $\cM$ is a  planar curve (the
specific case that $\cM$ is a non-degenerate, rational quadric
will be shown in full). Regarding Problem 1, our main result is
the following. As usual, $C^{(n)}(I)$ will denote the set of
$n$--times continuously differentiable functions defined on some
interval $I$ of $\R$.

\begin{theorem}\label{thm1}
Let $\psi$ be an approximating function with
$\sum_{h=1}^\infty\psi(h)^2=\infty$. Let $f\in C^{(3)}(I_0)$,
where $I_0$ is an interval, and $f''(x)\not=0$ for almost all
$x\in I_0$. Then for almost all $x\in I_0$ the point $(x,f(x))$ is
simultaneously $\psi$--approximable.
\end{theorem}

\begin{corollary}\label{corol1}
Any $C^{(3)}$ non--degenerate planar curve is of Khintchine type
for divergence.
\end{corollary}

To complete the `Khintchine theory' for $C^{(3)}$ non--degenerate
planar curves we need to show that any such curve is of Khintchine
type for convergence. We are currently  able to prove this in the
special case that the planar curve is a non-degenerate, rational
quadric. 
However, the truth of Conjecture 1 in \S\ref{rp} regarding the
distribution of rational points `near' planar curves would yield
the complete convergence theory.

\mysubsubsection{The Khintchine theory for  rational quadrics
\label{krq}} As above, let $\Veronese_2:= \{(x_1,x_2) \in \R^2 : x_2
= x_1^2 \}$ denote the standard parabola and let $\cC_1 :=
\{(x_1,x_2) \in \R^2 : x_1^2 + x_2^2 =1  \} $  and $\cC_1^* :=
\{(x_1,x_2) \in \R^2 : x_1^2 - x_2^2 =1  \} $ denote the unit circle
and standard hyperbola respectively. Next, let $\cQ$ denote a
non-degenerate, rational quadric in the plane. By this we mean that
$\cQ$ is the image of either the circle $\cC_1$, the hyperbola
$\cC_1^*$ or the parabola $\Veronese_2$ under a rational affine
transformation of the plane. Furthermore, for an approximating
function $\psi$ let
$$\cQ(\psi) \ := \ \cQ  \cap \cS_2(\psi)$$

\noindent In view of Corollary \ref{corol1} we have that $\cQ $ is
in $  \cK_{=\infty}$. The following result shows  that any
non-degenerate, rational quadric  is in fact in $\cK$ and provides
a complete criteria for the size of $\cQ(\psi)$  expressed in
terms of Lebesgue measure. Clearly, it contains the only
previously known result that the parabola  is in $\cK_{<\infty}$.

\begin{theorem}\label{corol2C}
Let $\psi$ be an approximating function. Then

 $$\big|
\cQ(\psi)\big|_{\cQ}= \left\{\begin{array}{ll} \mbox{\rm
Z{\scriptsize ERO}} & {\rm if} \;\;\; \sum \; \psi(h)^2 \; < \;
\infty
\\ &
\\ \mbox{\rm F{\scriptsize ULL}} & {\rm if} \;\;\; \sum  \;  \psi(h)^2 \;
 = \; \infty \; \;
\end{array}\right..$$
\end{theorem}

\medskip

\subsection{The Hausdorff measure/dimension theory \label{hmhdtheory}}

The aim is to obtain an analogue of Jarn\'{\i}k's theorem for the
set $ {\cal M} \cap \cS_n(\psi) $ of simultaneously
$\psi$--approximable points lying on ${\cal M}$. In the dual case,
the analogue of the divergent part of Jarn\'{\i}k's theorem has
recently been established  for any non-degenerate manifold
\cite[\S12.7]{BDV03}. Prior to this,  a general lower bound for
the Hausdorff dimension of the dual  set of $v$--approximable
points lying on any extremal manifold had been obtained
\cite{DD00}. Also in the dual case,  exact formulae for the
dimension of the dual $v$-approximating  sets are known for  the
case of  the Veronese curve \cite{BS70, Ber83} and for any  planar
curve with curvature non-zero except for a set of dimension zero
\cite{Bak78}.

As with the Khintchine theory, very little is currently known
regarding  the Hausdorff measure/dimension theory for the
simultaneous case.  Contrary to the dual case,
$\dim\cM\cap\cS_n(v)$  behaves in a rather complicated way and
appears to depend on the arithmetic properties of $\cM$. For
example, let $\cC_R=\{x^2+y^2=R^2\}$ be the circle of radius $R$
centered at the origin. It is easy to verify that $\cC_{\sqrt 3}$
contains no rational points $(s/q,t/q)$. On the other hand, any
Pythagorean triple $(s,t,q)$ gives rise to a rational point on the
unit circle $\cC_1$ and so there are  plenty of rational points on
$\cC_1$. For $v>1$, these facts regarding the distribution of
rational points on the circle under consideration  lead to
$\dim\cC_{\sqrt3}\cap\cS_2(v)=0$ whereas $\dim\cC_1\cap\cS_2(v)=
1/(1+v)$ \cite{BDV03,DD01}. The point is that for $v>1$, the
rational points of interest must lie on the associated circle.
Further evidence for the complicated behavior of the dimension can
be found in \cite{Ryn}. Recently, $\dim\cM\cap\cS_n(v)$ has been
calculated for  large values of $v$ when the  manifold $\cM$ is
parameterized by polynomials with integer coefficients \cite{DL} and
for $v >1$ when  the manifold is a non-degenerate, rational quadric
in $\R^n$ \cite{drutu}. Also, as a consequence of Wiles' theorem
\cite{wiles}, $\dim\cM\cap\cS_2(v)=0$ for the curve $x^k+y^k=1$ with
$k>2$ and $v>k-1$ \cite[p.~94]{BD99}.

The above examples illustrate that in the simultaneous case there is
no hope of establishing a single, general formula for
$\dim\cM\cap\cS_n(v)$. Recall, that for $v=1/n$ we have that
$\dim\cM\cap\cS_n(v) = \dim \cM := m $ for any manifold embedded in
$\R^n$  since $\cS_n(v) = \R^n$ by Dirichlet's theorem. Now notice
that in the various examples considered above the varying behaviour
of $ \dim\cM\cap\cS_n(v)$ is exhibited  for values of $v$ bounded
away from the Dirichlet  exponent $1/n$. Nevertheless, it is
believed that when $v$ lies in a critical range near the Dirichlet
exponent $1/n$ then, for a wide class of manifolds (including
non--degenerate manifolds), the behaviour of $\dim\cM\cap\cS_n(v)$
can be captured by a single, general formula.  That is to say, that
$ \dim\cM\cap\cS_n(v)$ is independent  of the arithmetic properties
of $\cM$ for $v$ close to $1/n$. We shall  prove that this is indeed
the case for planar curves. Note that for  planar curves the
Dirichlet exponent is $1/2$ and that the above `circles example'
shows that any critical range for $v$ is a subset of $[1/2,1]$. In
general, the critical range is governed by the dimension of the
ambient space and the dimension of the manifold.

Before stating our results we introduce the notion of lower order.
Given an approximating function $\p$, the {\em lower order} $\lp$
of $1/\p$ is defined by $$ \lp \ := \ \liminf_{h\to\infty}\;
\frac{- \log \p(h)}{\log h} \ , $$ and indicates the growth of the
function $1/\psi$ `near' infinity. Note that $\lp $  is
non-negative since $\psi$ is a decreasing function. Regarding
Problem 2, our main results are as follows.

\begin{theorem}\label{thm2HM}
Let $f\in C^{(3)}(I_0)$, where $I_0$ is an interval and $\cC_f :
=\set{(x,f(x)):x\in I_0}$.  Assume that there exists at least one
point on the curve $\cC_f $ which is non-degenerate. Let $s \in
(1/2, 1)$ and $\p$ be an approximating function.
Then $$ {\cal H}^s
(\cC_f\cap\cS_2(\p)) \ = \ \infty \ \hspace{5mm} \mbox{if}
\hspace{5mm} \sum_{h =1}^{ \infty} \ h^{1-s} \, \p(h)^{s+1}  \ = \
\infty \ . $$
\end{theorem}


\begin{theorem}\label{thm2a}
Let $f\in C^{(3)}(I_0)$, where $I_0$ is an interval and $\cC_f :
=\set{(x,f(x)):x\in I_0}$. Let $\p$ be an approximating function
with $\lp \in [1/2, 1)$. Assume that
\begin{equation} \dim\set{x\in I_0:
f''(x)=0}\le\frac{2-\lp}{1+\lp} \ . \label{dimcond} \end{equation}
Then $$ \dim\cC_f\cap\cS_2(\p) \ = \ d \ := \  \frac{2-\lp}{1+\lp}
\ . $$

\noindent Furthermore, suppose that  $\lp \in (1/2, 1)$. Then $$
{\cal H}^d (\cC_f\cap\cS_2(\p)) \ = \ \infty \ \hspace{5mm}
\mbox{if} \hspace{5mm}\limsup_{h \to \infty}  h^{2-s} \p(h)^{s+1}
> 0 \ . $$
\end{theorem}

By considering the function $\p : h \to h^{-v} $, an immediate
consequence of the theorems is the following corollary.

\begin{corollary}\label{thm2}
Let $f\in C^{(3)}(I_0)$, where $I_o$ is an interval and $\cC_f :
=\set{(x,f(x)):x\in I_0}$. Let $v \in [1/2,1) $ and assume that $
\dim\set{x\in I_0: f''(x)=0}\le (2-v)/(1+v)  $.  Then $$
\dim\cC_f\cap\cS_2(v) \; = \; d \; := \; \frac{2-v}{1+v}  \ . $$
Moreover, if $v \in (1/2,1) $  then  
${\cal H}^d (\cC_f\cap\cS_2(v))  = \infty $.
\end{corollary}

\vspace{0mm}

\noindent{{\bf Remark.}} Regarding Theorem \ref{thm2a}, the
hypothesis (\ref{dimcond}) on the set $\{x\in I_0: f''(x)=0\}$ is
stronger than simply assuming that the curve $\cC_f $ is
non-degenerate. It requires the curve to be non--degenerate
everywhere except on a set of Hausdorff dimension no larger than
$(2-\lp)/(1+\lp)$ -- rather than just measure zero. Note that the
hypothesis can be made independent of the lower order $\lp$ (or
indeed of $v$ in the case of the corollary) by assuming that $
\dim\{x\in I_0: f''(x)=0\} \leq 1/2 $. The proof of Theorem
\ref{thm2a} follows on establishing the upper and lower bounds for
$\dim\cC_f\cap\cS_2(\p)$ separately. Regarding the lower bound
statement, all that is required is that there exists at least one
point on the curve $\cC_f $ which is non-degenerate. This is not
at all surprising since  the lower bound statement can be viewed
as a simple consequence of
Theorem \ref{thm2HM}. 
The hypothesis (\ref{dimcond}) is required to obtain the upper bound
dimension statement. Even for non-degenerate curves, without such a
hypothesis the statement of Theorem \ref{thm2a} is clearly false as
the following example shows.
\medskip

\noindent{{\bf Example:} {\em The Cantor curve. \ }} Let $K$
denote the standard middle third Cantor set obtained by removing
the middle third of the unit interval $[0,1]$ and then inductively
repeating the process on each of the remaining intervals. For our
purpose, a convenient expression for $K$ is the following:  $$
\textstyle{ \bigcap_{i=1}^\infty (
[0,1]\setminus\bigcup_{\,j=1}^{\,2^{i-1}} I_{i,j})  =
[0,1]\setminus\bigcup_{i=1}^\infty\bigcup_{j=1}^{2^{i-1}} I_{i,j}}
\ \ ,  $$ where $I_{i,j}$ is  the $j^{\rm th}$  interval of the
$2^{i-1}$ open intervals of length $3^{-i}$ removed at the $i^{\rm
th}$-level of the Cantor construction. Note that the intervals
$I_{i,j}$ are pair wise disjoint. Give a pair $(i,j)$, define the
function $$ f_{i,j} \, : \, x \to  f_{i,j}(x) \ := \ \left\{
\begin{array}{cl}
e^{-i\textstyle-\frac{1}{(x-a)(b-x)}} &  \ \ {\rm if \ }  \ \ x\in
I_{i,j}
\\[1ex] 0  & \ \  {\rm if \ }  \ \ x\in [0,1]\setminus I_{i,j}
\end{array}
\right.  ,  $$ where $a$ and $b$ are the end points of the
interval $I_{i,j}$. Now set $$ f : x \to f(x) \ := \
\sum_{i=1}^\infty \ \sum_{j=1}^{2^{i-1}} f_{i,j}(x) \ \ . $$ Note
that the function $f$ is obviously $C^{(\infty)}$ as the sum
converges uniformly. Also,  for $x\in K$ and  $m\in\N$ we have
that $f_{i,j}^{(m)}(x)=0$ and so $$ f^{(m)}(x) \ = \
\sum_{i=1}^\infty \ \sum_{j=1}^{2^{i-1}} f_{i,j}^{(m)}(x) \ = \ 0
\ \ . $$  On the other hand, for $x\in[0,1]\smallsetminus K$ we
have that $f^{(m)}(x)>0$. Thus the curve
$C_{K}=\{(x,f(x)):x\in(0,1)\}$ is exactly degenerate on $K$ and
non-degenerate elsewhere. Note that $C_{K}$ is a non-degenerate
curve since $K$ is of Lebesgue measure zero.  The upshot of this
is that for any $x\in K$ the point $(x,f(x))$ is $1$-approximable;
i.e. there exists infinitely many $q\in\N$ such that  $$ \|qx\|<
q^{-1} \qquad \ {\rm and } \qquad \|qf(x)\|<q^{-1} \ \ . $$ The
second inequality is trivial as $f(x)=0$ and  the first inequality
is a consequence of Dirichlet's theorem. Thus, $$\dim
\cC_K\cap\cS_2(v) \ \geq \ \dim K  \, = \, \log 2/\log 3  $$
irrespective of $v\in(1/2,1)$. Obviously, by choosing  Cantor sets
$K$  with dimension close to one, we can ensure that $\dim
\cC_K\cap\cS_2(v) $ is close to one irrespective of $v\in(1/2,1)$.

\vspace{5mm}

For simultaneous Diophantine approximation on planar curves,
Theorem \ref{thm2HM} is the precise analogue of the divergent part
of Jarn\'{\i}k's theorem and Theorem \ref{thm2a} establishes a
complete Hausdorff dimension theory.

Note that the measure part of Theorem \ref{thm2a} is substantially
weaker than Theorem \ref{thm2HM} -- the general measure statement.
For example,  with $v \in (1/2,1)$ and $\alpha =1/(d+1)$ consider
the approximating function $\p$ given by
$$ \p: h \to h^{-v} (\log h)^{-\alpha}  \ . $$ Then $\lp = v $ and
assuming that (\ref{dimcond}) is satisfied, the dimension part of
Theorem \ref{thm2a} implies that $$ \dim\cC_f\cap\cS_2(\p) \ = \ d
\ := \ \frac{2-v}{1+v} \ . $$ However,   $$ \limsup_{h \to \infty}
h^{2-d} \p(h)^{d+1} \ = \ \lim_{h \to \infty} (\log h)^{ -1} \ = \
0 \  $$ and so the  measure part of Theorem \ref{thm2a} is not
applicable. Nevertheless,
$$\sum h^{1-d} \p(h)^{d+1} = \sum (h\log h)^{-1} = \infty $$ and
Theorem \ref{thm2HM}  implies that $ {\cal H}^d
(\cC_f\cap\cS_2(\p)) \ = \ \infty \ . $

\medskip

 Theorem \ref{thm2HM}   falls short of establishing a
complete Hausdorff measure theory for simultaneous Diophantine
approximation on planar curves. In its simplest form, it should be
possible to summarize the Hausdorff measure theory by a clear cut
statement of the following type.


\smallskip

\noindent{{\bf Conjecture  H }} {\em Let $s \in (1/2,1)$ and  $\p$
be an approximating function.
Let $f\in C^{(3)}(I_0)$, where $I_o$ is an interval and $\cC_f :
=\set{(x,f(x)):x\in I_0}$. Assume that $ \dim\{x\in I_0:
f''(x)=0\} \leq 1/2 $.  Then $$
\hs\left(\cC_f\cap\cS_2(\p)\right)=\left\{\begin{array}{ll} 0 &
{\rm if} \;\;\; \sum \; h^{1-s} \,  \p(h)^{s+1} \;\;
 <\infty\\ &
\\ \infty & {\rm if} \;\;\; \sum \;
 h^{1-s} \,  \p(h)^{s+1}  \;\;  =\infty
\end{array}\right..$$  }

\smallskip

The divergent part of the above statement is Theorem \ref{thm2HM}.
As with the `Khintchine theory', the above convergent part  would
follow on proving Conjecture 1 of \S\ref{rp}. However, for
rational quadrics   we are able to prove the convergent result
independently of any conjecture.

\begin{theorem}\label{circ2}
Let $s \in (1/2,1)$ and  $\p$ be an approximating function. Then
for any non-degenerate, rational quadric $\cQ$ we have that
$$\hs\left(\cQ\cap\cS_2(\p)\right) = 0
 \hspace{8mm} {\rm if} \;\;\; \sum \; h^{1-s} \,  \p(h)^{s+1} \;\;
 <\infty \ .  $$
\end{theorem}

\subsection{Rational points close to a curve \label{rp}}

First some useful notation. For any point $\vv r\in\Q^n$ there
exists the smallest $q\in\N$ such that $q\vv r\in\Z^n$. Thus, every
point $\vv r\in\Q^n$ has a unique representation in the form
$$\frac{\bp}{q}=\frac{(p_1,\dots,p_n)}{q}=
\left(\frac{p_1}{q},\dots,\frac{p_n}{q}\right)$$ with
$(p_1,\dots,p_n)\in\Z^n$. Henceforth, we will only consider points
of $\Q^n$ in this form.

Understanding the distribution of rational points close to a
reasonably defined  curve is absolutely crucial towards  making
any progress with the main problems considered in this paper. More
precisely, the behaviour of the following counting function will
play a central role.

\noindent{\em  The function $ N_f(Q,\p,I)$. \ }Let $I_0$ denote a
finite, open interval of $\R$ and let $f$ be a function in
$C^{(3)}(I_0)$  such that
\begin{equation}
\label{infsup} 0\ < \ c_1 \ := \ \inf_{x \in I_0} |f''(x)|  \ \leq
\ c_2 \ := \ \sup_{x \in I_0} |f''(x)| \ < \ \infty  \ .
\end{equation} Given an interval $I\subseteq I_0$, an
approximating function $\psi$ and $Q \in \R^+$ consider the
counting function $N_f(Q,\p,I)$ given by $$ N_f(Q,\p,I) \ := \
\#\{\bp/q\in\Q^2 \, : \, q\le Q,\, p_1/q\in I,\,
|f(p_1/q)-p_2/q|<\p(Q)/Q\}. $$

\noindent In short, the function $N_f(Q,\p,I)$ counts `locally'
the number of rational points  with bounded  denominator lying
within a specified  neighbourhood of the curve parameterized by
$f$. In \cite{Hux96}, Huxley obtains a reasonably sharp upper
bound for $N_f(Q,\p,I)$. We will obtain an exact lower bound and
also prove that the rational points under consideration are
`evenly' distributed. The proofs of the Khintchine type and
Hausdorff measure/dimension theorems  stated in this paper rely
heavily on this information.  In particular, the exact upper bound
in Theorem~\ref{thm2a} is easily established  in view of Huxley's
result \cite[Theorem 4.2.4]{Hux96} which we  state in a simplified
form.

\noindent {\bf Huxley's estimate}: {\it Let $\psi $ be an
approximating function such that  $t\p(t)\to \infty$ as $t\to
\infty$.  For $\ve>0$ and $Q$ sufficiently large
\begin{equation}\label{e:001a}
N_f(Q,\p,I_0) \ \le \  \p(Q) \, Q^{2+\ve}.
\end{equation}
}

  The complementary lower bound is the substance of our next
result.

\begin{theorem}\label{thm3}
Let $\psi$ be an approximating  function satisfying
\begin{equation}\label{e:002}
\lim_{t\to+\infty}\psi(t)=\lim_{t\to+\infty}\
\frac{1}{t\psi(t)}=0.
\end{equation}
There exists a constant $c>0$,
 depending on $I$, such that for $Q$ sufficiently large
$$ N_f(Q,\p,I) \ \ge \  c \,  Q^2 \, \psi(Q) \, |I|  \ \ . $$
\end{theorem}

We suspect that the lower bound given by Theorem~\ref{thm3} is
best possible up to a constant multiple. It is plausible that for
compact curves, the constant $c$ is independent of $I$.

Regarding Huxley's estimate, the presence of the `$\ve$' factor
prevents us from proving the desired `convergent' measure
theoretic results.  We suspect that a result of the following type
is in fact true -- proving it is another matter.

\begin{conjecture}\label{conj1}\sl
Let $\psi $ be an approximating function such that  $t\p(t)\to
\infty$ as $t\to \infty$.  There exists a  constant $\hat c>0$
such that for $Q$ sufficiently large $$ N_f(Q,\p,I_0) \ \le \ \hat
c \,  Q^2 \, \psi(Q)  \ .  $$
\end{conjecture}

Conjecture~\ref{conj1} has immediate consequences  for the main
problems considered in this paper. In particular, it would imply
the following.

\begin{conjecture}\label{conj2}\sl
Any $C^{(3)}$ non--degenerate planar curve is of Khintchine type
for convergence.
\end{conjecture}

Conjecture~\ref{conj2} would naturally complement
Theorem~\ref{thm1} of this paper. The implication
Conjecture~\ref{conj1} $\Longrightarrow$ Conjecture~\ref{conj2} is
reasonably  straightforward -- simply modify  the argument set out
in the proof of Theorem \ref{corol2C}. Also, it is not difficult
to verify that Conjecture~\ref{conj1} implies the `convergent'
part of Conjecture H -- simply modify  the argument set out in the
proof of Theorem \ref{circ2}. An intriguing  problem  is to
determine whether or not the two conjectures stated above are in
fact equivalent.

\mysection{Proof of the  rational quadric statements}

\subsection{Proof of Theorem \ref{corol2C} } The divergence part of the theorem
is a trivial consequence of Corollary \ref{corol1} to Theorem
\ref{thm1}. To establish the convergence part  we proceed as
follows.

Let $\psi$ be an approximating function such that $\sum \psi(h)^2
<\infty $. The claim is that $\big| \cQ(\psi)\big|_{\cQ} = 0 $. We
begin by introducing  an auxiliary function $\Psi $ given by $$
\Psi(h) \; := \; \max \left\{ \psi(h), \ h^{-\frac12} \; (\log
h)^{-1} \right\} \ \ . $$
 Clearly, $\Psi$ is an approximating
function and furthermore $$ \sum \Psi(h)^2 <\infty  \hspace{1cm}
{\rm and \ } \hspace{1cm} \Psi(h) \ \geq \ \psi(h)   \  \ . $$ Thus
$\cQ(\psi) \subset \cQ(\Psi)$  and the claim will follow on showing
that  $\big| \cQ(\Psi)\big|_{\cQ} = 0 $. It is easily verified that
such a `zero' statement is invariant under rational affine
transformations of the plane. In view of this, it suffices to
consider the  curves $\cC_1$, $\cC_1^*$ and $\Veronese_2$ -- see
\S\ref{krq}.

In the following,  $C(q;s,t)$ will denote the square with centre
at the rational point  $(s/q,t/q)$ and of side length $2 \Psi(q)/q
$.

\noindent{\bf Case (a): $\cQ = \cC_1$. } For $m \in \N$, let $$
W_m(\Psi;\cC_1) \ := \ \bigcup_{2^{m} < q \leq 2^{m+1} }
\bigcup_{(s,t)\in \Z^2 } \cC_1 \cap C(q;s,t) \ . $$ Then $
\cC_1(\Psi) = \limsup_{m \to \infty} W_m(\Psi;\cC_1) $ and in view
of the Borel-Cantelli lemma $\big| \cC_1(\Psi)\big|_{\cC_1} = 0 $
if $ \sum \big| W_m(\Psi;\cC_1) \big|_{\cC_1} < \infty $. Next,
note that if $\cC_1 \cap C(q;s,t) \neq \emptyset$ then $(q -
2\sqrt2 \Psi(q))^2 \, \leq \, s^2 + t^2 \, \leq \, (q + 2\sqrt2
\Psi(q))^2 $ and $\big|\cC_1 \cap C(q;s,t)\big|_{\cC_1}   \ll
\Psi(q)/q $. It follows that
\begin{eqnarray}
\big| W_m(\Psi;\cC_1) \big|_{\cC_1}  & \ll & \sum_{2^{m} < q \leq
2^{m+1} } \sum_{\substack{(s,t)\in \Z^2: \\ (q - 2\sqrt2
\Psi(q))^2 \, \leq s^2 + t^2 \, \leq \, (q + 2\sqrt2 \Psi(q))^2 }}
\!\!\!\!\!\!\!\!\!\!\!\!\!\!\! \big|\cC_1 \cap
C(q;s,t)\big|_{\cC_1} \nonumber \\ & & \nonumber \\ & \ll &
\frac{\Psi(2^{m})}{2^{m}} \sum_{2^{m} < q \leq 2^{m+1} }
\sum_{\substack{n: \\ |q - \sqrt{n}| < 4\Psi(q) }} \!\!\!\!\! r(n)
 \ \ \ ,
 \label{upbdwm}
\end{eqnarray}
where $r(n)$ denotes the number of representations of $n$ as the
sum of two squares.

With reference to Theorem A of Appendix II, with $\psi := 4 \Psi$,
$Q:=2^m$ and $N := [Q/\Psi(Q)] $ it is easily verified that the
error term associated with $\sum_{Q<q\le 2Q} {\sum_n}\!\!\!' \,
r(n) $ is $$ \ll Q^{\frac{15}{8} } (\log Q)^{65} \Psi(Q) \ \ .
$$ Here we use the trivial fact that $\Psi(Q^*):= \Psi(Q+1) \leq
\Psi(Q) $ since $\Psi $ is decreasing. On the other hand, for the
main term we have that $$  Q^2 \Psi(2Q) \ \ll \ \sum_{Q<q\le 2Q}
q\Psi(q)\ \ll \  Q^2 \Psi(Q)  \ \ . $$ Thus, Theorem A implies
that \begin{equation}
 \sum_{2^{m} < q \leq
2^{m+1} } \sum_{\substack{n: \\ |q - \sqrt{n}| < 4\Psi(q) }}
\!\!\!\!\! r(n) \ \ll \ 2^{2m} \ \Psi(2^m)  \  \ .
\label{thebigone} \end{equation}
 This estimate together with
(\ref{upbdwm}) implies that $ \big| W_m(\Psi;\cC_1) \big|_{\cC_1}
 \ll
2^m \ \Psi(2^{m})^2   $. In turn, we obtain that $$ \sum_{m \in
\N} \big| W_m(\Psi;\cC_1) \big|_{\cC_1} \ \ll \ \sum_{m \in \N}
2^m \ \Psi(2^{m})^2 \ \asymp \ \sum_{h \in \N} \Psi(h)^2 \ < \
\infty \ .  $$ This completes the proof of the theorem in  the
case that $\cQ$ is the image of the unit circle  $\cC_1$ under a
rational affine transformation of the plane. The other two cases
are similar. The key is to bring (\ref{thebigone}) into play.

\medskip

\noindent{\bf Case (b): $\cQ = \cC_1^*$. } For $k \in \N $, let
$\cC_{1;k}^* := \{(x_1,x_2) \in \R^2 : x_1^2 - x_2^2 =1 {\rm \
with \ } |x_1| \leq 2^k  \} $. Thus, $\cC_{1;k}^* $ is the
hyperbola  $\cC_{1}^*$ with the first co-ordinate bounded above by
$2^k$. For $m \in \N $, let $$ W_{m}(\Psi;\cC_{1;k}^* ) \ := \
\bigcup_{2^{m} < q \leq 2^{m+1} } \bigcup_{(s,t)\in \Z^2 }
\cC_{1;k}^* \cap C(q;s,t) \  $$ and let  $ \cC_{1;k}^*(\Psi) :=
\limsup_{m \to \infty} W_{m}(\Psi;\cC_{1;k}^*) $. Clearly,
$\cC_{1}^*(\Psi) = \bigcup_{k=1}^{\infty} \cC_{1;k}^*(\Psi) $ and
so $\big| \cC_{1}^*(\Psi)\big|_{\cC_{1}^*} = 0 $ if $\big|
\cC_{1;k}^*(\Psi)\big|_{\cC_{1}^*} = 0 $ for each $k \in \N$. The
latter follows on showing that $  \sum \big|
W_{m}(\Psi;\cC_{1;k}^*) \big|_{\cC_{1}^*} < \infty $.

\noindent It is easily verified that if $\cC_{1;k}^* \cap C(q;s,t)
\neq \emptyset$ then $1/2 < |s|/q < a:=2^{k+1} $, $|t| < |s| $ and
\begin{eqnarray*}
|q^2 + t^2 - s^2| \ < \  8 \, |s| \, \Psi(q) + 8 \,  \Psi(q)^2 \ <
\ 8 \, |s|\, \Psi( |s|/a) + 8 \, \Psi(|s|/a)^2 \ .
\end{eqnarray*}
Here we have used that fact that the function  $\Psi$ is
decreasing. It follows via (\ref{thebigone}), that for $m$
sufficiently large
\begin{eqnarray*}
\big| W_{m}(\Psi;\cC_{1;k}^*) \big|_{\cC_1}  & \ll &
\frac{\Psi(2^{m})}{2^{m}}  \sum_{2^{m} < q \leq 2^{m+1} }
\sum_{\substack{(s,t)\in \Z^2 \; :  \ q/2 < s < a q  \\ (s - 8
\Psi(s/a))^2 \, \leq q^2 + t^2 \, \leq \, (s + 8 \Psi(s/a))^2 }}
\!\!\!\!\!\!\!\!\!\!\!\!\!\!\! 1  \nonumber \\ & & \nonumber \\ &
\leq & \frac{\Psi(2^{m})}{2^{m}} \sum_{2^{m-1} < s \leq a 2^{m+1}
} \ \ \  \sum_{\substack{n: \\ |s - \sqrt{n}| < 8\Psi(s/a) }}
\!\!\!\!\! r(n) \nonumber \\ & & \nonumber \\ & \leq &
\frac{\Psi(2^{m})}{2^{m}} \ \ \sum_{i=0}^{k+2} \ \ \ \
\sum_{2^{m+i-1} < s \leq 2^{m+i} } \ \ \ \sum_{\substack{n: \\
|s - \sqrt{n}| < 8\Psi(s/a) }} \!\!\!\!\! r(n) \nonumber \\ & &
\nonumber \\ & \ll &  k \ \frac{\Psi(2^{m})}{2^{m}} \ 2^{2(m+k+1)}
\  \Psi(2^{m-k-2}) \ \ll \ k \ 2^{3k} \ 2^{m-k-2} \ \Psi
(2^{m-k-2})^2 \ .
\end{eqnarray*}
Thus, $  \sum \big| W_{m}(\Psi;\cC_{1;k}^*) \big|_{\cC_{1}^*} \ll
\sum 2^m \ \Psi(2^{m})^2 \ \asymp \ \sum \Psi(h)^2 \ < \infty $
and we are done.

\medskip

\noindent{\bf Case (c): $\cQ = \Veronese_2$. } For $k \in \N $,
let $\Veronese_{2;k} := \{(x_1,x_2) \in \R^2 : x_2 = x_1^2  {\rm \
with \ } |x_1| \leq 2^k  \} $. For $m \in \N $, let $$
W_{m}(\Psi;\Veronese_{2;k} ) \ := \ \bigcup_{2^{m} < q \leq
2^{m+1} } \bigcup_{(s,t)\in \Z^2 } \Veronese_{2;k} \cap C(q;s,t) \
 \ . $$
We need to show that  $  \sum \big| W_{m}(\Psi;\Veronese_{2;k})
\big|_{\Veronese_{2}} < \infty $.
 It is easily verified that if $\Veronese_{2;k} \cap
C(q;s,t) \neq \emptyset$ then $0 \leq  |s|/q < a:=2^{k+1} $, $ -1
< t/q < a^2 $ and
$|s^2 \, - \, tq | \, < \,   2\, \Psi(q) (2\,  |s| + |t|) \, + \,
4 \,  \Psi(q)^2 \, < \, 6 \, a^2 q \Psi(q) + 4 \, \Psi(q)^2
$; 
that is,
\begin{eqnarray}
|(2s)^2 \, - \, 4tq | \ < \    24 \, a^2 q \Psi(q) + 16 \,
\Psi(q)^2 \ . \label{4x}
\end{eqnarray}
Let $w:=q+t$ and $z:=q-t$. Then, $2q = w+z $,  $2t=w-z$  and $q-1
< w < q(a^2+1) $.  Furthermore, (\ref{4x}) becomes
\begin{eqnarray}
|(2s)^2 \, + \, z^2 - w^2| \ < \    24 \, a^2 q \Psi(q) + 16 \,
\Psi(q)^2 \ < \ 48 \, a^2 \, w \,
\Psi\left(\textstyle{\frac{w}{(a^2+1)}} \right) \, + \, 16 \,
\Psi\left(\textstyle{\frac{w}{(a^2+1)}} \right)^2 \ . \label{4xx}
\end{eqnarray}
It follows, that for $m$ sufficiently large
\begin{eqnarray*}
\big| W_{m}(\Psi;\Veronese_{2;k}) \big|_{\Veronese_2} & \ll &
\frac{\Psi(2^{m})}{2^{m}}  \sum_{2^{m} < q \leq 2^{m+1} } \ \ \
\sum_{\substack{(s,t)\in \Z^2 \; :  \ -q < t < a^2 q  \\ {\rm
(\ref{4x}) \ holds }
 }}
\!\!\!\!\!\!\!\!\!\!\!\!\!\!\! 1 \nonumber \\ & & \nonumber \\ &
\leq &  \frac{\Psi(2^{m})}{2^{m}} \sum_{2^{m-1} < w \leq (a^2+1)
2^{m+1} } \ \ \  \sum_{\substack{ (s,z)\in \Z^2 \; : \;  {\rm
(\ref{4xx}) \ holds} }
 } \!\!\!\!\! 1 \nonumber \\ & & \nonumber \\ &
\leq & \frac{\Psi(2^{m})}{2^{m}} \ \sum_{2^{m-1} < w \leq a^2
2^{m+2} } \ \ \ \sum_{\substack{n: \\
|w - \sqrt{n}| < 48\Psi(w/(2a^2)) }} \!\!\!\!\! r(n)  \ .
\end{eqnarray*}
As in  case (b), the desired statement now follows on using
(\ref{thebigone}) to estimate the  double sum. \hfill $\spadesuit$

\medskip

Before moving onto the proof of Theorem \ref{circ2}, we  define
Hausdorff measure and dimension for the sake of  completeness and
in order to establish some notation.

\subsection{Hausdorff measure and dimension \label{HM}}

The Hausdorff dimension of a non--empty subset $X$ of
$n$--dimensional Euclidean space $\R^{n}$, is an aspect of the
size of $X$ that can discriminate between sets of Lebesgue measure
zero.

\smallskip

For $\rho > 0$, a countable collection $ \left\{C_{i} \right\} $
of Euclidean cubes in  $\R^{n}$ with side length $l(C_i)  \leq
\rho $ for each $i$  such that $X \subset \bigcup_{i} C_{i} $ is
called a $ \rho $-cover for $X$.  Let $s$ be a non-negative number
and define $$
 {\cal H}^{s}_{ \rho } (X)
  \; = \; \inf \left\{ \sum_{i} l^s (C_i)
\ :   \{ C_{i} \}  {\rm \  is\ a\  } \rho {\rm -cover\  of\ } X
\right\} \; , $$ where the infimum is taken over all possible $
\rho $-covers of $X$. The {\it s-dimensional Hausdorff measure}
${\cal H}^{s} (X)$ of $X$ is defined by $$ {\cal H}^{s} (X) =
\lim_{ \rho \rightarrow 0} {\cal H}^{s}_{ \rho } (X) = \sup_{ \rho
> 0} {\cal H}^{s}_{ \rho } (X)
$$ \noindent and the {\it Hausdorff dimension} dim $X$ of $X$ by
$$ \dim \, X = \inf \left\{ s : {\cal H}^{s} (X) =0 \right\} =
\sup \left\{ s : {\cal H}^{s} (X) = \infty \right\} \, . $$ \vskip
9pt

Strictly speaking, in the standard definition of Hausdorff measure
the  $\rho$--cover by cubes is replaced by non--empty subsets in
$\R^n$ with diameter at most $\rho \,$. It is easy to check that
the resulting measure is comparable to ${\cal H}^{s} $ defined
above and thus the Hausdorff dimension is the same in both cases.
For our purpose using cubes is just more convenient. Moreover, if
${\cal H}^{s} $ is zero or infinity then there is no loss of
generality by restricting to cubes.   Further details and
alternative definitions of Hausdorff measure and dimension can be
found in \cite{falc,MAT}.

\subsection{Proof of Theorem \ref{circ2}} To a certain degree the proof
follows the same line of argument as the proof of the convergent
part of Theorem \ref{corol2C}.  In particular, it suffices  to
consider the  rational quadrics $\cC_1$, $\cC_1^*$ and
$\Veronese_2$. Below, we consider the case of the unit circle
$\cC_1$ and leave the hyperbola $\cC_1^*$ and parabola $\Veronese_2$
to the reader. The required modifications are obvious.

Let $\psi$ be an approximating function such that $\sum h^{1-s} \,
\psi(h)^{s+1} <\infty $ and consider the auxiliary function $\Psi
$ given by $$ \Psi(h) \; := \; \max \left\{ \psi(h), \ h^{-1} \;
(\log h)^{260} \right\} \ \ .
$$ Clearly, $\Psi$ is an approximating function and since $s >
1/2$ we have that $ \sum h^{1-s} \, \Psi(h)^{s+1} <\infty   $.
With the same notation as in the proof of Theorem \ref{corol2C},
for each $l \in \N$ $$ \left\{ W_m(\Psi,\cC_1) : m=l, \, l+1,
\ldots  \, \right\}
$$ is a cover for $ \cC_1(\Psi) := \cC_1\cap\cS_2(\p)$ by squares $C(q;s,t)$ of maximal
side length $2 \Psi(2^l)/2^l $. It follows from the definition of
$s$--dimensional Hausdorff measure  that with $\rho := 2
\Psi(2^l)/2^l $
\begin{eqnarray*}
{\cal H}^s_{\rho} ( \cC_1(\Psi)) & \leq  & \sum_{m=l}^{\infty} \ \
\  \sum_{2^{m} < q \leq 2^{m+1} } \sum_{\substack{(s,t)\in \Z^2:
\\ (q - 2\sqrt2 \Psi(q))^2 \, \leq s^2 + t^2 \, (q + 2\sqrt2
\Psi(q))^2 }} \!\!\!\!\!\!\!\!\!\!\!\!\!\!\!
\left(\frac{2\Psi(2^m)}{2^m} \right)^s \nonumber
\\ & & \nonumber \\ & \ll & \sum_{m=l}^{\infty}
\left(\frac{\Psi(2^{m})}{2^{m}}\right)^s \sum_{2^{m} < q \leq
2^{m+1} } \sum_{\substack{n: \\ |q - \sqrt{n}| < 4\Psi(q) }}
\!\!\!\!\! r(n)
 \ \ \ .
\end{eqnarray*}

\noindent In view of  Theorem A of Appendix II, the contribution
from the two inner sums is $\ll   2^{2m} \Psi(2^m) $. Thus, $$
{\cal H}^s_{\rho} ( \cC_1(\Psi)) \ \ll \
 \sum_{m=l}^{\infty}  2^{m(2-s)} \ \Psi(2^{m})^{1+s} \ \to \ 0 $$
 as $\rho \to 0$; or equivalently at $l \to \infty$.   Hence,
 $ {\cal H}^s ( \cC_1(\psi)) \leq  {\cal H}^s ( \cC_1(\Psi)) = 0 $
  as required.
\hfill $\spadesuit$


\mysection{Ubiquitous systems \label{ub}}

In \cite{BDV03}, a general framework is developed for establishing
divergent results analogous to those of Khintchine and Jarn\'{\i}k
(see \S\ref{intro}) for a natural class of $\limsup $ sets.  The
framework is based on the notion of `ubiquity', which goes back to
\cite{BS70} and \cite{DRV} and captures the key measure theoretic
structure necessary to prove such measure theoretic laws. The
`ubiquity' introduced below is  a much simplified version  of that
in \cite{BDV03} and  takes into consideration the specific
applications that we have in mind.

\subsection{Ubiquitous systems in $\R$ \label{ubin1}}

Let $I_0$ be an interval in $\R$ and 
$\cR:=(\ra)_{\al\in\cJ}$  be a family of {\em resonant points}
$\ra$ of $I_0$ indexed by an infinite, countable  set $\cJ$. Next
let $\beta:\cJ\to \R^+:\alpha\mapsto\ba$ be a positive function on
$\cJ$.  Thus, the function $\beta$ attaches a `weight' $\ba$  to
the resonant point $\ra$. Also, for $t \in \N$ let  $J(t):=\{\al
\in \cJ:\ba\le 2^t\}$ and assume that  $\#J(t)$ is always finite.
Given an approximating  function $\Psi$ let
\begin{equation*}\label{e:004}
\Lambda(\cR,\beta,\Psi) \ := \ \{x\in I_0 :|x-\ra|<\Psi(\ba) \
\mbox{for\ infinitely\ many\ }\al \in \cJ \} \ .
\end{equation*}
The set $\Lambda(\cR,\beta,\Psi)$ is easily seen to be a $\limsup$
set. The general theory of ubiquitous systems developed in
\cite{BDV03}, provides a natural measure theoretic condition for
establishing divergent results analogous to those of Khintchine
and Jarn\'{\i}k  for $\Lambda(\cR,\beta,\Psi)$. Since
$\Lambda(\cR,\beta,\Psi)$ is a subset of $I_0$, any Khintchine
type result would naturally be with respect to one-dimensional
Lebesgue measure $|\, .\, |$.

Throughout, $\rho: \Rp \to\Rp$ will denote a function  satisfying
$\lim_{t\to\infty}\rho(t)=0$ and is usually referred to as the {\em
ubiquitous function}. Also $B(x,r)$ will denote the ball (or rather
the interval) centred at $x$ or radius $r$.

\begin{definition}[Ubiquitous systems on the real line]\sl
\label{US} Suppose there exists a function $\rho$ and an absolute
constant $\kappa > 0$ such that for any interval $I\subseteq I_0$
\begin{equation*}
\liminf_{t\to\infty}\,\,\left|\,{\textstyle\bigcup_{\al\in J(t)}}
\left(B(\ra,\r(2^t)\right)\cap I)\right| \ \ge \  \ka \, |I| \ .
\end{equation*}
Then the system $(\cR;\beta)$ is called {\em locally ubiquitous in
$I_0$ with respect to $\rho$}.
\end{definition}

The consequences of this definition of ubiquity are the following
key results.

\begin{lemma}\label{KT}
Suppose that $(\cR,\beta)$ is a local ubiquitous system in  $I_0$
with respect to $\rho$ and let $\Psi$ be an approximating function
such that $\Psi(2^{t+1})\le \frac12 \Psi(2^t)$ for $t$ sufficiently
large. Then
\begin{equation*}
|\Lambda(\cR,\beta,\Psi)| \ = \ \mbox{\rm F{\scriptsize ULL}}
 \ := \ |I_0| \hspace{10mm} {\rm if } \hspace{10mm}
\sum_{t=1}^\infty\frac{\Psi(2^t)}{\rho(2^t)}=\infty \ \ .
\end{equation*}
\end{lemma}


\begin{lemma}\label{JT}
Suppose that $(\cR,\beta)$ is a local ubiquitous system in  $I_0$
with respect to $\rho$ and let $\Psi$ be an approximating
function. Let $s \in (0,1)$ and let $$  G \, := \, \limsup_{t \to
\infty} \, \frac{\Psi(2^t)^s}{\rho(2^t)}
\; \ . $$

\noindent {\bf (i)}  Suppose that $G = 0$  and  that
$\Psi(2^{t+1})\le \frac12 \Psi(2^t)$ for $t$ sufficiently large.
Then, $$
 {\cal H}^s \! \left( \Lambda(\cR,\beta,\Psi) \right) \ =
\ \infty \hspace{10mm} {\rm if } \hspace{10mm} \sum_{t=1}^{\infty}
\frac{\Psi(2^t)^s}{\rho(2^t)}\ = \ \infty \ \ . $$

\noindent {\bf (ii)} Suppose that $G >0 $. Then,  ${\cal H}^s \!
\left( \Lambda(\cR,\beta,\Psi) \right) \ = \ \infty$.

\end{lemma}



\begin{corollary}\label{LB}
Suppose that $(\cR,\beta)$ is a local ubiquitous system in  $I_0$
with respect to $\rho$ and let $\Psi$ be an approximating
function. Then $$ \dim(\Lambda(\cR,\beta,\Psi)) \ \ge \   d \, :=
\, \min\left\{1, \left| \limsup_{t\to\infty} \frac{\log \r(2^t)}{
\log \Psi(2^t)} \right| \right\}. $$ Moreover,  if $d <1$ and
 $ \limsup_{t \to \infty} \Psi(2^{t})^d/ \rho(2^{t})  > 0 $,
 then $ {\cal H}^d(\Lambda(\cR,\beta,\Psi)) = \infty $.
\end{corollary}

The concept  of ubiquity was originally formulated by Dodson, Rynne
$\&$ Vickers \cite{DRV} to obtain lower bounds for the Hausdorff
dimension of $\limsup$ sets. In the one-dimensional setting
considered here, their `ubiquity result'  essentially corresponds to
Corollary \ref{LB} above. Furthermore, the ubiquitous systems of
\cite{DRV} essentially coincide with the regular systems of Baker
$\&$ Schmidt \cite{BS70} and both have proved very useful in
obtaining lower bounds  for the Hausdorff dimension of $\limsup$
sets. However, unlike the framework developed in \cite{BDV03},  both
\cite{BS70} and \cite{DRV}  fail to shed any light on establishing
the more desirable divergent Khintchine and Jarn\'{\i}k type
results. The latter, clearly implies lower bounds for the Hausdorff
dimension. For further details regarding regular systems and the
original formulation of ubiquitous systems see \cite{BDV03,BD99}.

\smallskip

Lemmas \ref{KT} and \ref{JT} follow directly from Corollaries 2
and 4 in \cite{BDV03}. Note that in Lemma \ref{JT}, if $G >0$ then
the divergent sum condition of part (i) is trivially satisfied.
The
 dimension statement (Corollary \ref{LB}) is a consequence of part (ii) of
Lemma \ref{JT} and so the regularity condition $ \textstyle 2 \,
\Psi(2^{t+1})\le
 \Psi(2^t)$  on the function
$\Psi$ is not necessary -- see \cite[Corollary 6]{BDV03}.

The framework and results of \cite{BDV03} are  abstract and
general unlike the concrete  situation described above.  In view
of this and for the sake of completeness we retraced the argument
of  \cite{BDV03} in the above simple setting  at the end of the
paper \S\ref{pf12}-\ref{proof3}. This has the effect of making the
paper self-contained and more importantly should help the
interested reader with understanding the abstract approach
undertaken in \cite{BDV03}. The direct proofs of Lemmas \ref{KT}
and \ref{JT} are substantially easier (both technically and
conceptionally) than the general statements of \cite{BDV03}.

\subsection{Ubiquitous systems close to a curve in $\R^n$ \label{usnc} }

In this section we develop the theory of ubiquity to incorporate
the situation in which the resonant points of interest lie within
some specified neighborhood of a given curve in $\R^n$.

With $n \geq 2$, let  $\cR:=(\ra)_{\al\in\cJ}$ be  a family of {\em
resonant points} $\ra$ of $\R^n$ indexed by an infinite set $\cJ$.
As before,  $\beta:\cJ\to \R^+:\alpha\mapsto\ba$ is  a positive
function on $\cJ$. For a point  $\ra$ in $\cR$,  let $\rak$
represent the $k$'th coordinate of $\ra$.  Thus, $\ra :=
(R_{\alpha,1}, R_{\alpha,2}, \ldots, R_{\alpha,n})$. Throughout this
section and the remainder of the paper we will use the notation
$\cR_{\cC}(\Phi)$ to denote the sub-family of resonant points $\ra$
in $\cR$ which are ``$\Phi$--close'' to the curve $\cC=\cC_{\vv
f}:=\{(x,f_2(x),\dots,f_n(x)):x\in I_0\}$ where $\Phi$ is an
approximating function, $\vv f=(f_1,\dots,f_n):I_0\to\R^n$ is a
continuous map with $f_1(x)=x$ and $I_0$ is an interval in $\R$.
Formally, and more precisely  $$ \cRp := ( \ra )_{ \alpha\in
\cJ_\cC(\Phi) } \hspace{7mm} \text{where} \hspace{5mm}
\cJ_\cC(\Phi):= \{\alpha\in\cJ: \max\limits_{1\le k\le n}|
f_k(\rao)-\rak|<\Phi(\ba)\} \ . $$

\noindent Finally, we  will  denote by $\cR_1$ the family of first
co-ordinates of the points in $\cRp $; that is $$ \cR_1 \ := \ (
R_{\alpha,1} )_{\alpha\in\cJp } \ \  . $$ \noindent By definition,
$\cR_1 $ is a subset of the interval $I_0$ and can therefore be
regarded as a set of resonant points for the theory of ubiquitous
systems in $\R$.  This leads us naturally to the following
definition in which the ubiquity function 
$\rho$ is as  in \S\ref{ubin1}.

\begin{definition}[Ubiquitous systems near curves]\sl \label{USNC}
The system $(\cRp,\beta)$ is called locally  ubiquitous with
respect to $\rho$
if the system $(\cR_1,\beta)$ is locally ubiquitous in $I_0$ with
respect to $\rho$.
\end{definition}

Next, given an approximating function $\Psi$ let
$\Lambda(\cRp,\beta,\Psi)$ denote the  the set  $x\in I_0$ for
which the system of inequalities
\[
\left\{\begin{array}{rcl} |x-\rao|&<&\Psi(\ba)  \\[0.5ex]
\max\limits_{2\le k\le
n}|f_k(x)-R_{\alpha,k}|&<&\Psi(\ba)+\Phi(\ba) \ \ ,
\end{array}\right.
\]
is simultaneously satisfied  for infinitely many $\alpha\in\cJ$.
The following two lemmas are the analogues of Lemmas \ref{KT} and
\ref{JT} for the case of ubiquitous systems close to a curve.
Similarly, Corollary \ref{hdl} is the analogue of Corollary
\ref{LB}.

\begin{lemma}\label{ktl}
Consider  the curve $\cC:=\set{(x,f_2(x),\dots,f_n(x)):x\in I_0}$,
where $f_2,\dots,f_n$ are locally Lipshitz in a finite interval
$I_0$. Let $\Phi$ and $\Psi$ be approximating functions. Suppose
that $(\cRp,\beta)$ is a locally ubiquitous system with respect to
$\rho$. If $\Psi$ and $\r$ satisfy the conditions of Lemma
\ref{KT} then $$|\, \Lambda\left(\cRp,\beta,\Psi\right) \, | \ = \
|I_0| \ \ .$$
\end{lemma}

\begin{lemma}\label{JTl}
Consider  the curve $\cC:=\set{(x,f_2(x),\dots,f_n(x)):x\in I_0}$,
where $f_2,\dots,f_n$ are locally Lipshitz in a finite interval
$I_0$. Let $\Phi$ and $\Psi$ be approximating functions. Suppose
that  $(\cRp,\beta)$ is a locally ubiquitous system with respect
to $\rho$. Let $s \in (0,1)$ and let $$  G \, := \, \limsup_{t \to
\infty} \, \frac{\Psi(2^t)^s}{\rho(2^t)} \; \ . $$ \noindent {\bf
(i)}  Suppose that $G =0$  and that $\Psi(2^{t+1})\le \frac12
\Psi(2^t)$ for $t$ sufficiently large. Then, $$
 {\cal H}^s \! \left( \Lambda(\cRp,\beta,\Psi) \right) \ =
\ \infty \hspace{10mm} {\rm if } \hspace{10mm} \sum_{t=1}^{\infty}
\frac{\Psi(2^t)^s}{\rho(2^t)}\ = \ \infty \ \ . $$
\noindent {\bf
(ii)} Suppose that $G >0$. Then,  ${\cal H}^s \! \left(
\Lambda(\cRp,\beta,\Psi) \right) \ = \ \infty$.
\end{lemma}

\begin{corollary}\label{hdl}
Consider  the curve $\cC:=\set{(x,f_2(x),\dots,f_n(x)):x\in I_0}$,
where $f_2,\dots,f_n$ are locally Lipshitz in a finite interval
$I_0$.  Let $\Phi$ and $\Psi$ be  approximating functions. Suppose
that $(\cRp,\beta)$ is a locally ubiquitous system with respect to
$\rho$.
Then $$ \dim\Lambda \left(\cRp,\beta,\Psi \right) \ \ge \   d \,
:= \, \min\left\{1, \left| \limsup_{t\to\infty} \frac{\log
\r(2^t)}{ \log \Psi(2^t)} \right| \right\}. $$ Moreover,  if $d
<1$  and $ \limsup_{t \to \infty} \Psi(2^{t})^d/ \rho(2^{t})  > 0
$, then $ {\cal H}^d(\Lambda \left(\cRp,\beta,\Psi \right)) =
\infty $.
\end{corollary}

\noindent{\bf Proof of Lemmas~\ref{ktl} \& \ref{JTl}  and
Corollary \ref{hdl}.} It suffices to prove the lemmas for a
sufficiently small neighborhood of a fixed point in $I_0$.
Therefore, there is no loss of generality in assuming that
$f_2,\dots,f_n$ satisfy the Lipshitz condition on $I_0$. Thus, we
can fix a constant $c_3\ge 1$ such that for  $k\in\set{2,\dots,n}$
and $x,y\in I_0$
\begin{equation}\label{e:010}
  |f_k(x)-f_k(y)|\le c_3|x-y|.
\end{equation}

Since $(\cRp,\beta)$ is a locally ubiquitous system with respect
to $\rho$, by definition $(\cR_1,\beta)$ is a locally ubiquitous
system in $I_0$ with respect to $\rho$. The set
$\Lambda(\cR_1,\beta,\Psi/c_3)$ consists of $x\in I_0$ for which
the inequality
\begin{equation}\label{e:011}
  |x-\rao|<\Psi(\ba)/c_3\le \Psi(\ba)
\end{equation}
is satisfied  for infinitely many $\alpha\in\cJp$. Suppose $x$
satisfies (\ref{e:011}) for some $\alpha\in\cJp$. In view of
(\ref{e:010}), $|f_k(x)-f_k(R_{\alpha,1})|\le c_3|x-R_{\alpha,1}|$
which implies that
\begin{equation*}
\begin{array}[t]{rcl}
  |f_k(x)-\rak| & = & |f_k(x)-f_k(\rao)+f_k(\rao)-\rak|
  \\[1ex]
   & \le & |f_k(x)-f_k(\rao)|+|f_k(\rao)-\rak|
   \\[1ex]
   & \le & c_3|x-\rao|+\Phi(\ba) \\[1ex]
   & < & c_3\cdot\Psi(\ba)/c_3+\Phi(\ba)= \Psi(\ba)+\Phi(\ba).
\end{array}
\end{equation*}
Thus
$\Lambda(\cR_1,\beta,\Psi/c_3)\subset\Lambda(\cR,\beta,\Psi)$.
Applying Lemmas~\ref{KT} \& \ref{JT} and Corollary \ref{LB} to the
set $\Lambda(\cR_1,\beta,\Psi/c_3)$ gives the desired statements
concerning the set  $\Lambda(\cRp,\beta,\Psi)$. \hfill
$\spadesuit$

\mysection{Proof of Theorem \ref{thm3}}

We begin by stating a key result which not only implies
Theorem~\ref{thm3} but gives rise to a ubiquitous system that will
be required in proving Theorems~\ref{thm1} and \ref{thm2a}.

\subsection{The ubiquity version of Theorem~\ref{thm3}}

\begin{theorem}\label{thm4}
Let $I_0$ denote a finite, open interval of $\R$ and let $f$ be a
function in $C^{(3)}(I_0)$ satisfying (\ref{infsup}). Let $\psi$
be an approximating function satisfing {\rm(\ref{e:002})}. Then
for any interval $I\subseteq I_0$ there exist constants
$\delta_0,C_1>0$ such that for $Q$ sufficiently large$$
\left|\bigcup_{\bp/q\in A_Q(I)}
\left(B\left(\frac{p_1}{q},\frac{C_1}{Q^2 \p(Q)}\right)\cap
I\right)\right| \ \ge \  \frac{1}{2} \, |I| \ \ ,$$ where
$$A_Q(I)\ :=\ \left\{\bp/q\in\Q^2 \, : \, \delta_0Q<q\le Q, \
p_1/q\in I \, ,  \  |f(p_1/q) - p_2/q|< \p(Q)/Q \right\} \ \ . $$
\end{theorem}

\noindent{\bf Proof of Theorem~\ref{thm3} modulo
Thorem~\ref{thm4}.  } This is trivial. Given the hypotheses of
Theorem \ref{thm4}, the hypotheses of Theorem \ref{thm3} are
clearly satisfied. Fix an interval $I\subseteq I_0$.  By Theorem
\ref{thm4}, there exist constants $\delta_0$ and $C_1$ so that for
all $Q$ sufficiently large $$\#A_Q(I) \cdot
\frac{2C_1}{Q^2\psi(Q)} \geq \sum_{\bp/q\in A_Q(I)} \left|
B\left(\frac{p_1}{q},\frac{C_1}{Q^2\p(Q)}\right) 
\right|\geq \left|\bigcup_{\bp/q\in A_Q(I)} \!\!
\left(B\left(\frac{p_1}{q},\frac{C_1}{Q^2\p(Q)}\right) \cap I\right)
\right|\ge\frac{|I|}{2}  .  $$  We have that $ N_f(Q,\psi,I)\ge
\#A_Q(I)$ and Theorem \ref{thm3} follows. \hfill $\spadesuit$

\medskip

The following corollary of Theorem~\ref{thm4} is crucial for
proving Theorems~\ref{thm1} and \ref{thm2a}.

\begin{corollary}\label{cor1}
Let $\psi$ and $f$ be  as in Theorem~\ref{thm4} and
$\cC:=\{(x,f(x)) :x\in I_0\}$.   With reference to the ubiquitous
framework of \S\ref{usnc}, set
\begin{equation}\label{beta}
\beta:\ \cJ := \Z^2 \times \N \to\N : (\bp,q) \to q  \ ,
\hspace{5mm} \Phi: t \to t^{-1}\psi(t) \hspace{5mm} \mbox{and}
\hspace{5mm} \rho :t \to u(t)/(t^2\psi(t))
\end{equation}
where $u:\Rp\to\Rp$ is any function such that
$\lim_{t\to\infty}u(t)=\infty$. Then the system
$(\Q^2_{\cC}(\Phi),\beta)$ is  locally  ubiquitous with respect to
$\rho$. 
\end{corollary}

\noindent{\bf Remark. \ } Given $\alpha=(\bp,q) \in \cJ$, the
associated resonant point $\ra$ in the above ubiquitous system is
simply the rational point $ \bp/q$ in the plane.  Furthermore, $\cR
:=  \Q^2$.

\vspace{2mm}

\noindent{\bf Proof of Corollary \ref{cor1}. \ } 
For an interval $I\subseteq I_0$, let $$\waq \ := \ \{\bp/q\in\Q^2
\, : \, Q/u(Q)<q\le Q \, , \ p_1/q\in I  \, ,  \
|f(p_1/q)-p_2/q|<\p(Q)/Q\} \ \ . $$ For any $\delta_0 \in (0,1)$,
we have that  $1/u(Q) < \delta_0$   for $Q$ sufficiently large
since $\lim_{t\to\infty}u(t)=\infty$. Thus, for $Q$ sufficiently
large,  $A_Q(I)\subset \waq$ and  Theorem \ref{thm4} implies that
\[\left|\bigcup_{\bp/q\in\waq}\left(B
\left(\frac{p_1}{q},\frac{u(Q)}{Q^2\p(Q)}\right)\cap
I\right)\right| \ge \left|\bigcup_{\bp/q\in A_Q(I)}\left(B
\left(\frac{p_1}{q},\frac{C_1}{Q^2\p(Q)}\right)\cap
I\right)\right| \ge \frac{|I|}{2}.\]

\noindent This establishes the corollary. \hfill $\spadesuit$

\subsection{An auxiluary lemma}

The following lemma is an immediate consequence of Theorem 1.4 in
\cite{BKM01}.

\begin{lemma}\label{BKM0}
Let $\vv g :=(g_1,g_2):I_0\to\R^2$ be a $C^{(2)}$ map such that
$(g_1'g_2''-g_2'g_1'')(x_0)\neq0$  for some point $x_0\in I_0$.
Given positive real numbers $\delta,K,T$ and an interval
$I\subseteq I_0$, let $B(I,\delta,K,T)$ denote the set  of $x\in
I$ for which there exists $(q,p_1,p_2)\in\Z^3\smallsetminus\{0\}$
satisfying the following system of inequalities:
\begin{equation*}
\left\{
\begin{array}{l}
|q \, g_1(x) \, + \, p_1 \, g_2(x)+ p_2| \ \le \ \delta
 \\ \\
|q \, g_1'(x) \, + \, p_1 \, g_2'(x)| \ \le \ K
 \\  \\
|q| \ \le \  T \ \  .
\end{array}
\right.
\end{equation*}
Then there is a sufficiently small $\eta=\eta(x_0)>0$ so that for
any interval $I\subset(x_0-\eta,x_0+\eta)$ there exists a constant
$C>0$ such that for
\begin{equation}\label{dkt}
0<\delta\le 1,\quad T\ge1, \quad K>0\quad \text{and} \quad
 \delta KT\le1
\end{equation}
one has
\begin{equation}\label{est}
|{B(I,\delta,K,T)}|\le C \max\left(\delta^{1/3}, \left(\delta K T
\right)^{1/9}\right)|I|.
\end{equation}
\end{lemma}

Note that the constant  $C$ depends on the interval $I$. We now
show that under the assumption that $\vv g$ is non-degenerate
everywhere,  the above  lemma can  be extended to a  global
statement in which $I$ is any sub-interval of $I_0$.

\begin{lemma}\label{BKM}
Assume that the conditions of Lemma~\ref{BKM0} are satisfied and
that $(g_1'g_2''-g_2'g_1'')(x)\neq0$ for all $x\in I_0$. Then for
any finite interval $I\subseteq I_0$ there is a constant $C>0$
such that for any $\delta,K,T$ satisfying\/ {\rm(\ref{dkt})}\/ one
has the estimate\/ {\rm(\ref{est})}.
\end{lemma}

\noindent {\bf Proof of Lemma \ref{BKM}. \ } As $I$ is a finite
interval, its closure $\overline I$ is compact. By
Lemma~\ref{BKM0}, for every point $x\in\overline I$ there is an
interval $B(x,\eta(x))$ centred at $x$ such that for any
sub-interval $J$ of $B(x,\eta(x))$ there is a constant $C=C_J$
(dependent  on $J$) satisfying (\ref{est}) with $\delta,K,T$
satisfying (\ref{dkt}). Since  $\overline I$ is compact, there is
a finite cover $\{I_i:= B(x_i,\eta(x_i)) : i=1,\dots,n \} $ of
$\overline I$. Choose this cover so that $n$ is minimal. Then any
interval in this cover is not contained in the union of the
others. Otherwise, we would be able to choose another cover with
smaller $n$. We show that any three intervals of this minimal
cover do not intersect. Assume the contrary. So there is an
$x\in(a_1,b_2)\cap(a_2,b_2)\cap(a_3,b_3)$,  where $(a_i,b_i)$,
$i=1,2,3$ are  intervals of the  minimal cover. Then $a_i< x<b_i$
for each $i$. Without loss of generality, assume that $a_1\le
a_2\le a_3$. If $b_2<b_3$ then
$(a_2,b_2)\subset(a_1,b_3)=(a_1,b_1)\cap(a_3,b_3)$, which
contradicts the minimality of the cover. Similarly, if $b_3\le
b_2$ then $(a_3,b_3)\subset(a_1,b_2)=(a_1,b_1)\cap(a_2,b_2)$, a
contradiction. This means that the multiplicity of the cover is at
most $2$. Hence $ \sum_{i=1}^n|I_i|\le 2|I|, $ where $I_i :=
B(x_i,\eta(x_i)$. This together with  Lemma~\ref{BKM0} implies
that
\begin{eqnarray*}
|B(I,\delta,K,T)|& = &   |\; \textstyle{\bigcup_{i=1}^{n}}
B(I_i,\delta,K,T) \; |
 \ \leq \ \sum_{i=1}^n
|B(I_i,\delta,K,T)| \\ &\le & \textstyle{\sum_{i=1}^n} C_{I_i} \
\max\left(\delta^{1/3}, \left(\delta K T \right)^{1/9}\right)|I_i|
\\ & \le & \max_{i=1,\dots,n} C_{I_i}\cdot\max\left(\delta^{1/3},
\left(\delta K T \right)^{1/9}\right) \ \textstyle{\sum_{i=1}^n}
\; |I_i| \\ & \le & 2\max_{i=1,\dots,n}
C_{I_i}\cdot\max\left(\delta^{1/3}, \left(\delta K T
\right)^{1/9}\right)|I|  \ ,
\end{eqnarray*}
as required. \hfill $\spadesuit$

\subsection{Proof of Theorem~\ref{thm4}}

Define $\vv g(x) :=(g_1(x),g_2(x))$ by setting $g_1(x):=x
f'(x)-f(x)$ and $g_2(x):=-f'(x)$. Then $\vv g\in C^{(2)}$. Also,
note that
\begin{equation}\label{e:013}
\vv g'(x)=(xf''(x),\,-f''(x)) \;  ,\qquad \vv
g''(x)=(f''(x)+xf'''(x),\,-f'''(x))
\end{equation}
and $$ (g_1'g_2''-g_2'g_1'')(x)=f''(x)^2 \ . $$ As $f''(x)\neq0$
everywhere, Lemma~\ref{BKM} is applicable to this $\vv g$. In view
of  the conditions on the theorem,
\begin{equation}\label{e:015}
\sup_{x\in I_0}|g_2'(x)|=\sup_{x\in I_0}|f''(x)|\le c_2.
\end{equation}
Define $\delta_0:=\min\{1,(2^{19}c_2C^9)^{-1}\}$, where $C$ is the
constant appearing in Lemma~\ref{BKM}. Without loss of generality,
assume that $C>1$.

 Next,  fix an interval
$I\subseteq I_0$.  By Minkowski's linear forms theorem in the
geometry of numbers, for any $x\in I$ and $Q \in \N$ there is a
solution $(q,p_1,p_2)\in\Z^3\smallsetminus\{0\}$ to the system
\begin{equation}\label{e:016}
\left\{
\begin{array}{l}
 |q^{}g_1(x)+p_1g_2(x)+p_2|\le\delta_0\psi(Q)
 \\ \\
 |q^{}g_1'(x)+p_1g_2'(x)|\le c_2(\delta_0 Q\psi(Q))^{-1}
 \\ \\
 0\le q\le Q \ .
\end{array}
\right.
\end{equation}

\noindent By definition, the set $B(I,\delta,K,T)$ with
\begin{equation}\label{e:017}
\delta :=\delta_0\psi(Q),\ \  \  K:=c_2(\delta_0 Q\psi(Q))^{-1},\
\ \ T:=2\delta_0 Q
\end{equation}
consists of points $x\in I$ such that there exists a non-zero
integer solution $(q,p_1,p_2)$ to the system  (\ref{e:016}) with
$q\le2\delta_0 Q$. 
By Lemma~\ref{BKM}, for sufficiently large $Q$ we have that
\begin{eqnarray*}
|B(I,\delta,K,T)| & \le & C\, |I|\
\max\big\{(\delta_0\psi(Q))^{1/3} \, , \
\big(\delta_0\psi(Q)c_2(\delta_0 Q\psi(Q))^{-1}2\delta_0
Q\big)^{1/9}\big\} \\   & = & C\, (2c_2\delta_0)^{1/9}|I| \ \le \
|I|/4 \  \ .
\end{eqnarray*}
Therefore,  with $\delta,K,T$ given by (\ref{e:017}) and  $Q$
sufficiently large
\begin{equation}\label{e:018}
|{\textstyle\frac34}I\setminus B(I,\delta,K,T)|\ge|I|/2  \ ,
\end{equation}
where $\frac34I$ is the interval $I$ scaled  by $\frac34$. Notice,
that for  $x \in \frac34I\setminus B(I,\delta,K,T)$ we have that
\begin{equation}\label{e:019}
q>2\delta_0 Q
\end{equation}
for any solution $(q,p_1,p_2)$ of (\ref{e:016}). 
From now on, assume that $x \in \frac34I\setminus
B(I,\delta,K,T)$. In view of (\ref{e:013}) and the second
inequality of (\ref{e:016}) we have that
\begin{equation*}
 |q x f''(x)-p_1f''(x)|<c_2(\delta_0 Q\psi(Q))^{-1}.
\end{equation*}
This together with (\ref{e:019}) and the fact that $|f''(x)|>c_1$,
implies that
\begin{equation}\label{e:021}
\left|x-\frac{p_1}{q}\right|\ \le \ \frac{c_2}{q|f''(x)|\delta_0
Q\psi(Q)} \ < \ \frac{c_2}{c_1\delta_0^2
Q^2\psi(Q)}=\frac{C_1}{Q^2\psi(Q)} \ \ ,
\end{equation}
where $C_1: =\frac{c_2}{c_1\delta_0^2}$. In view of (\ref{e:002})
and the fact that $x \in \frac34I$, we have that $p_1/q\in I$ for
$Q$ is sufficiently large. By Taylor's formula, $$ \textstyle
f\big(\frac{p_1}{q}\big)=f(x)+f'(x)\big(\frac{p_1}{q}-x\big)+
\frac12f''(\tilde x)\big(\frac{p_1}{q}-x\big)^2  $$ for some
$\tilde x$ between $x$ and $p_1/q$. Thus $\tilde x \in I$. Now the
expression on the left hand side of  the first inequality of
(\ref{e:016}) is equal to $$
\begin{array}{l}
|q(xf'(x)-f(x))-p_1f'(x)+p_2|=|(qx-p_1)f'(x)+p_2-qf(x)|\\[2ex] =
\big|(qx-p_1)f'(x)+p_2-q\Big(
f\big(\frac{p_1}{q}\big)-f'(x)\big(\frac{p_1}{q}-x\big)-
\frac12f''(\tilde x)\big(\frac{p_1}{q}-x\big)^2 \Big)\big|
\\[2ex] =
\big|p_2-qf\big(\frac{p_1}{q}\big)+\frac{q}{2}f''(\tilde
x)(x-\frac{p_1}{q})^2\big|  \ .  
\end{array}
$$ It follows from (\ref{e:002}), (\ref{e:015}), (\ref{e:016}) and
(\ref{e:021}) that for $Q$ sufficiently large
\begin{eqnarray*}
\textstyle \big|qf\big(\frac{p_1}{q}\big)-p_2\big| &\le&
\big|p_2-qf\big(\frac{p_1}{q}\big)+\frac{q}{2}f''(\tilde
x)(x-\frac{p_1}{q})^2\big|+\big|\frac{q}{2}f''(\tilde
x)(x-\frac{p_1}{q})^2\big| \\ \\ & \le &
\delta_0\psi(Q)+\frac{Q}{2}c_2\Big(\frac{C_1}{Q^2\psi(Q)}\Big)^2 \
< \  2\delta_0\psi(Q)  \ \ \ .
\end{eqnarray*}
This inequality together with (\ref{e:019}) implies that
\begin{equation}\label{e:022}
\textstyle \big|f\big(\frac{p_1}{q}\big)-\frac{p_2}{q}\big| \ < \
\frac{2\delta_0\psi(Q)}{q} \ < \ \frac{\psi(Q)}{Q}  \ .
\end{equation}
Thus, for any $x\in \frac34I\smallsetminus B(I,\delta,K,T)$
conditions  (\ref{e:021}) and (\ref{e:022}) hold for some
$(p_1,p_2)/q$ with $2\delta_0 \, q <   q \leq Q$.  Thus, ${\bf
p}/q := (p_1,p_2)/q \in A_Q(I)  $ and moreover, in view of
(\ref{e:018}) we have that
\begin{equation*}
\Big|\textstyle{\bigcup_{{\bf p}/q \in A_Q(I)}}
\textstyle\big\{x\in
I:\big|x-\frac{p_1}q\big|<\frac{C_1}{Q^2\psi(Q)}\big\}\Big| \ \ge\
|I|/2 \ \ ,
\end{equation*}
for all sufficiently large $Q$. The statement of the theorem now
follows. \hfill $\spadesuit$

\mysection{Proof of Theorem \ref{thm2a} \label{thm2aproof} }

\noindent Throughout, $\psi$ is an approximating function with
$\lp := \liminf_{t\to\infty}\frac{-\log\p(t)}{\log t} \in
(1/2,1)$.  It is readily verified that for any $\ve>0$
\begin{equation}\label{p1}
\p(t)\le t^{-\lp+\ve}\qquad\text{for all but finitely many }t\in\N
\ \ ,
\end{equation}
and that  there exists a  strictly increasing sequence of natural
numbers $t_i$ such that
\begin{equation}\label{p2}
\p(t_i)\ge t_i^{-\lp-\ve}\qquad\text{for all }i\in\N \ .
\end{equation}

\noindent  The dimension part of Theorem \ref{thm2a} is obtained
by considering upper and lower bounds separately.

\medskip

\noindent{\bf The upper bound. \ } First notice that since $f$ is
continuously differentiable the map $x\mapsto(x,f(x))$ is locally
bi-Lipshitz and thus preserves Hausdorff dimension
\cite{falc,MAT}. Hence, we will investigate $\dim\Omega_{f,\p}$
instead of $\dim\cC_f\cap\cS_2(\p)$, where $\Omega_{f,\p}$ is
defined to be the set of $x\in I_0$ such that the system of
inequalities
\begin{equation}\label{e:032a}
\left\{
\begin{array}{l}
  \big|x-\frac{p_1}{q}\big|<\frac{\p(q)}{q},  \\[1ex]
  \big|f(x)-\frac{p_2}{q}\big|<\frac{\p(q)}{q}
\end{array}
\right.
\end{equation}
is satisfied for infinitely many $\bp/q\in\Q^2$. Furthermore,
there is  no loss of generality in assuming that $p_1/q\in I_0$
for solutions $\bp/q$ of (\ref{e:032a}).  

Next, without loss of generality, we can   assume that $I_0$ is
open in $\R$. Notice that the set $B:=\{ x \in I_0 : |f''(x)| = 0
\}$ is closed in $I_0$.  Thus the set  $G:= I_0\setminus B := \{ x
\in I_0 : |f''(x)| \neq 0 \}$ is open and a standard argument
allows one to write $G$ as a countable union of intervals $I_i$ on
which $f$ satisfies (\ref{infsup}) with $I_0$ replaced by $I_i$.
Of course, the constants $c_1$ and $c_2$ appearing in
(\ref{infsup}) will depend on the particular interval $I_i$. The
upper bound result will follow on showing that $\dim\Omega_{f,\p}
\cap I_i \leq d $, since by the conditions imposed on the theorem
$\dim B \leq d $ and so $$\dim\Omega_{f,\p} \leq \dim \left(B \cup
\mbox{ \small $\bigcup\limits_{i=1}^{\infty}$ } (\Omega_{f,\p}\cap
I_i ) \right) \ \leq \ d \ . $$

\noindent Without loss of generality, and  for the sake of clarity
we 
assume that $f$ satisfies (\ref{infsup}) on $I_0$.

For a point $\bp/q\in\Q^2$, denote by $\sigma(\bp/q)$ the set of
$x\in I_0$ satisfying (\ref{e:032a}). Trivially,
$|\sigma(\bp/q)|\le 2 \p(q)/q $.  Assume  that
$\sigma(\bp/q)\not=\emptyset$ and let $x \in \sigma(\bp/q)$. By
the mean value theorem, $f(x)=f(p_1/q)+f'(\tilde x)(x-p_1/q)$ for
some $\tilde x \in I_0$. We can assume that $f'$ is bounded on
$I_0$ since $f''$ is bounded and $I_0$ is a  bounded interval.
Suppose  $2^t\le q< 2^{t+1}$. By (\ref{e:032a}), $$\textstyle{
\big|f(\frac{p_1}{q})-\frac{p_2}{q}\big| \ \le \
\big|f(x)-\frac{p_2}{q}\big| \ + \ \big|f'(\tilde
x)\big(x-\frac{p_1}{q}\big)\big| \ \le \  c_4 \p(q)/q \ \le \
c_4\p(2^t)/2^t} $$ where  $c_4>0$ is a constant. In view of
(\ref{p1}), this implies that for any $\ve \in (0, 1) $ and $t$
sufficiently large
$$\textstyle{\big|f(\frac{p_1}{q})-\frac{p_2}{q}\big| \ \le \ 4
c_4  \ 2^{(t+1)(-\lp + \ve)} / 2^{t+1}  } \ . $$
 By (\ref{e:001a}), for $t$ sufficiently large
the number of $\bp/q\in\Q^2$ with $2^t\le q<2^{t+1}$ and
$\sigma(\bp/q)\not=\emptyset$ is at most $2^{t(2- \lp +3\ve)}$.
Therefore, with  $\eta  := (2-\lp+4\ve)/(\lp+1-\ve) $ it follows
that
\begin{eqnarray*}
\sum_{\bp/q\in\Q^2 \, :
\,\sigma(\bp/q)\not=\emptyset}|\sigma(\bp/q)|^\eta
 & = & \sum_{t=0}^\infty \ \ \ \ \sum_{\bp/q\in\Q^2,\,\sigma(\bp/q)
\not=\emptyset,\,2^t\le q<2^{t+1}}
\!\!\!\!\!\!\!\!\!\!\!\!\!\!\!\!\!\! |\sigma(\bp/q)|^\eta
\\ \\ & \le &  c' \,
\sum_{t=0}^\infty 2^{t(-\lp-1+\ve)\eta}\cdot2^{t(-\lp +2+3\ve)} \
= \ c'\sum_{t=0}^\infty 2^{- t\ve} \ < \ \infty \ \ ,
\end{eqnarray*} where $c'$ is a positive constant. By the
Hausdorff--Cantelli Lemma \cite[p.~68]{BD99},
$\dim\Omega_{f,\p}\le\eta $. As $\ve>0$ is arbitrary,
\begin{equation}\label{e:033a}
\dim\cC_f\cap\cS_2(\p) \ =  \ \dim\Omega_{f,\lp} \ \le \
d:=\frac{2-\lp}{\lp+1} \ \ .
\end{equation}

\noindent {\bf The lower bound (modulo Theorem \ref{thm2HM}). \ }
This is a simple consequence of Theorem \ref{thm2HM} and so all
that is required is that the curve is non-degenerate at a single
point.

\noindent Fix  $\epsilon>0$ such that $\lp+\epsilon<1$ and let $$
s \ := \ \frac{2-\lp-\epsilon}{1+\lp+\epsilon} \ < \ d \ \ . $$
Clearly, $s \in (1/2,1)$. In view of (\ref{p2}) and the fact that
$\psi$ is decreasing, there exists a strictly increasing sequence
$m_i$ of natural numbers  such that
\begin{equation}\label{ohgod}
 \psi(2^{m_i}) \ \geq \  2^{-(\lp+\epsilon)}  \
2^{-m_i(\lp+\epsilon)}  \ \  . \end{equation}
 To see that this is the case,
notice that for each $t_i$ there exists a natural number $m_i$
such that $2^{m_i} < t_i \leq 2^{m_i + 1 }$. It follows that $
\psi(2^{m_i}) \geq \psi(t_i) \geq t_i^{-(\lp+\epsilon)}  \geq
2^{-(m_i+1)(\lp+\epsilon)}$ and to ensure that $m_{i-1} < m_{i}$
simply  choose a suitable subsequence. By (\ref{ohgod}) and the
fact that $\psi$ is decreasing, we obtain that
\begin{eqnarray*}
\sum_{h=1}^{\infty}  h^{1-s} \p(h)^{s+1}   & =  &
\sum_{t=1}^\infty \ \ \sum_{2^{t-1}\le h<2^{t}} h^{1-s} \,
\p(h)^{s+1}  \  \gg  \
 \sum_{t=1}^\infty  2^{t(2-s)} \; \p(2^t)^{s+1} \\
 & \gg &  \sum_{i=1}^\infty  2^{m_i(2-s)} \  2^{-m_i(\lp+\epsilon)(s+1)}
 \ = \ \infty \ \ .
 \end{eqnarray*}
Hence,  Theorem \ref{thm2HM} implies that $ {\cal H}^s(
\cC_f\cap\cS_2(\p))  = \infty $ and so $ \dim\cC_f\cap\cS_2(\p)
\geq s $. As $\epsilon>0$ can be made arbitrarily small, we obtain
the required lower bound result.

\medskip

The Hausdorff measure part of Theorem \ref{thm2a} is a direct
consequence of Theorem \ref{thm2HM}.  Simply note that if  $
\limsup_{h \to \infty}  h^{2-d} \p(h)^{d+1} >  0 $ then $\sum
h^{1-d} \p(h)^{d+1} = \infty $ and also that if $ \lp \in (1/2,1)
$ then $d \in (1/2,1)$. The latter is obvious.  The former follows
by first observing that if $ \limsup_{h \to \infty}  h^{2-d}
\p(h)^{d+1} >  0 $, then there exists a strictly increasing
sequence $m_i$ of natural numbers such that $2^{m_i(2-d)} \;
\p(2^{m_i})^{d+1} \geq \eta >  0 $. It follows that
\begin{eqnarray*}
\sum_{h=1}^{\infty}  h^{1-d} \p(h)^{d+1}   & \gg   &
\sum_{t=1}^\infty  2^{t(2-d)} \ \p(2^t)^{d+1}
 \ \gg  \  \sum_{i=1}^\infty \   2^{m_i(2-d)} \p(2^{m_i})^{d+1}
 \ = \ \infty \ \ ,
 \end{eqnarray*}
 as required.
\hfill $\spadesuit$

\smallskip

Alternatively, the lower bound result for $
\dim\cC_f\cap\cS_2(\p)$ and  the Hausdorff measure part of Theorem
\ref{thm2a} can be deduced independently of Theorem \ref{thm2HM}
via Corollary~\ref{hdl}. Note that the upper bound result  is
complete. It  has been established without reference to any other
result.

 \mysection{Proof of Theorem~\ref{thm1}}

As $\cC:=\cC_f$ is non-degenerate almost everywhere, we can
restrict our attention  to a sufficiently small patch of $\cC$,
which can be written as $\{(x,f(x)):x\in I\}$ where  $I$ is a
sub-interval of $I_0$  and $f $ satisfies  (\ref{infsup}) with
$I_0$ replaced by $I$. Clearly, Theorem \ref{thm2a} is applicable
to $f$ restricted to $I$. However, without loss of generality and
for clarity, we assume that $f $ satisfies (\ref{infsup}) on
$I_0$.

Throughout this section, $\psi$ will be an approximating function
such  that
\begin{equation}\label{pdiv}
\sum_{h=1}^\infty\p^2(h) \ = \ \infty \ \ . \end{equation}


\noindent{\bf Step 1. \ } We show that there is no loss of
generality in assuming that
\begin{equation}\label{e:025}
\psi(h)\le h^{-1/2}\qquad\text{ for all $h$.}
\end{equation}
Define the auxiliary  function $\tilde\psi : h \to \tilde\psi(h)
:= \min\{h^{-1/2},\psi(h)\}$. Clearly $\tilde\psi$ is  an
approximating function. First we show that
\begin{equation}\label{e:026}
\sum_{h=1}^\infty\tilde\psi^2(h) \ = \ \infty \ .
\end{equation}
Assume  that (\ref{e:026}) is false. Then using the fact that
$\tilde\psi$ is decreasing, we obtain $$ 0\
\begin{array}[t]{c}\longleftarrow\\[-1.5ex]
\scriptstyle l\to\infty
\end{array}
\ \sum_{[l/2]\le h < l}\tilde\psi^2(h)\ge \sum_{[l/2]\le h <
l}\tilde\psi^2(l)\ge \tilde\psi^2(l)l/3 \ . $$ Thus,
$\tilde\psi(l)l^{1/2}\to0$ as $l\to\infty$. It follows that
$\tilde\psi(l)=o(l^{-1/2})$ and so  $\tilde\psi(l)=\psi(l)$ for
all but finitely many $l$. This together with (\ref{pdiv}) implies
(\ref{e:026}), a contradiction.

By definition,  $\cS_2(\tpsi)\subseteq\cS_2(\psi)$. Thus to
complete the proof of Theorem~\ref{thm1} it suffices to prove the
result with $\psi$ replaced by $\tpsi$. Hence, without loss of
generality, (\ref{e:025}) can be assumed.

\medskip

\noindent{\bf Step 2. \ } We show that there is no loss of
generality in assuming that
\begin{equation}\label{e:025b}
\psi(h)\ge h^{-2/3}\text{ for all $h$.}
\end{equation}
To this end, define $\hat\psi: h \to
\hat\psi(h):=\max\{\psi(h),h^{-2/3}\}$. It is readily verified
that $$ \cS_2(\hat\psi)=\cS_2(\psi)\cup\cS_2(h\mapsto h^{-2/3}).
$$ By the upper bound result
established in \S\ref{thm2aproof}, we have that $
\dim \cC_f \cap \cS_2(h\mapsto h^{-2/3}) \leq 4/5 < 1 $.  It
follows from the definition of Hausdorff dimension that
$\cH^1(\cC_f \cap \cS_2(h\mapsto h^{-2/3})) = 0$; i.e. for almost
all $x\in I_0$
 $$(x,f(x))\ \not\in\ \cS_2(h\mapsto
h^{-2/3}) \  \ . $$ Thus,  $$ \big|\{x\in
I_0:(x,f(x))\in\cS_2(\hat\psi)\}\big| \ = \ \big| \{x\in
I_0:(x,f(x))\in\cS_2(\psi)\}\big| $$ and to complete the proof of
Theorem~\ref{thm1} it suffices to prove that the set on the left
has full measure. Hence, without loss of generality,
(\ref{e:025b}) can be assumed.


\noindent{\bf Step 3. \   } In view of Steps 1 and 2 above, the
function $\psi$ satisfies (\ref{e:002}) and Corollary~\ref{cor1}
is applicable to $\p$. By (\ref{pdiv}) and the fact that $\psi$ is
decreasing,  we obtain that $$ \infty \ = \ \sum_{t=0}^\infty \ \
\sum_{2^t\le h<2^{t+1}}\p(h)^2 \ \le \ \sum_{t=0}^\infty \ \
\sum_{2^t\le h<2^{t+1}}\p(2^t)^2 \ = \
\sum_{t=0}^\infty2^t\p(2^t)^2 \ \ . $$ Hence
\begin{equation*}
\sum_{t=0}^\infty2^t\p(2^t)^2 \ = \ \infty.
\end{equation*}

\noindent Next, define the increasing function $u: \R^+ \to \R^+$
as follows $$ u(h) \ := \ \sum_{t=0}^{[h]} \; 2^t \; \p(2^t)^2 \ \
. $$ Trivially, $\lim_{t\to\infty}u(t)=\infty$. Let
$a_t=2^t\p(2^t)^2$ and $u_t=u(t)$. Fix $k\in\N$. Then $$
\sum_{t=k}^m\frac{a_t}{u_t}\ \ge \ \sum_{t=k}^m\frac{a_t}{u_m} \ =
\ \frac{u_m-u_{k-1}}{u_m} \to1\text{ \ as \  } m\to\infty. $$
Hence $$ \sum_{t=k}^\infty\frac{a_t}{u_t} \ \ge\ 1\ \text{ for all
}k. $$ This implies  that the sum $\sum_{t=1}^\infty a_t/u_t$
diverges; i.e.
\begin{equation}\label{e:027b}
\sum_{t=0}^\infty \ \frac{2^t \ \p(2^t)^2}{u(t)} \ = \ \infty \ \
.
\end{equation}

\noindent Now let $\Psi(t) = \Phi(t):=\p(t)/t$  and
$\rho(t):=u(\log_2t)/(t^2\p(t))$. By Corollary~\ref{cor1},
$(\Q^2_{\cC}(\Phi),\beta)$ is locally ubiquitous relative to
$\rho$,  where the function $\beta$ is given by (\ref{beta}). In
view of (\ref{e:027b}), $$
\sum_{t=1}^\infty\frac{\Psi(2^t)}{\rho(2^t)} \ := \
\sum_{t=1}^\infty\frac{\quad
\frac{\p(2^t)}{2^t}\quad}{\frac{u(t)}{2^{2t}\p(2^t)}} \ = \
\sum_{t=1}^\infty  \ \frac{2^t \ \p(2^t)^2}{u(t)} \ = \ \infty \ \
. $$

\noindent Since $\psi$  is decreasing,  $$ \Psi(2^{t+1})\ := \
\frac{\p(2^{t+1})}{2^{t+1}} \ \le \
\frac12\cdot\frac{\p(2^t)}{2^t} \ :=  \ \frac12 \ \Psi (2^t) \ \ .
$$

\noindent Thus the conditions of Lemma~\ref{ktl} are satisfied and
it follows that the set $\Lambda(\Q^2_{\cC}(\Phi),\beta,\Psi)$ has
full measure. By definition, the set
$\Lambda(\Q^2_{\cC}(\Phi),\beta,\Psi)$ consists of points  $x \in
I_0 $ such that the system of inequalities $$ \left\{
\begin{array}{l}
  \big|x-\frac{p_1}{q}\big|<\Psi(q)=
\frac{\p(q)}{q}<\frac{2\p(q)}{q}  \\[1ex]
\big|f(x)-\frac{p_2}{q}\big|<\Psi(q)+\Phi(q)=
\frac{\p(q)}{q}+\frac{\p(q)}{q}\le\frac{2\p(q)}{q}
\end{array}
\right. $$ is satisfied for infinitely many $\bp/q\in\Q^2$.
Obviously, for $x\in\Lambda(\Q^2_{\cC}(\Phi),\beta,\Psi)$ the
point $(x,f(x))$ is in $\cS_2(2\p)$. In order to complete the
proof of  Theorem~\ref{thm1}, simply apply what has already been
proved to the approximating function $\frac12\p$. \hfill
$\spadesuit$

\mysection{Proof of Theorem \ref{thm2HM} \label{thm2HMproof}}

We are assuming that there exists at least one point on the  curve
$\cC_f$ which is non-degenerate. Thus, there exists a sufficiently
small patch of $\cC_f$, which can be written as $\{(x,f(x)):x\in
I\}$ where  $I$ is a sub-interval of $I_0$ and $f $ satisfies
(\ref{infsup}) with $I_0$ replaced by $I$. Clearly, Theorems
\ref{thm1} and \ref{thm2a} are applicable to $f$ restricted to
$I$.  However, without loss of generality and for the sake of
clarity, we assume that  $f $ satisfies (\ref{infsup}) on $I_0$.

 Throughout this section, $s \in (1/2,1)$ and $\psi$ will
be an approximating function such  that
\begin{equation}\label{hmdiv}
\sum_{h=1}^{\infty} \ h^{1-s} \, \p(h)^{s+1} \ = \ \infty \ \ .
\end{equation}


\noindent{\bf Step 1. \ } We show that there is no loss of
generality in assuming that
\begin{equation}\label{e:095}
\lim_{t \to \infty} \psi(t) \ = \ 0 \ \ .
\end{equation}
Suppose on the contrary that $ \limsup_{t \to \infty} \psi(t) >
0$. Then for any $ s \leq 1 $, we have that (\ref{hmdiv}) holds.
In particular, $\sum_{h=1}^{\infty}  \p^2(h)  =  \infty $ and so
Theorem \ref{thm1} implies that $\cH^1(\cC_f \cap \cS_2(\p))
> 0 $. It follows that $\cH^s(\cC_f \cap \cS_2(\p)) = \infty
$ for any $s < 1$.   Hence, (\ref{e:095}) can be assumed.

\medskip

\noindent{\bf Step 2. \ } Since $s > 1/2$, there exists $\eta > 0
$ such that $s= \frac{1}{2}  +  \eta$.  We show that there is no
loss of generality in assuming that for all $h \in \N$,
\begin{equation}\label{e:095b}
\psi(h)\ge h^{-(1-\epsilon)} \hspace{6mm} \text{  where
}\hspace{6mm} 0 < \epsilon < 4 \eta/(3+2\eta)   \ \ .
\end{equation}
To this end, define $\hat\psi: h \to
\hat\psi(h):=\max\{\psi(h),h^{-(1-\epsilon)}\}$. It is readily
verified that $$ \cS_2(\hat\psi)=\cS_2(\psi)\cup\cS_2(h\mapsto
h^{-(1-\epsilon)}). $$ By the upper bound result
established in \S\ref{thm2aproof}, $ \dim \cC_f \cap
\cS_2(h\mapsto h^{-(1-\epsilon)}) \leq (1+\epsilon)/(2-\epsilon) <
s$ and so  $\cH^s(\cC_f \cap \cS_2(h\mapsto h^{-(1-\epsilon)})) =
0$.  Thus, $$ \cH^s(\cC_f \cap \cS_2(\hat\psi)) = \cH^s(\cC_f \cap
\cS_2(\psi))  $$ and to complete the proof of Theorem~\ref{thm2HM}
it suffices to prove that  $\cH^s(\cC_f \cap \cS_2(\hat\psi))=
\infty $. Hence, without loss of generality, (\ref{e:095b}) can be
assumed.


\noindent{\bf Step 3.  \ } In view of Steps 1 and 2 above, the
function $\psi$ satisfies (\ref{e:002}) and Corollary \ref{cor1}
is applicable to $\p$. In view of (\ref{hmdiv}), we can find a
strictly increasing sequence of positive integers
$\{h_i\}_{i\in\N}$ such that $$ \sum_{h_{i-1}<h \leq \ h_i} \
h^{1-s} \, \p(h)^{s+1} \
> \ 1 \  \hspace{6mm} (h_0 :=0) \  .
 $$ 
 Now simply  define the increasing function  $u$ as follows:
$$ u \, : \, h \, \to \, u(h) \, := \, i   \hspace{6mm} {\rm  for
} \ \ \ \  h_{i-1}<h\leq h_i \ \ . $$
Note  that $$ \sum_{h=1}^{\infty}  \ h^{1-s} \, \p(h)^{s+1} \,
u(h)^{-1}  \ = \ \sum_{i=1}^{\infty} \ \ \sum_{h_{i-1}<h \leq \
h_i} \ h^{1-s} \, \p(h)^{s+1} \, u(h)^{-1} \ > \
\sum_{i=1}^{\infty} i^{-1}   \ = \ \infty \ . $$ In particular,
since the function $\psi^{s+1}/u $ is decreasing we have that
$$\infty \  =  \ \sum_{t=0}^{\infty}  \sum_{2^t \leq h < 2^{t+1}}
\ h^{1-s} \, \p(h)^{s+1} \, u(h)^{-1} \ \leq \ 2^{2-s}
\sum_{t=0}^{\infty} \ 2^{t(2-s)} \ \p(2^t)^{1+s} \  u(2^t)^{-1} \
\ . $$ Hence
\begin{equation}\label{e:097b}
\sum_{t=0}^\infty \ 2^{t(2-s)} \ \p(2^t)^{1+s} \  u(2^t)^{-1} \ =
\ \infty \ \ .
\end{equation}

\noindent Now let $\Psi(t) = \Phi(t):=\p(t)/t$  and
$\rho(t):=u(t)/(t^2\p(t))$. By Corollary~\ref{cor1},
$(\Q^2_{\cC}(\Phi),\beta)$ is locally ubiquitous relative to
$\rho$,  where the function $\beta$ is given by (\ref{beta}). In
view of (\ref{e:097b}), $$
\sum_{t=1}^\infty\frac{\Psi(2^t)^s}{\rho(2^t)} \ := \
\sum_{t=1}^\infty \ 2^{t(2-s)} \ \p(2^t)^{1+s} \  u(2^t)^{-1} \ =
\ \infty \ \ . $$

\noindent Since $\psi$  is decreasing,  $ \Psi(2^{t+1})  \le
 \frac12 \ \Psi (2^t) $.
Thus the conditions of Lemma~\ref{JTl} are satisfied and it
follows that the set $\Lambda(\Q^2_{\cC}(\Phi),\beta,\Psi)$ is of
infinite $s$-dimensional Hausdorff measure. The statement of
Theorem~\ref{thm2HM} now follows on repeating verbatim the
argument given towards the end of the proof of Theorem \ref{thm1}.
\hfill $\spadesuit$

\mysection{Various Generalizations}

\subsection{Theorem \ref{thm2HM} for general Hausdorff measures
\label{GHMTHM}}

A {\em dimension function} $h \, : \, \R^+ \to \R^+ $ is an
increasing, continuous  function such that $h(r)\to 0$ as $r\to 0
\, $. Let ${\cal H}^{h}$ denote the  Hausdorff $h$--measure with
respect to the  dimension function $h$.  With reference to
\S\ref{HM},  this measure is defined by replacing $l^s(C_i)$ in
the definition of $s$--dimensional Hausdorff measure $\hs $ by the
quantity  $h(l(C_i))$ -- see \cite{falc,MAT} for further details.
In the case that $h: r \to r^s$ ($s \geq 0$), the measure ${\cal
H}^{h}$ is precisely $\hs $. For Hausdorff $h$--measures,
Jarn\'{\i}k's Theorem can be stated as follows -- see \cite[\S1.2
and \S12.1]{BDV03}.

\begin{thjarg}
 Let $h$ be a dimension function such that
$r^{-1} \, h(r)\to \infty$ as $r\to 0 \, $ and  $r^{-1} \, h(r) $
is decreasing.
Let $\p$ be  an approximating function.
Then $$ {\cal H}^{h}
\left(\cS_n(\p)\right)=\left\{\begin{array}{ll} 0 & {\rm if}
\;\;\; \sum_{r=1}^\infty \; r^n \,  h\left(\p(r)/r\right)  \;\;
 <\infty\\ &
\\ \infty & {\rm if} \;\;\; \sum_{r=1}^\infty \;
r^n \, h\left(\p(r)/r\right)   \;\;  =\infty
\end{array}\right..$$
\end{thjarg}

In the most simplest form, the following statement is the
Hausdorff $h$--measure analogue  of Theorem \ref{thm2HM}.

\begin{theorem} Let $f\in C^{(3)}(I_0)$, where $I_0$ is an interval and $\cC_f :
=\set{(x,f(x)):x\in I_0}$.  Assume that there exists at least one
point on the curve $\cC_f $ which is non-degenerate. Let  $\psi$
be an approximation function and let $h$ be a dimension function
such that $r^{-1} \, h(r) \to \infty$ as $r \to 0$, $r^{-1} \,
h(r)$ is decreasing and $r^{-(1/2 + \epsilon)} \, h(r) \to 0 $ as
$r \to 0$ for $\epsilon > 0 $ sufficiently small.
Furthermore, suppose $h$ satisfies the following growth condition:
there exist constants $r_0, \lambda_1,\lambda_2 \in (0,1) $ such
that $ h(\lambda_1 r) \leq \lambda_2 \, h(r) $ for $ r \in (0,
r_0)$. Then, $$
 {\cal H}^h \! \left((\cC_f \cap \cS_2(\p)) \right) \ =
\ \infty \hspace{10mm} {\rm if } \hspace{10mm} \sum_{r=1}^{\infty}
\; r \; \p(r) \; h(\p(r)/r)     \ = \ \infty \ \ . $$
\label{thm3**}
\end{theorem}

Apart from the growth condition imposed on the dimension function,
Theorem \ref{thm3**}  is the precise analogue of the divergent
part of Jarn\'{\i}k's General Theorem  for simultaneous
Diophantine approximation on planar curves. The growth condition
is not particularly restrictive and can be completely removed from
the statement of the theorem in the case that $ G:= \limsup_{r \to
\infty} \; h(\p(r)/r) \; \p(r) \, r^{2} > 0 $. Furthermore, when
$G=0$, if there exists a constant $\lambda \in (0,1)$ such that
$\p(2r) > \lambda \p(r) $ for all sufficiently large $r$ then the
growth condition on $h$ is again redundant.

\medskip

 Notice that if $h: r \to
r^s$ ($s \geq 0$), then the growth condition is trivially
satisfied and the above theorem reduces to  Theorem \ref{thm2HM}.

\medskip

\noindent{\bf  Remark on the proof of Theorem \ref{thm3**}. \ }
The first step is to obtain the analogue of Lemma \ref{JTl} for
general Hausdorff measures. This is easy, it follows directly from
Corollary 3 of \cite[\S5]{BDV03} in the same way that Lemma
\ref{JTl} is deduced from  Lemma \ref{JT}. The proof of Theorem
\ref{thm3**} then  follows on  modifying the  argument used to
prove Theorem \ref{thm2HM} in \S\ref{thm2HMproof}. Note that
Corollary \ref{cor1}, the important local ubiquity statement which
gives the `optimal' ubiquitous function $\rho$,  is independent of
any dimension function. The following useful fact concerning
dimension functions is also required: if  $ \, f$ and $g$ are two
dimension functions such that $f(r)/g(r) \to 0 $ as $ r \to 0 $,
then ${\cal H}^{f} (F) =0 $ whenever ${\cal H}^{g} (F) < \infty $.
We leave the details to the reader.

\subsection{The Multiplicative Problems/Theory
\label{Multi}}

Given an approximating function $\psi$, a point $\vv
y=(y_1,\dots,y_n)\in\R^n$ is called {\it simultaneously
multiplicatively $\psi$--approximable} if there are infinitely
many $q\in\N$ such that $$ \prod_{1\le i\le n}\|q y_i\| \ < \
\psi(q)^n . $$ Thus, the maximum  in the definition of
simultaneously $\psi$--approximable is replaced by the product.
Denote by $\cS_n^{\rm M} (\psi)$ the set of simultaneously
multiplicatively $\p$--approximable points.  Trivially, we have
that $$ \cS_n(\psi) \ \subset  \  \cS_n^{\rm M} (\psi)  \ . $$

 The two
fundamental problems posed in the introduction can obviously be
reinstated for the multiplicative setup. In a forthcoming paper
\cite{BVe}, the first and third authors develop the simultaneous
multiplicative theory for metric Diophantine approximation on
planar curves. As an illustration of the type of results
established in \cite{BVe}, we mention the following analogue  of
Theorem \ref{thm2a}. With the same notation and hypotheses of
Theorem \ref{thm2a},  $$ \dim \cC_f \cap \cS_2^{\rm M}(\p) \ = \
\frac{2-\lp}{\lp+1} \ \ .  $$



\newpage
\begin{appendix}

\begin{center}
{\Large \bf Appendix I:  Proof of Ubiquity Lemmas}
 \addcontentsline{toc}{section}{{\bf APPENDIX ~I:~~~  Proof of Ubiquity Lemmas}}
\end{center}



\vspace{-5mm}


\myappsection{Ubiquity with respect to sequences  \label{pf12}}

In this appendix we prove the  ubiquity lemmas  of \S\ref{ubin1}
which are the key towards establishing the divergent results of this
paper. It is both more convenient and no more difficult to consider
a slightly  more general setup in which the sequence $\{2^n\}$ is
replaced by an arbitrary increasing sequence $u$. Apart from this
the setup remains unchanged.

Let $u:=\{u_n\}$  be a  positive increasing sequence such that $
 \lim_{n \to
\infty} u_n = \infty $ and  let $ J^u(n):=\{\alpha\in J \, : \,
\beta_{\alpha} \le u_n\} $. Given a function $\rho : \R^+ \to \R^+
$  such that  $\lim_{t \to \infty} \rho(t) = 0$, let
$$\Delta^u(\r,n) := \bigcup_{\alpha\in J^u(n)} B(\ra,\r(u_n)) \
.$$

\begin{definition}
\label{uub}  Suppose there exists a function $\r$, a sequence $u$
and an absolute constant $\ka>0$ such that for any interval
$I\subseteq I_0$
\begin{equation}\label{xub}
  | \Delta^u(\r,n) \cap I | \ \ge \ \ka \ |I|  \qquad \mbox{for $n \geq n_o(I)$} .
\end{equation}
Then the pair $(\cR,\beta)$ is said to be {\em locally ubiquitous
in $I_0$ relative to $(\r,u)$}.
\end{definition}

 Notice that any subsequence $v$ of $u$  will also do in the above definition;
i.e. (\ref{xub})  is satisfied for $\Delta^v(\r,n)$. In order to
state  the consequences of this slightly more general definition
of ubiquity  we introduce the following notion. Given a sequence
$u$, a function $h$ will be said to be {\bf ${\bf u}$-regular} if
there exists a strictly positive constant $\lambda < 1$ such that
for $n$ sufficiently large
\begin{equation} h(u_{n+1}) \leq \lambda \, h(u_n) \ .  \label{afmh}
\end{equation}
The constant $\lambda$ is independent of $n$ but may depend on
$u$. Clearly,  if $h$ is $u$-regular then it is $v$--regular for
any subsequence $v$ of $u$.

\begin{theorem} \label{xthm1}
Suppose that $(\cR,\beta)$ is locally ubiquitous in $I_0$ relative
to $(\rho,u)$ and let  $ \Psi$ be  an approximating function such
that $\Psi$ is $u$-regular. Then $$ |\Lambda(\cR,\beta,\Psi)| \ =
\ \mbox{\rm F{\scriptsize ULL}}
 \ := \ |I_0| \hspace{10mm} {\rm if } \hspace{10mm}
\sum_{n=1}^\infty\frac{\Psi(u_n)}{\rho(u_n)}=\infty \ \ .$$
\end{theorem}

\begin{theorem}\label{xthm3}
  Suppose that $(\cR,\beta)$ is
locally ubiquitous in $I_0$ relative to $(\rho,u)$ and let $\Psi$
be an approximation function. Let $ s \in (0,1) $ and let
\begin{equation}
 G \, := \, \limsup_{t \to \infty} \,
\frac{\Psi(u_n)^s}{\rho(u_n)} \; \ . \label{my5}
\end{equation}
\noindent {\bf (i)}  Suppose that $G = 0$  and  that $\Psi$ is
$u$-regular. Then, $$
 {\cal H}^s \! \left( \Lambda(\cR,\beta,\Psi) \right) \ =
\ \infty \hspace{10mm} {\rm if } \hspace{10mm} \sum_{n=1}^{\infty}
\frac{\Psi(u_n)^s}{\rho(u_n)}\ = \ \infty \ \ . $$

\noindent {\bf (ii)} Suppose that $G >0 $. Then,  ${\cal H}^s \!
\left( \Lambda(\cR,\beta,\Psi) \right) \ = \ \infty$.

\end{theorem}

In the case $ u= \{u_n\}  := \{2^n\} $,  these theorems clearly
reduce to Lemmas \ref{KT} and \ref{JT} of \S\ref{ubin1}.

\myappsubsection{Prerequisites}

\myappsubsubsection{The Mass Distribution Principle and a covering
lemma \label{hmd} }

 A general and classical  method for obtaining a
lower bound for the $s$--dimensional Hausdorff measure of an
arbitrary set $F$ is the following mass distribution principle.

\begin{lemma}[Mass Distribution Principle]
Let $\mu$ be a probability measure supported on a subset $F$ of
$\R$. Suppose there are  positive constants $c$ and $r_o$ such
that $ \mu(B)  \le   c \,   r^s  $ for any ball  $B$ with radius
$r \leq r_o$.  Then,  $    {\cal H}^{s} (F) \geq 1/c  $.
\end{lemma}

 The following  covering result will be used at various
stages during  the proof of our theorems.
\begin{lemma}[Covering lemma]
Let  ${\cal B}$ be a finite collection of balls  in $\R$ with
common radius $r>0$. Then there exists a disjoint sub-collection
$\{B_i\}$ such that $$ \bigcup_{B\in {\cal B}} \, B \ \subset \
\bigcup_i \  {3B}_i \ . $$
\end{lemma}

These lemmas are easily established and relatively standard, see
\cite{falc,MAT} and \cite[\S7]{BDV03}.

\myappsubsubsection{Positive and full measure sets}

\begin{proposition}\label{lem1b}
Let $E \subset \R$ be a measurable set and let $I_0 \subset \R$ be
an interval. Assume that there is a constant $c>0$ such that for
any finite interval $I \subset I_0$  we have that $|E\cap I|\ge c\
|I|  $.  Then $E$ has full measure in $I_0$, i.e.\ $|I_0\setminus
E|=0$.
\end{proposition}

For the proof see \cite[Lemma 2]{vb1} and \cite[\S8]{BDV03}.

\begin{proposition}\label{lem2}
Let  $E_n \in \R$ be a sequence of measurable sets such that
$\cup_{n=1}^{\infty} E_n $ is bounded and $\sum_{n=1}^\infty
|E_n|=\infty $. Then $$|\limsup_{n \to \infty} E_n | \; \geq \;
\limsup_{Q \to \infty} \frac{ \left( \sum_{s=1}^{Q} |E_s|
\right)^2 }{ \sum_{s, t = 1}^{Q} |E_s \cap E_t | }  \  \  \ . $$
\end{proposition}

This result is a generalization of the divergent part of the
standard Borel--Cantelli lemma. For the proof see Lemma 5 in
\cite{Spr79}.

\begin{proposition}\label{lem0}
Let $E_n \in \R$ be a  sequence of measurable sets and let  $I
\subset \R$ be a bounded interval. Suppose there  exists a
constant $c> 0$ such that $\limsup_{n \to \infty} |I \, \cap \,
E_n| \ge c \; |I|$. Then, $ |I \, \cap \, \limsup_{n\to\infty}E_n
| \ \ge \ c^2 \, |I| $.
\end{proposition}

\medskip

\noindent{\bf  Proof of Proposition \ref{lem0}.} For any $0< \ve <
c$, there is a subsequence $E_{n_i}$ with $n_i$ strictly
increasing such that $|I \cap E_{n_i})\ge (c-\ve) \: |I|$. Clearly
$\textstyle \big(\sum_{i=1}^N |I \cap E_{n_i}|\,\big)^2\ge
\big(\sum_{i=1}^N (c-\ve) \; |I| \,\big)^2= (c-\ve)^2N^2 \, |I|^2
$ and $\textstyle{ \sum_{n,m=1}^N |I \cap E_n\cap E_m|\le
\sum_{m,n=1}^N |I|= |I| \, N^2 } . $ Also notice  that
$\sum_{i=1}^\infty |I \cap E_{n_i}| \geq |I| \,
\sum_{i=1}^\infty(c-\ve) =\infty$. Thus on applying Proposition
\ref{lem2} and observing that $ I \, \cap \,
\limsup_{n\to\infty}E_n   \supseteq     I \, \cap  \,
\limsup_{i\to\infty}E_{n_i} $ we have that $$ \Big|I \cap
\limsup_{n\to\infty}E_n\Big|\ge
\limsup_{N\to\infty}\frac{(c-\ve)^2N^2 |I|^2}{|I|N^2} =(c-\ve)^2\,
|I| \ . $$ As $\ve>0$ is arbitrary,  this completes the proof of
the proposition. \hfill $\spadesuit$


\myappsection{ Proof of Theorem \ref{xthm1} \label{proof1} }

Let $B$ be an arbitrary ball
 in $I_0$ and let $r(B)$ denote its radius.
 In view of Proposition~\ref{lem1b},  the aim is to show that
\begin{equation}\label{my1}
|\Lambda(\cR,\beta,\Psi)\cap B| \ \ge \  |B|/C \ ,
\end{equation}
where  $C>0$  is a constant independent of $B$.

\myappsubsection{The subset $A(\Psi,B)$ of
$\Lambda(\cR,\beta,\Psi)\cap B $\label{proof1.1}}

Consider the  collection  of balls  $ 
\{B(\ra,2\r(u_n)):\alpha \in J^u(n) \}$.    By the covering lemma,
there exists a disjoint sub-collection
$
\{B(\ra,2\r(u_n)):\alpha \in G^u(n) \}$, where $G^u(n)$ is a
subset of  $J^u(n)$, 
 such that
\begin{equation}\label{my2}
 \bigcup_{\alpha \, \in
G^u(n)}^\circ B(\ra,\r(u_n)) \ \subset \ \Delta^u(\rho,n) \
\subset \ \bigcup_{\alpha \, \in  G^u(n)} B(\ra,6\r(u_n)) \ .
\end{equation}
The left hand side follows from the fact that the balls $B(\ra, 2
\r(u_n))$ with $\alpha\in G^u(n) $ are disjoint.  Choose $n$
sufficiently large so that $24\r(u_n) < r(B) $ (by definition,  $
\r(u_n) \to  0$ as $n \to \infty$) and let $$ \gb(n) \, := \,
\left\{ \alpha \in G^u(n) \ : \ \ra \in \mbox{\small
$\frac{1}{2}$} \,  B \right\} \ \ . $$ In view of  (\ref{my2}),
\begin{equation}
 \bigcup_{\alpha \, \in \gb(n)}^\circ
B(\ra,\r(u_n)) \ \subset \ \Delta^u(\rho,n) \ \cap \ B \label{r<p}
\end{equation}  and
$$ \bigcup_{\alpha \, \in \gb(n)} B(\ra,6\r(u_n))) \ \supset \
\Delta^u(\rho,n) \ \cap \ \mbox{\small $\frac{1}{4}$} \, B  \ \ .
$$ We now estimate the cardinality of $\gb(n)$. By (\ref{xub}),
for $n$ sufficiently large
\begin{eqnarray*}
\#\gb(n)  \ \r(u_n) \   \gg \  | \textstyle{\bigcup_{\alpha\in
\gb(n)} B(\ra,6\r(u_n)) }|  \ \geq
 \ |\Delta^u(\rho,n)  \cap  \mbox{\small $\frac{1}{4}$} \,
B | \ \
   \gg \
\ka \, |B|
 \ .
\end{eqnarray*}
  On the other hand, $ |B|
 \ge   \big|\textstyle{\bigcup\limits^{\circ}}_{\alpha\in
\gb(n)} B(\ra,\r(u_n))\big|  \gg  \#\gb(n) \,  \r(u_n)  $. The
upshot  is that
 \begin{equation}\label{NGn2}
 \#\gb(n) \ \asymp \
 \frac{|B|}{ \r(u_n)} \ \ .
\end{equation}

\medskip

\noindent Suppose for the moment that for some sufficiently large
$n \in \N$ we have that $\Psi(u_n) \geq  \r(u_n)$.  Then
(\ref{xub})  implies that $ |\Delta^u(\Psi,n) \cap B | \, \geq \,
|\Delta^u(\rho,n) \cap B | \, \geq \, \ka \, |B| $.  Thus, if
$\Psi(u_n) \geq  \, \r(u_n)$ for infinitely many $n \in \N$,
Proposition \ref{lem0}  implies (\ref{my1}) and we are done.
Hence, without loss of generality,  we can assume  that for $n$
sufficiently large
\begin{equation}\label{my3}
\r(u_n) >   \; \Psi(u_n)\; .
\end{equation}
Now let $$ A_n(\Psi,B) \, := \, \bigcup_{\alpha\in \gb(n)}^{\circ}
B(\ra,\Psi(u_n)) \ \ .$$ The disjointness is a consequence
(\ref{my3}). Indeed, for $\alpha \in \gb(n)$ the balls $B(\ra,
2\r(u_n)) $ are disjoint. Therefore, $|A_n(\Psi,B)| \asymp
\Psi(u_n) \  \# \gb(n)$ and in view of (\ref{NGn2}) we have that
\begin{equation}\label{mAncomp2}
|A_n(\Psi,B)|  \ \asymp \
 |B| \ \times \  \frac{\Psi(u_n)}{\r(u_n)} \ .
\end{equation}
Finally, let $$ A(\Psi,B) \, := \, \limsup_{n \to \infty} \,
A_n(\Psi,B) \, := \, \bigcap_{m=1}^{\infty} \
\bigcup_{n=m}^{\infty} A_n(\Psi,B ) \, . $$ By construction,  we
have $ A_n(\Psi,B)  \subset \Delta^u(\Psi,n) \cap B  $ and  it
follows that  $A(\Psi,B)\setminus{\cal R} $ is a subset of
$\Lambda(\cR,\beta,\Psi) \cap B$. Now in view of (\ref{my1}) and
the fact that ${\cal R} $ is countable and therefore of measure
zero, the proof of Theorem~\ref{xthm1} will be completed on
showing that
\begin{equation}\label{my1+}
|A(\Psi,B)| = |A(\Psi,B) \cap B  |\ \ge \  m(B)/C.
\end{equation}
Notice that  (\ref{mAncomp2})  together with the divergent sum
hypothesis  of the theorem implies that
\begin{equation}
\sum_{n = 1}^{\infty} |A_n(\Psi,B)| = \infty \ . \label{div}
\end{equation}

\noindent In view of  Proposition \ref{lem2}, this together with
the following quasi-independence on average result implies
(\ref{my1+}) and thereby  completes the proof of Theorem
\ref{xthm1}.

\begin{lemma}[Quasi--independence on average]
\label{In2}
 $ \ \; $ There exists a constant $C >1$ such that for
$Q$ sufficiently large, $$  \sum_{s, t = 1}^{Q} |A_s(\Psi,B) \cap
A_t(\Psi,B) |  \; \leq \; \frac{C}{|B|} \,  \left( \sum_{s=1}^{Q}
|A_s(\Psi,B)| \right)^2  \; . $$
\end{lemma}

\noindent{\bf Proof of Lemma  \ref{In2}. } Throughout,  write
$A_t(\Psi)$ for $ A_t(\Psi,B)$.  Also, let $s<t$ and  note that
\begin{eqnarray}
|A_s(\Psi)\cap A_t(\Psi)| & = & \sum_{\alpha\in G_{\mbox{\tiny
B}}^u(s)}
  |\, B(\ra,\Psi(u_s))\cap
 A_t(\Psi)\, |  \ .
  \label{asintat}
  \end{eqnarray}
Let $B_s(\Psi)$ denote a generic ball $B(\ra,\Psi(u_s))$ with
$\alpha\in G_{\mbox{\tiny B}}^u(s)$.  We now obtain an upper bound
for $|B_s(\Psi)\cap A_t(\Psi)|$. Trivially,
\begin{eqnarray}
|B_s(\Psi)\cap A_t(\Psi)|  & := & |B_s(\Psi) \, \cap \,
\textstyle{\bigcup\limits^{\circ}}_{\alpha\in G^u_{\mbox{\tiny
B}}(t)}  \ B(\ra,\Psi(u_t)) |  \nonumber \\ & & \nonumber \\  & =
& \!\!\! \sum_{\alpha \in G_{\mbox{\tiny B}}^u(t)}
 |B_s(\Psi)\cap B(\ra,\Psi(u_t))| \  \ll \ N(t,s) \ \Psi(u_t) \
 \label{mbsintat}
\end{eqnarray}
where $N(t,s) := \# \{ \alpha \in G_{\mbox{\tiny B}}^u(t) :
B_s(\Psi)\cap B(\ra,\Psi(u_t)) \neq \emptyset \} $. We proceed by
considering  two  cases. 

\bigskip

\noindent{\bf Case (i):} $t > s $ such that $ \Psi(u_s)<\r(u_t)$.
\ \ \ On using the fact that  the balls $B(\ra,2\r(u_t))$ with
$\alpha \in G_{\mbox{\tiny B}}^u(t)$  are disjoint, it is easily
verified that $ N(t,s) \leq 1 $.  This together with (\ref{NGn2}),
(\ref{asintat})  and (\ref{mbsintat}) implies that
\begin{eqnarray*}
|A_s(\Psi)\cap A_t(\Psi)| \ \ll \ \# G_{\mbox{\tiny B}}^u(s) \
\Psi(u_t) \ \ll \  |B| \ \times \ \frac{\Psi(u_t)}{\r(u_s)}
 \ .
\end{eqnarray*}

\noindent {\bf Case (ii):} $t > s $ such that $
\Psi(u_s)\ge\r(u_t)$. \ \ \ First note that if $B_s(\Psi)\cap
B(\ra,\r(u_t)) \neq \emptyset $, then  $ B(\ra,\Psi(u_t)) \subset
3 B_s(\Psi) $.  The  balls $B(\ra,\r(u_t))$ with $\alpha \in
G_{\mbox{\tiny B}}^u(t)$ are disjoint and so   $ N(t,s) \ll
\Psi(u_s)/ \r(u_t) $. It now follows, via (\ref{NGn2}),
(\ref{mAncomp2}), (\ref{asintat}) and (\ref{mbsintat}), that $$
|A_s(\Psi)\cap A_t(\Psi)| \  \ll \ \frac{1}{|B|} \ |A_s(\Psi)| \
|A_t(\Psi)| \ . $$

The upshot of these two cases, is that for $Q$ sufficiently large
$$
\begin{array}{l}
 \displaystyle
 \sum_{s, t = 1}^{Q} |A_s(\Psi) \cap A_t(\Psi) |   \displaystyle
 = \
\sum_{s=1}^{Q} |A_s(\Psi)| \ + \ 2\,  \sum_{s=1}^{Q-1} \
 \sum_{\substack{s+1 \leq t \leq Q \\ { \rm case (i)} } }
 |A_s(\Psi) \cap A_t(\Psi)|
\\[2ex] \hspace*{29ex} \displaystyle + \  2 \, \sum_{s=1}^{Q-1} \
\sum_{\substack{s+1 \leq t \leq Q \\ { \rm case (ii)} } }
|A_s(\Psi) \cap A_t(\Psi)| \\[2ex] \hspace*{21ex} \displaystyle
 \ll  \ \sum_{s=1}^{Q}
|A_s(\Psi)| \quad  +  \quad \frac{1}{|B|}\Big( \textstyle{
\sum_{s=1}^{Q} } |A_s(\Psi)| \Big)^2 \\[2ex] \hspace*{29ex}
\displaystyle + \ \
 |B|\ \sum\limits_{s=1}^{Q-1}  \sum_{\substack{s+1 \leq t \leq Q
 \\
\Psi(u_s)<\r(u_t) } } \frac{\Psi(u_t)}{\r(u_s)}
  \ .
\end{array}
$$ We now make use of the fact that $\Psi$ is $u$-regular. For $t
> s$ with $s$ sufficiently large,  we have
that $ \Psi(u_t) \leq \lambda^{t-s}   \Psi(u_s) $ for some
$0<\lambda < 1$.  This together with (\ref{mAncomp2}), implies
that for $Q$ sufficiently large
 $$ |B| \ \sum\limits_{s=1}^{Q-1}  \sum_{\substack{s+1 \leq t \leq Q
 \\
\Psi(u_s)<\r(u_t) } } \frac{\Psi(u_t)}{\r(u_s)}  \ll \, |B|
\sum_{s=1}^{Q-1} \ \frac{\Psi(u_s)}{\r(u_s)} \ \sum_{s \, < \, t
\leq \, Q}   \lambda^{t-s} \ll \ \sum_{s=1}^{Q} \ |A_s(\Psi)|  \ .
 $$
By (\ref{div}), for $Q$ sufficiently large $\sum_{s=1}^{Q}
|A_s(\Psi)|\le |B|^{-1}(\sum_{s=1}^{Q} |A_s(\Psi)|)^2$.  The
statement of Lemma \ref{In2} now readily follows. This completes
the proof of Theorem \ref{thm1}. \hfill $\spadesuit$

\myappsection{Proof  of Theorem \ref{xthm3} \label{proof3}}

To prove Theorem \ref{xthm3} we proceed as follows. For any fixed
$\eta \gg 1$ we construct a Cantor subset ${\bf K_{\eta}} $ of
$\Lambda(\cR, \beta, \Psi)$ and a probability measure $\mu$
supported on ${\bf K_{\eta}}$ satisfying the condition that for an
arbitrary ball $A$ of sufficiently small radius $r(A) $
\begin{equation}
\mu(A) \; \ll \; \frac{r(A)^s}{ \eta } \;   , \label{task}
\end{equation}
where the implied constant is absolute. By the Mass Distribution
Principle, the above inequality implies that $ {\cal H}^s({\bf
K_{\eta}}) \gg \eta $. Since $ {\bf K_{\eta}}\subset
\Lambda(\cR,\beta, \Psi)$, we obtain that ${\cal H}^s\left(
\Lambda(\cR,\beta, \Psi) \right)\gg \eta$. However, $\eta \gg 1  $
is arbitrarily large whence ${\cal H}^s\left( \Lambda(\cR,\beta,
\Psi) \right)=\infty$ and this proves Theorem \ref{xthm3}.

In view of the above outline, the whole strategy of our proof is
centred around the  construction of a `right type' of Cantor set
${\bf K_{\eta}}$ which  supports a measure $\mu$ with the desired
property. The actual nature of the construction of ${\bf
K_{\eta}}$  depends heavily on whether  $G$ defined by (\ref{my5})
is finite  or infinite. We first    deal with the case that $0\leq
G < \infty$. The case that $G = \infty$ is substantially easier.

\myappsubsection{Preliminaries}
In this section we group together for clarity and convenience
various concepts and results which will be required in
constructing the Cantor set ${\bf K_{\eta}}$. Throughout,  $g$
will denote the function given by $$ 
g(r) := \Psi(r)^s    \r(r)^{ -1}   \hspace{6mm}  \mbox{and so } \
\ \  G \, := \, \limsup_{n \to \infty} \, g(u_n) \ .  $$
\myappsubsubsection{{\bf The sets $\gb(n)$ again.\label{pregb}}}
Let $B$ be an arbitrary ball in $I_0$ with radius $r(B)$.  Relabel
the sets $\gb(n)$
 constructed in \S\ref{proof1}  by $\tilde{G}_{\mbox{\!\tiny \em B}}^{u}(n)$.
 By keeping track of constants, the
estimate (\ref{NGn2}) for $\#\tilde{G}_{\mbox{\!\tiny \em
B}}^{u}(n)$ is explicitly as follows: $$
 \frac{ \ka}{24} \
\frac{r(B)}{\r(u_n)}  \ \leq \  \# \tilde{G}_{\mbox{\!\tiny \em
B}}^{u}(n) \ \leq \ \frac{r(B)}{\r(u_n)}
 \ ,  $$
where $\ka$ is as in (\ref{xub}). Now let $ 0 \, < \, c_1:= \frac{
\ka}{24} \, < 1 $
 and  define $\gb(n)$ to be any
sub-collection of $ \tilde{G}_{\mbox{\!\tiny \em B}}^{u}(n) $ such
that
$$ \#\gb(n) \ = \left[c_1 \ \frac{r(B)}{\r(u_n)} \right] \ ,  $$
where $[x]$ denotes the integer part of a real number $x$.
 Thus, for $n$ sufficiently large
\begin{equation}
\mbox{\small $ \frac{1}{2} $ }  \, c_1 \ \frac{r(B)}{\r(u_n)} \
\leq \ \#\gb(n) \ \leq  \ c_1 \ \frac{r(B)}{\r(u_n)}
 \ .  \label{NGnB}
\end{equation}

\noindent{\bf Remark. } Recall, that by construction the balls
$B(\ra,2\r(u_n))$ with $\alpha \in G_{\mbox{\tiny B}}^u(n)$  are
disjoint. Also note, that  we
 can assume  that
$\rho(u_n)^{-1} \, \Psi(u_n) \to 0$ as $n \to \infty $. If this
was not the case then $\limsup \rho(u_n)^{-1} \, \Psi(u_n) > 0$ as
$n \to \infty $ and Theorem \ref{thm1} implies that
$|\Lambda(\cR,\beta,\Psi)| = |I_0| > 0$. In turn, ${\cal H}^{s}
(\Lambda(\cR,\beta,\Psi)) = \infty $ for any $s \in [0,1)$ and we
are done. Hence, without loss of generality, we can assume that
for $n$ sufficiently large
\begin{equation} 2 \, \Psi(u_n) \  <  \  \rho(u_n)  \hspace{9mm} {\rm and }
\hspace{9mm} \lim_{n \to \infty }\Psi(u_n)= 0    \  \ .
\label{mythm33}
\end{equation}

\myappsubsubsection{{\bf Working on a  subsequence of $u$ and the
ubiquity function $\rho$.\label{subseqwork}}} The proof of Theorem
\ref{xthm3} in the case that  $G$ is finite  relies on the fact that
the ubiquity function $\rho$ can be taken to be $u$-regular with
constant $\lambda$ as small as we please. The fact that we have
assumed that the approximating function $\Psi$ is $u$-regular in the
hypothesis of the theorem is purely for convenience with the
application to planar curves in mind. To begin with recall the
following simple facts: (i) if we have local ubiquity for a
particular sequence $u$ then we automatically have local ubiquity
for any subsequence $v$ and (ii) if a function $h$ is $u$-regular
then it is $v$-regular for any subsequence $v$. Also note that if
$G$ is finite, then $\limsup_{n \to \infty}g(v_n) < \infty $ for any
subsequence $v$ of $u$.

Suppose $G$ is finite and fix some   $\lambda  \in (0,1)$.  We now
prove the existence of an appropriate subsequence $v$ of $u$ on
which $\r$ is $v$-regular with constant $\lambda$ and $\sum g(v_n)
= \infty $. In the case  $G=0$ (part (i)  of Theorem \ref{xthm3}),
we have that  $\Psi$ is $u$-regular and so there exists a constant
$\lambda_* \in (0,1) $ such that $ \Psi(u_{n+1}) \leq \lambda_* \,
\Psi(u_n) $ for all $n$ sufficiently large. It follows that for $n
$ sufficiently large $$ x_{n+1} := \Psi(u_{n+1})^s \leq
\lambda_*^s \Psi(u_{n})^s = \lambda_*^s \,  x_{n} \ . $$ Next, fix
some sufficiently large $n_1$ and for $k \geq 2 $ let $n_k$ be the
least integer strictly greater than $n_{k-1}$ such that $
\rho(u_{n_k}) \leq \lambda \, \rho(u_{n_{k-1}}) $. This is
possible since $\r(r) \to 0$ as $ r \to \infty $. By construction,
$\rho(u_{m}) \geq  \lambda \,  \rho(u_{n_{k-1}}) $ for any integer
$ m \in [n_{k-1}, n_{k}-1 ] $. It follows that
\begin{eqnarray*}
\infty & = & \ \sum_{n=n_1}^{\infty} g(u_n) \  =  \
\sum_{n=n_1}^{\infty} x_n \; \rho(u_n)^{-1}   =
\sum_{k=2}^{\infty} \ \ \sum_{n_{k-1} \leq m < n_k } \!\!\!\! x_m
\; \rho(u_m)^{-1} \\ & \leq & \sum_{k=2}^{\infty} \sum_{n_{k-1}
\leq m < n_k } \!\!\!\! x_m \; \rho(u_{n_{k-1}})^{-1} \,
\lambda^{-1}  = \lambda^{-1} \sum_{k=2}^{\infty}
\rho(u_{n_{k-1}})^{-1} \!\!\!\! \sum_{n_{k-1} \leq m < n_k }
\!\!\!\!\!\!\!\! x_m \\ & \ll & \sum_{k=2}^{\infty}
\rho(u_{n_{k-1}})^{-1} x_{n_{k-1}} \ \sum_{i=0}^{\infty}
(\lambda_*^s)^i
\ll \sum_{k=1}^{\infty} \rho(u_{n_{k}})^{-1}  x_{n_{k }}
 :=  \sum_{k=1}^{\infty} g(u_{n_{k}}) .
\end{eqnarray*}
Now set $v:= \{u_{n_k} \}$. By construction, $\rho$ is $v$-regular
with constant $\lambda$ and $\sum g(v_n)= \infty $. Next notice
that if $0 < G < \infty $, then there exists a strictly increasing
sequence $\{ n_i\} $ such that  $g(u_{n_i}) \geq G/2
> 0 $. Since $\lim_{r \to \infty}  \r(r) = 0$,  it follows
that for any $\lambda < 1 $ there exists a subsequence $v$  of $\{
u_{n_i} \} $ such that $\r(v_{t+1})   <   \lambda \, \r(v_{t}) $
and $  \sum g(v_t) = \infty $. The upshot is that  in establishing
Theorem \ref{xthm3} for the case that $0 \leq G < \infty $ we can
assume that $\r$ is $u$ regular with constant $\lambda$ as small
as we please.


\myappsubsection{Proof  of Theorem \ref{xthm3}: $ 0\leq G < \infty$
}

\myappsubsubsection{The Cantor set ${\bf K_{\eta}}$ }

 Let $G^*:= \max \{2,
\sup_{n \in \N} g(u_n) \} $  and fix a real number $\eta >  G^*$.
Thus $$ g(u_n) \ < \ G^*  \ < \ \eta \hspace{8mm} {\rm for \ all \
} n  \ \ . $$ To avoid cumbersome expressions, let $\varpi$ denote
the following repeatedly occurring constant \begin{equation}
\varpi \ := \ \textstyle{ \frac{c_1}{96} } \  < \ 1 \ \  .
\label{varpi}
\end{equation}
In view of the  discussion of  \S\ref{subseqwork},  we can assume
that for $n$ sufficiently large
\begin{equation}
\r(u_{n+1})  \ \leq \  \lambda \,   \r(u_n) \hspace{8mm} {\rm with
\ } \hspace{6mm} 0 < \lambda  \, <  \, \textstyle{ \frac{1}{9} } \
. \label{lambda}
\end{equation}
\vskip 6pt \noindent{\bf Constructing the first level ${\bf K ( }
1{\bf )}$.} \vskip 6pt \noindent Choose $t_1$ large enough so that
\begin{equation}
g(\kt) \  <  \ G^* \  <  \frac{\eta}{24\varpi}  \ ,
\label{t1choiceb}
\end{equation}
\begin{equation} \p(u_{t_1})^{s-1}> \frac{\eta}{\varpi}  \ ,    \label{t1choiced}
\end{equation}
and so that the counting estimate (\ref{NGnB}) is valid for the
set  $G_{I_0}^u(t_1)$; i.e. with $B=I_0$.  Note that the first of
these inequalities is possible since   $g(u_n) < G^* < \eta $. The
latter inequality is possible  in view of (\ref{mythm33}) and
since $s<1$. Let $k_1 \geq 1 $ be the unique integer  such that
\begin{eqnarray}
\frac{6\varpi}{\eta} \ \sum_{i=0}^{k_1-1} g(u_{t_1+i}) \   \le \
\frac{1}{4}  \ < \ \frac{6\varpi}{\eta} \ \sum_{i=0}^{k_1}
g(u_{t_1+i})  \ .  \label{k1choicea}
\end{eqnarray}
Note, the fact that $k_1 \geq 1 $ is a consequence of
(\ref{t1choiceb}). The first level ${\bf K ( } 1{\bf )}$ of the
Cantor set ${\bf K_{\eta}}$ will  consist of sub-levels $K(t_1+
i)$  where $0 \, \leq \, i \, \leq k_1$.

\medskip

\noindent $\bullet$ \ {\bf The  sub-level $K(t_1)$ : \ } This
consists of balls of common radius $\p(u_{t_1})$ defined as
follows:- $$K(t_1)\ := \ \bigcup_{\alpha\in
V_{I_0}^u(t_1)}^{\circ} B(\ra,\p(\kt))   \hspace{5mm} {\rm where }
\hspace{5mm}
 V_{I_0}^u(t_1) := G_{I_0}^u(t_1) \ . $$
\noindent $\bullet$ \ {\bf The  sub-levels $K(t_1+i)$ for $1 \leq
 i  \leq k_1$  : \ } These are constructed inductively. The key to
 the whole  procedure is the existence of  `special'
 subsets $V_{I_0}^u(t_1+i)$ of $G_{I_0}^u(t_1+i)$. 
Suppose for $0 \leq j \leq i-1$  we have constructed the
sub-levels $$K(t_1+j) \ =\bigcup_{\alpha \in
V_{I_0}^u(t_1+j)}^{\circ}
 B(\ra,\p(u_{t_1+j}))  \ .  $$
  We proceed to construct $K(t_1+i)$ -- equivalently $V_{I_0}^u(t_1+i )$. Let
$$h(u_{t_1+j}) \ :=   \ \frac{\varpi
   }{\eta} \Psi(u_{t_1+j})^s  \ . $$ Note that in view of
(\ref{t1choiceb}) and (\ref{t1choiced}) we have that $
\p(u_{t_1+j})  <  h(u_{t_1+j})< \r(u_{t_1+j})$.
 Define $$T(t_1+j) \ := \ \left\{
B(\ra,h(u_{t_1+j})) :\alpha \in V_{I_0}^u(t_1+j) \right\} \ . $$
Now for each $\alpha \in G_{I_0}^u (t_1+i)$ construct the ball
$B(\ra,\r(u_{t_1+i}))$. Clearly, the balls in this collection are
also disjoint and we proceed by disregarding any of  those which
lie too close to balls from any of the previous sub-levels
$K(t_1+j)$. To make this precise, we introduce the sets
\begin{eqnarray*}
U_{I_0}^u (t_1+i) & := & \{\alpha \in G_{I_0}^u
(t_1+i):B(\ra,\r(\ktk))\cap \bigcup_{j=0}^{i-1}T(t_1+j)
\neq\emptyset\} \\ V_{I_0}^u(t_1+i)& := &G_{I_0}^u
(t_1+i)\setminus U_{I_0}^u (t_1+i) \ .
\end{eqnarray*}
\noindent By construction,  $V_{I_0}^u(t_1+j) \subseteq G_{I_0}^u
(t_1+j)$ for $0\le j \le i $. In particular,   the balls in $
T(t_1+j)$ are disjoint. Thus, $ \#T(t_1+j) =   \# V_{I_0}^u(t_1+j)
\leq \#G_{I_0}^u (t_1+j) $. We claim that $\#V_{I_0}^u(t_1+i) \geq
\frac{1}{2} \, G_{I_0}^u (t_1+i)  \ . $ There are two cases to
consider.

\vspace{3mm} \noindent{\rm \underline{Case} (i):  $ 0 \leq j \leq
i-1$ such that $ \r(u_{t_1+i}) < h(u_{t_1+j})$.  } The number of
disjoint balls of radius $\r(u_{t_1+i})$ that can possibly
intersect a ball in  $T(t_1+j)$  is $ \le 3h(u_{t_1+j})/
\r(u_{t_1+i})$.

\vspace{2mm} \noindent{\rm \underline{Case} (ii):  $ 0 \leq j \leq
i-1$ such that $\r(u_{t_1+i}) \geq h(t_1+j)$.  } In this case,
each ball in $T(t_1+j)$ can intersect at most one ball
$B(\ra,\r(\ktk))$ with $\alpha \in G_{I_0}^u (t_1+i)$. This makes
use of the fact that the corresponding enlarged  balls
$B(\ra,2\r(\ktk))$ are disjoint.

\vspace{3mm}

\noindent It follows that $$ \#U_{I_0}^u (t_1+i)  \le
\sum_{\rm{case\ (i)\ }} \frac{3\varpi}{\eta} \
\frac{\p(u_{t_1+j})^s}{ \r(u_{t_1+i})} \ \#T(t_1+j) +
\sum_{\rm{case\ (ii)\ }} \#T(t_1+j). $$ Recall that $ \#T(t_1+j)
\leq \#G_{I_0}^u (t_1+j) $. Thus, the contribution from the sum
over case (i) is: $$ \leq  \sum_{j=0}^{k_1-1} \frac{6\varpi}{\eta}
\ g(u_{t_1+j}) \ \#G_{I_0}^u (t_1+i) \ \le \ \frac{1}{4} \
\#G_{I_0}^u (t_1+i) \ , $$  by (\ref{NGnB}) and the choice of
$k_1$ -- see (\ref{k1choicea}). The contribution from the sum over
case (ii) is:
\begin{eqnarray*}
& \leq  &  \sum_{\rm{case\ (ii)\ }} \#G_{I_0}^u (t_1+j)
 \ \leq \
2  \ \sum_{j=0}^{i-1} \frac{\r(u_{t_1+i})}{\r(u_{t_1+j})}
\#G_{I_0}^u (t_1+i)
\\  & \leq &
2 \ \#G_{I_0}^u (t_1+i)\sum_{j=0}^{i-1} \lambda^{i-j} \  < \  2 \,
\#G_{I_0}^u (t_1+i)\sum_{s=1}^{\infty} \lambda^s  \ < \
\frac{1}{4} \ \#G_{I_0}^u (t_1+i) \ ,
\end{eqnarray*}
by (\ref{NGnB}) and the choice of $\lambda$ -- see (\ref{lambda}).
Hence,  $ \#U_{I_0}^u (t_1+i) < \frac{1}{2} \ \#G_{I_0}^u (t_1+i)
$ so that \begin{equation}
 \#V_{I_0}^u(t_1+i) \ \geq \
\textstyle{\frac{1}{2} } \ \#G_{I_0}^u (t_1+i) . \label{ohno1}
\end{equation} The
 sub-level $K(t_1+i)$ is defined to be:
$$K(t_1+i) \ := \ \bigcup_{\ra\in V_{I_0}^u(t_1+i)}
 B(\ra,\p(\ktk)) \ . $$
Also, note that by construction $K(t_1+i)  \cap  K(t_1+j)  =
\emptyset  $  for $ 0 \leq i \neq j \leq k_1 $. The first level
${\bf K ( }1{\bf )}$ of the Cantor set is defined to be $$ {\bf K
( }1{\bf )} \ := \ \bigcup_{i=0}^{k_1} K(t_1+i)  \ . $$

\vskip 6pt \noindent{\bf Higher levels ${\bf K(}n{\bf)} $ and the
Cantor set ${\bf K_{\eta}}$. } \vskip 6pt \noindent  For any
integer $n \geq 2$, the n'th level ${\bf K(}n{\bf)}$ will be
defined recursively in terms of local levels ${\rm K}(n,B)$
associated with balls $B$ from the previous level ${\bf K (
}n-1{\bf )}$:
$$ {\bf K ( }n{\bf )} \ := \ \bigcup_{B \in {\bf K ( }n-1{\bf )} }^{\circ}
{\rm K}(n,B)  \ , $$ where $${\rm K}(n,B) \ := \
\bigcup_{i=0}^{k_n(B)}K(t_n+i,B)  \ .$$ To start with, choose $t_n
> t_{n-1} $ sufficiently large so that for any ball $B\in{\bf K ( }n-1{\bf )}$
the counting estimate (\ref{NGnB}) is valid and so that
\begin{equation}
\p (u_{t_n})^{s-1}  > \frac{r(B)^{s-1}}{\varpi}.
 \label{tnchoice}
\end{equation}
\noindent 
In view of (\ref{mythm33}), (\ref{t1choiced}), the fact that
$g(u_n) < G^*$ for all $n$  and that $s<1$,  we have that
\begin{equation}
g(u_{t_n}) \ < \ G^*  \ < \frac{r(B)^{s-1}}{24\varpi}
\hspace{10mm} \forall \ \ \ \ B \in {\bf K ( }n-1{\bf )} .
\label{tngchoice}
\end{equation}
Fix a ball $B$  in  ${\bf K ( }n-1{\bf )}$ and let $k_n(B) \geq 1
$ be the unique integer such that
\begin{eqnarray}
\frac{ 6\varpi}{r(B)^{s-1}}
 \ \sum_{i=0}^{k_n(B)-1} g(u_{t_n+i}) \  \le   \
\frac{1}{4}  \ < \ \frac{6\varpi}{r(B)^{s-1}}
 \ \sum_{i=0}^{k_n(B)} g(u_{t_n+i}) \  .  \ \label{knchoicea}
\end{eqnarray}
The fact that $k_n(B) \geq 1 $ is a consequence of
(\ref{tngchoice}). We now construct the local level ${\rm
K}(n,B)$.
\medskip

\noindent $\bullet$ \ {\bf The local sub-level $K(t_n,B)$ : \ }
Let $$ K(t_n,B) \ := \ \bigcup_{\alpha \in V_B^u(t_n)}
B(\ra,\p(u_{t_n})) \hspace{5mm} {\rm where } \hspace{5mm}
 V_{B}^u(t_n) := G_{B}^u(t_n) \ . $$

 By construction, $ K(t_n,B) \subset B $ -- see (\ref{r<p}).
\medskip

 \noindent $\bullet$ \ {\bf The local sub-levels
$K(t_n+i,B)$ for $1 \leq i \leq k_n(B)$ :   }  Suppose
 for $0 \leq j \leq i-1$  we have constructed the
local sub-levels $$K(t_n+j,B) \ =\bigcup_{\alpha \in
V_{B}^u(t_n+j)}^{\circ}
 B(\ra,\p(u_{t_n+j})) \   \ .  $$
 Let
 $$ h_B(u_{t_n+j}) := \frac{\varpi \,  \p(u_{t_n+j})^s}{r(B)^{s-1}} \ .
$$ In view of (\ref{tnchoice}) and  (\ref{tngchoice}) we have that
\begin{equation}
 \p(u_{t_n+j}) \ < \ h_B(u_{t_n+j}) <\r(u_{t_n+j}) .
\label{p<h}
\end{equation}
 Define
$$T(t_n+j,B) := \left\{B(\ra,h_B(u_{t_n+j})) : \alpha \in
V_B^u(t_n+j) \right\} \ . $$ Next, introduce the sets
\begin{eqnarray*}
 U_{B}^u(t_n+i)&:=& \{\alpha\in
G_{B}^u(t_n+i):B(\ra,\r(u_{t_n+i}))\cap
\textstyle\bigcup\limits_{j=0}^{i-1} T(t_n+j,B)\neq\emptyset \}\\
\\ V_{B}^u(t_n+i)&:=&G_{B}^u(t_n+i)\setminus U_{B}^u(t_n+i) \  \ \ .
\end{eqnarray*} By construction, $V_{B}^u(t_n+j) \subseteq G_{B}^u
(t_n+j)$ for $0\le j \le i $ and so the balls in $ T(t_n+j)$ are
disjoint. By adapting the argument used in establishing
(\ref{ohno1}),  it is easily verified that
\begin{equation}
\#V_{B}^u(t_n+i) \ \geq \ \textstyle{\frac{1}{2} } \ \#G_{B}^u
(t_n+i) . \label{ohno2}
\end{equation}
 Now let
$$ K(t_n+i,B) \ := \ \bigcup_{\alpha\in V_B^u(t_n+i)}
 B(\ra,\p(u_{t_n+i})) \ .
$$ This completes the inductive step and the construction of the
local level ${\rm K}(n,B)$ associated with $B \in {\bf K (
}n-1{\bf )} $. Clearly, for $0\leq i\neq j \leq k_n(B)$ we have
that $$ K(t_n+i,B) \ \cap \ K(t_n+j,B) \ = \ \emptyset \ . $$
Furthermore, by construction  ${\rm K}(n,B)$  is contained in $B$.
Therefore, ${\bf K ( }n{\bf )}  \subset  {\bf K ( }n-1{\bf )} $.
The Cantor set ${\bf K_{\eta} }$ is
 defined as $$ \textstyle {\bf K_{\eta} } \ := \ \bigcap\limits_{n=1}^{\infty}
{\bf K ( }n{\bf )} \ \ . $$ Strictly speaking, $ {\bf K_{\eta}
}\setminus {\cal R} \subset  \Lambda({\cal R}, \beta,\Psi) $ and
not  ${\bf K_{\eta} }\subset  \Lambda({\cal R}, \beta,\Psi) $.
However, this is irrelevant  since ${\cal R} $ is countable  and
so ${\cal H}^s({\bf K_{\eta} }\setminus {\cal R}) = {\cal
H}^s({\bf K_{\eta} })$.  Before  constructing  a measure on ${\bf
K_{\eta} }$, we state an important lemma.  The proof is a simple
consequence of (\ref{NGnB}), (\ref{ohno1}) and (\ref{ohno2}).

\begin{lemma}\label{lem:sums}
\begin{enumerate}
\item[{\rm (i)}]  For $0\le i\le k_1$, $$
\# V_{I_0}^u(t_1+i)   \ \p(u_{t_1+i} )^s  \ \ge \
\frac{c_1|I_0|}{8} \ g(u_{t_1+i})   \ \ . $$
\item[{\rm (ii)}] For $n\ge 2$, let $B$ be a ball in ${\bf K ( }n-1{\bf )}$.
Then, for
 $0\le i\le k_n(B)$
$$\#  V_B(t_n+i) \  \p (u_{t_n+i})^s  \ \ge
 \ \frac{c_1|B|}{8} \ g(u_{t_n+i})  \ \ . $$
\end{enumerate}
\end{lemma}

\myappsubsubsection{A measure on ${\bf K_{\eta} }$}

In this section, we construct a  probability measure $\mu$
supported on ${\bf K_{\eta} }$ satisfying (\ref{task}). Suppose $n
\geq 2$ and $B \in {\bf K ( }n{\bf )}$. For $1\leq m < n$, let
$B_m$ denote the unique ball in ${\bf K ( }m{\bf )}$ containing
the ball $B$. With this notation in mind we now define $\mu$. For
any $B \in {\bf K ( }n{\bf )}$, we attach a weight $\mu(B)$
defined recursively as follows: For  $n=1$, $$ \mu(B)\ := \
\frac{r(B)^s}{\sum_{B'\in {\bf K ( }1{\bf )}} r(B')^s} \ $$ and
for $n\ge 2$, $$ \mu(B) \ := \ \frac{r(B)^s}{\sum_{B'\in
K(n,B_{n-1})} r(B')^s} \ \times \ \mu(B_{n-1})  \ . $$
\medskip

\noindent This procedure defines inductively a mass on any ball
appearing in the construction of ${ \bf K_{\eta}} $. In fact a lot
more is true: The probability measure $\mu$ constructed above is
supported on ${\bf K_{\eta} }$ and for any Borel subset $F$ of
$I_0$
\[
\mu(F):= \mu(F \cap {\bf K_{\eta} })  \; = \;
\inf\;\sum_{B\in{\cal B}}\mu(B)  \ ,
\]
where the infimum is taken over all coverings $\cal B$ of $F \cap
{\bf K_{\eta} }$ by balls  $B\in \{{\bf K ( }n{\bf )} : n\in\N\}$.

\medskip

For further details see  \cite[Prop. 1.7]{falc} . It remains to
establish  (\ref{task}) for   $\mu \, $.

\vskip 5pt

\noindent{\bf Measure of a ball in the Cantor construction.} If
$B\in{\bf K ( }n{\bf )}$ for some $n\in \N$, then
\begin{eqnarray}
\mu(B) & := & \frac{r(B)^s}{\sum_{B'\in K(n,B_{n-1})}r(B')^s} \
\times  \ \mu\left(B_{n-1}\right) \nonumber \\ & & \nonumber \\ &
= & \frac{r(B)^s}{\sum_{B'\in{\bf K ( }1{\bf )}}r(B')^s} \;
\prod_{m=1}^{n-1} \frac{r(B_m)^s}{\sum_{B'\in
K(m+1,B_{m})}r(B')^s} \label{measball} \;.
\end{eqnarray}

\noindent The  product term is taken to be one  when $n=1$. To
proceed we require the following lemma which gives us a lower
bound on the terms in the denominator of the above expression.

\begin{lemma}

\[\textstyle \sum\limits_{B\in{\bf K ( }1{\bf )}}r(B)^s \ \ge \
\textstyle{\frac{\eta}{2}}   \; |I_0| \hspace{7mm} {\rm  and   }
\hspace{7mm} \sum\limits_{B\in K(n, B_{n-1})} \!\!\!\!\!  r(B)^s
\ \ge \ r(B_{n-1})^s \   \ \ \ ( n\ge2) . \] \label{triv}
\end{lemma}

\noindent{\bf Proof of Lemma \ref{triv}.}~ By Lemma
\ref{lem:sums}, the choice of $k_1$ (\ref{k1choicea}) and $\varpi$
(\ref{varpi}) it follows that $$ \sum_{B\in {\bf K ( }1{\bf )}}
r(B)^s \ =  \ \sum_{i=0}^{k_1}\ \ \# V_{I_0}^u(t_1+i) \  \p(
u_{t_1+i})^s \ \ge  \ \frac{c_1|I_0|}{8} \ \sum_{i=0}^{k_1}
g(u_{t_1+i}) \ > \ \frac{c_1|I_0|}{192} \frac{\eta}{\varpi} \ \ge
\ \textstyle{\frac{\eta}{2}} \, |I_0|  \ . $$

\noindent The statement for $n \geq 2$ follows in a similar
fashion --  use (\ref{knchoicea}) rather than (\ref{k1choicea}).
\hfill $\spadesuit$

 In view of the above lemma, it now follows from
(\ref{measball}) that for any ball $B\in{\bf K ( }n{\bf )}$

\begin{equation}
 \mu(B) \ \le \ \frac{2 \, r(B)^s}{|I_o| \ \eta}  \ \ll \ \frac{ r(B)^s}{ \eta} \ .
 \label{bink}
 \end{equation}
\vskip 4pt \noindent{\bf Measure of an arbitrary ball.}  The aim
is to show that $ \mu(A)  \ll  r(A)^s/\eta $ for an arbitrary ball
$A$ with radius $r(A) \leq r_o $. The measure $\mu$ is supported
on ${\bf K_{\eta}}$. Thus, without loss of generality we can
assume that $ A \cap {\bf K_{\eta}} \neq \emptyset $; otherwise
$\mu(A) = 0$ and there is nothing to prove. We can also assume
that for every $n$ large enough $A$ intersects at
 least two balls in  ${\bf K ( }n{\bf )}$; since if $B$ is the only ball in  ${\bf K ( }n{\bf )}$
which has non--empty intersection with  $A$, then in view of
(\ref{bink})
 $$
 \mu(A) \ \leq \ \mu(B) \ \ll  \  r(B)^s/\eta \ \to \ 0
 \hspace{8mm}  {\rm as } \hspace{5mm} n \to \infty  \
 $$
 ($r(B) \to 0$ as $n \to \infty$) and again there is nothing to
 prove. Thus we may assume that there exists an integer $n \geq 2$
 such that $A$ intersects only one ball $\widetilde{B}$ in ${\bf K ( }n-1{\bf )}$
 and at least two balls from ${\bf K ( }n{\bf )}$. The case that $A$
 intersects two or more balls from the first level can be excluded
 by choosing $r(A)$ sufficiently small. This follows from the fact
 that by construction balls in any one level are disjoint.
 Furthermore, we can assume that
 $$
 r(A) \ < \ r(\widetilde{B})  \ .
 $$
 Otherwise,
 $
 \mu(A) \ \leq \ \mu(\widetilde{B}) \ \leq \ r(\widetilde{B})^s/\eta
 \ \leq \ r(A)^s/\eta
 $
 and we are done.
Given that $A$ only intersects the ball $ \widetilde{B}$ in ${\bf
K ( }n-1{\bf )}$, the balls from level ${\bf K ( }n{\bf )}$  which
intersect $A$ must be contained in the local level $$ {\rm
K}(n,\widetilde{B}) \ := \
\bigcup_{i=0}^{k_n(\widetilde{B})}K(t_n+i,\widetilde{B})  \ .
$$ By construction, any ball $B(\ra,\p(u_{t_n+i}))$ in ${\rm
K}(n,\widetilde{B})$ is contained in the ball
$B(\ra,h_{\widetilde{B}}(u_{t_n+i}))$.
 Thus $A$ intersects at least one
ball in $T(t_n+i,\widetilde{B})$ for some $ 0 \leq i \leq
k_n(\widetilde{B})$.

 Let $K(t_n+i',\widetilde{B})$ be the first
local sub-level associated with $\widetilde{B}$ such that $$
K(t_n+i',\widetilde{B}) \ \cap \ A \ \neq \ \emptyset \  . $$
Thus, $A$ intersects at least one ball $B(\ra,\p(u_{t_n+i'}))$
from $K(t_n+i',\widetilde{B})$ and such balls are indeed the
largest balls from the $n$'th level ${\bf K ( }n{\bf )}$ that
intersect $A$. Clearly, $A$ intersects at least one ball $ B_*$ in
$T(t_n+i',\widetilde{B})$. We now prove a trivial but crucial
geometric lemma.

\begin{lemma}
\label{geo} For $i \geq i'$, if $A$ intersects
$B(\ra,\p(u_{t_n+i})) \subset B(\ra,h_{\widetilde{B}}(u_{t_n+i}))
\neq B_*$ then $$ r(A)  > \textstyle{ \frac{1}{2} } \;
\r(u_{t_n+i}) \ . $$
\end{lemma}

\noindent{\bf  Proof of Lemma \ref{geo}.} If $i=i'$ then as balls
of radius
 $\r(u_{t_n+i'})$ are disjoint we have that $r(A)>\r(u_{t_n+i'})$. Assume that
 $i>i'$, then
 by construction $ B(\ra,\r(u_{t_n+i})) \cap
B_* = \emptyset $.   Hence,
 $r(A)>\r(u_{t_n+i})-\p(u_{t_n+i})>\frac12\r(u_{t_n+i})$ -- see
 (\ref{mythm33}).
\hfill $\spadesuit$

\medskip

In view of the definition of $i'$ and (\ref{bink}), we have that
\begin{eqnarray}
\mu(A) & \leq  & \sum_{i \, = \, i'}^{k_n(\tilde{B})}
\sum_{\substack{ \alpha \in \, V_{\tilde{B}}^u(t_n+i) \, : \\
B\left(\ra,\p(u_{t_n+i})\right) \cap A \neq \emptyset} }
\mu(B(\ra,\p(u_{t_n+i}))) \nonumber \\ & & \nonumber
\\ & \leq & \frac{2}{\eta \, |I_0|} \ \sum_{i \, = \, i'}^{k_n(\tilde{B})}
\p(u_{t_n+i})^s \ \sum_{\substack{\alpha \in
V_{\tilde{B}}^u(t_n+i) \, : \\ B\left(\ra,\p(u_{t_n+i})\right)\cap
A\neq\emptyset }} \!\!\!\!\!\!\!\!\!\!\!\!\!\!\! 1     \label{ab}
\  \  .
\end{eqnarray}

In view of Lemma \ref{geo}, if $A$ intersects some ball
$B(\ra,h_{\widetilde{B}}(t_n+i))$ in $T(t_n+i,\widetilde{B})$ then
the ball $B(\ra,\r(u_{t_n+i}))$ which contains it is itself
contained in the ball $5  A$. Let $N_i$ denote the number of balls
$B(\ra,\r(u_{t_n+i}))$ with $\alpha \in   V_{\tilde{B}}^u(t_n+i) $
that can possibly intersect $A$. By construction these balls are
disjoint. Thus, $ 2 \r(u_{t_n+i}) \times N_i \ \leq \  |5A|   = 10
\; r(A)  $.
 This implies, via (\ref{ab}) that
$$ \mu(A) \ \leq \   \frac{2}{\eta \, |I_0| } \ \sum_{i \, = \,
i'}^{k_n(\tilde{B})} \p(u_{t_n+i})^s \   \;N_i
  \ \leq  \  \frac{10 }{\eta \, |I_0| } \
r(A)  \  \sum_{i \, = \, 0}^{k_n(\tilde{B})} g(u_{t_n+i}) \ . $$
By (\ref{knchoicea}), $$  \sum_{i=0}^{k_n(B)-1} g(u_{t_n+i}) \
\leq  \   \frac{r(\widetilde{B})^{s-1}}{ 24 \varpi } \  , $$ and
by (\ref{tngchoice}) together with the fact that $g(u_n) < G^*$
for all $n$  $$ g(u_{t_n+k_n(\tilde{B})})  \ < \
\frac{r(\widetilde{B})^{s-1}}{24\varpi} \ .$$

\noindent Hence
\begin{equation} \mu(A) \ \ll \ \frac{1}{\eta} \
r(A) \ r(\widetilde{B})^{s-1} \label{case1ab}
\end{equation}

\noindent However, $r(A) < r(\widetilde{B}) $ and $s<1$.  The
desired inequality, namely (\ref{task}) now readily follows. This
completes the proof of Theorem \ref{xthm3} in the case that $G$ is
finite. \hfill $\spadesuit$

\myappsubsection{Proof of  Theorem \ref{xthm3}:  $G= \infty$
\label{pfdimthm}}

The proof of  Theorem \ref{xthm3}  in the case that  $G$ is
infinite follows the same strategy as the proof when $G$ is
finite.  However, to execute the strategy is far simpler than in
the finite case.

\myappsubsubsection{The Cantor set {\bf K} and the measure $\mu$ }

We start by defining a Cantor set ${\bf K}$ which is dependent on
a certain, strictly increasing sequence of natural numbers $\{t_i
: i \in \N \}$. The main difference between this case and the
previous case is that we do not need to consider sublevels.

\medskip

\noindent{\em The Cantor set {\bf K}.  \  } Choose $t_1$
sufficiently large so that the counting estimate  (\ref{NGnB}) is
valid for the set $G_{I_0}^u(t_1)$ and define the first level
${\bf K ( }1{\bf )}$ of the Cantor set ${\bf K}$  as follows: $$
\textstyle {\bf K ( }1{\bf )}\ := \ \bigcup\limits_{\alpha\in
G_{I_0}^u(t_1)} B(\ra, \p(\kt)) \ . $$ For $n \geq 2 $ we define
the $n$'th level ${\bf K ( }n{\bf )}$ recursively as follows:
$$\textstyle {\bf K ( }n{\bf )} \ := \ \bigcup\limits_{B \in {\bf K ( }n-1{\bf )} }^{\circ} {\rm K}(n,B) \ , $$
where $$ \textstyle {\rm K}(n,B) \ := \ \bigcup\limits_{\alpha\in
G_{B}^u(t_n)} B(\ra,\p(u_{t_n})) $$ is the $n$'th local level
associated with the ball $B:= B(\ra, \p(u_{t_{n-1}})) \in {\bf K (
}n-1{\bf )} $. Here  $t_n
> t_{n-1} $ is chosen sufficiently large so that (\ref{NGnB}) is valid  for any ball
$B$ in ${\bf K ( }n-1{\bf )}$.  By construction ${\rm K}(n,B)
\subset B$ and so  ${\bf K ( }n{\bf )}  \subset {\bf K ( }n-1{\bf
)}$. The Cantor set ${\bf K}$ is simply given by
$$\textstyle {\bf K } \ := \ \bigcap\limits_{n=1}^{\infty} {\bf K
( }n{\bf )} \ \ . $$

By construction, $ \, {\bf K }\setminus {\cal R } \subset
\Lambda(\cR,\beta, \Psi)$ and since ${\cal R}$ is countable ${\cal
H}^s ( \Lambda(\cR,\beta, \Psi)) \geq {\cal H}^s ( {\bf K })$.

\medskip

\noindent{\em The measure $\mu$.   } Suppose $n \geq 2$ and $B \in
{\bf K ( }n{\bf )}$. For $1\leq m <  n$, let $B_m$ denote the
unique ball in ${\bf K ( }m{\bf )}$ containing  $B$. For any $B
\in {\bf K ( }n{\bf )}$, we attach a weight $\mu(B)$ defined
recursively as follows:
\medskip

For  $n=1$, $$\mu(B)\ := \ \frac{1}{ \#G_{I_0}^u(t_1) } $$

and for $n\ge 2$, $$\mu(B) \ := \ \frac{1}{ \#G_{ B_{n-1} }^u
(t_n) } \times \ \mu(B_{n-1}) \ . $$

\noindent By the definition of $\mu$ and the counting estimate
(\ref{NGnB}), it follows that
\begin{eqnarray}
\mu(B) & \leq  &  \mbox{\large $\frac{2}{|I_0|} $ } \ c_1^{-n} \
\r(u_{t_n})\
  \ \times \
\prod_{m=1}^{n-1} \frac{\r(u_{t_m})}{\p(u_{t_m})} \ .
\label{measballubexact}
\end{eqnarray}
\indent  The  product term is taken to be one when $n=1$.

\myappsubsubsection{Completion of the proof }

Fix $\eta \geq 1 $. Since $G = \infty $, the sequence $\{t_i\}$
associated with the construction of the Cantor set ${\bf K }$ can
clearly be chosen so that
\begin{equation}
\eta \ \times \ \mbox{\large $\frac{2}{|I_0|} $ } \ c_1^{-i} \
\times \ \prod_{j=1}^{i-1} \frac{\r(u_{t_j})}{\p(u_{t_j})} \leq \
g(u_{t_i}) \label{productchoice} \ \ . \end{equation} The product
term is one when $i=1$. It now immediately follows from
(\ref{measballubexact}) that for  any $B \in {\bf K ( }n{\bf )}$,
\begin{equation*}
\mu(B) \   \leq  \ r(B)^s \ \mbox{\large $\frac{2}{|I_0|} $ } \
c_1^{-n} \ \frac{1}{g(u_{t_n})} \ \times \ \prod_{m=1}^{n-1}
\frac{\r(u_{t_m})}{\p(u_{t_m})}
 \ \leq \ r(B)^s /\eta   \ .
\label{bdone}
\end{equation*}

We now show that $\mu(A) \ll  r(A)^s / \eta $ where $A$ is an
arbitrary ball. The same reasoning as before
 enables us to assume that $A\cap {\bf K } \neq
\emptyset $ and that
 there exists an integer $n \geq 2$
 such that $A$ intersects only one ball $\widetilde{B}$ in ${\bf K ( }n-1{\bf )}$
 and at least two balls from ${\rm K}(n,\widetilde{B}) \subset {\bf K ( }n{\bf )}$. Thus,
\begin{equation}
 \rho(u_{t_n})  \
\leq \ r(A) \ \leq \  r(\widetilde{B}) \, :=  \, \Psi(u_{t_{n-1}})
\ \ . \label{rineq}
\end{equation} The left
hand side of (\ref{rineq})  makes use of   the fact that $B(\ra,
\psi(u_{t_n})) \subset  B(\ra, \rho(u_{t_n})) $  and that the
 balls $ B(\ra, 2\rho(u_{t_n})) $ with  $\alpha \in
G_{\widetilde{B}}^u(t_n) $ are  disjoint.  A simple geometric
argument yields  that $ N := \# \{ \alpha \in G_{\widetilde{
B}}^u(t_n) :   B(\ra, \rho(u_{t_n})) \cap A \neq \emptyset \} \
\leq \  3 \, r(A)/\r(u_{t_n})  \ .  $ In view of
(\ref{measballubexact}), (\ref{productchoice}), (\ref{rineq}) and
the fact that $s<1$,  we obtain
\begin{eqnarray*}
\mu(A) & \leq &     N \   \mu(B(\ra,\p(u_{t_n})))
 \  \le \    r(A)  \ \mbox{\large
$\frac{6}{|I_0|} $ } \ c_1^{-n} \ \prod_{m=1}^{n-1}
\frac{\r(u_{t_m})}{\p(u_{t_m})}
\\ & \le &  r(A)^s   \  \p(u_{t_{n-1}})^{1-s}     \   \mbox{\large
$\frac{6}{|I_0|} $ } \ c_1^{-n}   \ \prod_{m=1}^{n-1}
\frac{\r(u_{t_m})}{\p(u_{t_m})}
\\ & \le &  r(A)^s  \ \ \mbox{\large
$\frac{6}{|I_0|} $ } \ c_1^{-n}  \ \frac{ 1}{g(u_{t_{n-1}}) }  \
 \prod_{m=1}^{n-2} \frac{\r(u_{t_m})}{\p(u_{t_m})} \ \ \leq \
3 \, c^{-1} \ \frac{r(A)^s}{\eta} \ .
\end{eqnarray*}

The upshot is that (\ref{task}) is satisfied and thereby completes
the proof. \hfill $\spadesuit$
\end{appendix}


\vspace{6mm}

\noindent{\bf Acknowledgements.} We would like to thank the
referee for making  many useful suggestions. In particular, it was
the referees comments which lead us  to the convergent statements
for rational quadrics -- thank you for sharing  your insight.
Previously, we had only obtained these statements for the unit
circle. \\ \indent As ever, SV would like to thank his old friend
Bridget and his new friends Ayesha and Iona for bringing so much
love and laughter into his life.

\vspace{5mm}

\noindent Victor V. Beresnevich:  Institute of Mathematics,
Academy of Sciences of Belarus,

\vspace{-1mm}
 ~ \hspace{30mm} 220072, Surganova 11, Minsk, Belarus.

\vspace{0mm}
~ \hspace{30mm} e-mail: beresnevich@im.bas-net.by

\vspace{4mm}

\noindent Detta Dickinson:  Department of Mathematics, National
University of Ireland,

\vspace{-1mm}
 ~\hspace{21mm}  Maynooth,   Co.Kildare, Ireland.

\vspace{0mm}
 ~ \hspace{21mm} e-mail: : ddickinson@maths.may.ie

\vspace{4mm}

\noindent Sanju L. Velani: Department of Mathematics, University
of York,

\vspace{-1mm}
 ~ \hspace{19mm}  Heslington, York, YO10 5DD, England.

 \vspace{0mm}
 ~ \hspace{19mm} e-mail: slv3@york.ac.uk
\newpage
\begin{appendix}

\begin{center} {\Large \bf Appendix II:  \  Sums of Two Squares Near  \\ ~ \hspace*{10ex} Perfect
Squares} \\ ~  \\   {\large R.C. Vaughan }
\addcontentsline{toc}{section}{{\bf APPENDIX ~II:~~ Sums of Two
Squares Near Perfect Squares}}
\end{center}

\vspace{7mm}

\noindent{\bf \large{A.1 \ \  The theorem}}

\noindent Let $r(n)$ denote the number of representations of a
number $n$ as the sum of two squares of integers and let
$\psi:{\Bbb N}\to{\Bbb R}$ be a non-negative decreasing function.
We prove the following theorem.

\begin{thdel}
 Let $Q^*$ denote the smallest integer with $Q^*>Q$.
Then for each real number $Q$ and natural number $N$ with $N\le
Q^3$,
\begin{eqnarray*}
\sum_{Q<q\le 2Q}&& \!\!\!\!\!\!\!\!\!\! {\sum_n}'r(n)  =
\sum_{Q<q\le 2Q} 4\pi q\psi(q)
\\  \\ &+ & O\bigg( Q\log Q + Q^{\frac32}\psi(Q^*)^{\frac12}(\log
Q)^{64} + Q^2\psi(Q^*)^{\frac12} (\log Q)^{64}N^{-\frac14} \\ \\ &
& \ \ \ \ \ \ +  \ \ N^{\frac34}(\log N)^3 Q^{\frac12} \psi(Q^*) +
N^{\frac14}(\log N)Q^{\frac12} \sum_{Q<q\le 2Q} \psi(q) + Q^2(\log
Q)^3N^{-\frac12} \bigg)
\end{eqnarray*}
where $\sum'$ indicates that the sum is over $n$ with
$|q-\sqrt{n}|\le\psi(q)$ and that any terms with
$|q-\sqrt{n}|=\psi(q)$ are counted with weight $\frac{1}{2} $.
\label{thmbob}
\end{thdel}

 When $\psi(Q^*)$ has the same order of magnitude as
$Q^{-1}\sum_{Q<q\le 2Q}\psi(q)$ and the sum $\sum_{Q<q\le
2Q}\psi(q)$ is large, a good choice for $N$ is $$ \textstyle{ Q^2
\left( \sum\limits_{Q<q\le 2Q}\psi(q) \right)^{-1} \  } .$$ This
leads to the error estimate $$\ll  \textstyle{ Q\left(
\sum\limits_{Q<q\le 2Q} \psi(q) \right)^{\frac34} (\log Q)^{64}  \
. } $$
 Then the main term will
dominate provided that $\sum_{Q<q\le 2Q}\psi(q)$ is large compared
with $(\log Q)^{256}$. A concomitant remark pertains if the
theorem is averaged over $Q$ with, say $R<Q\le 2R$. It may well be
possible to replace the $(\log Q)^{64}$ in the error term by a
smaller power of $\log Q$, but that some power of a logarithm has
to be present follows from either of the observations that
$$\textstyle\sum\limits_{q\le Q}r(q^2)\sim \frac4{\pi} \, Q\log
Q$$ (see (\ref{tag1.7}) below) and $$ \textstyle \sum\limits_{q\le
Q}r(q^2+1)\sim \frac{12}{\pi}Q\log Q.$$

\vspace{3mm}

\noindent{\bf \large{A.2 \ \  Proof of Theorem \ref{thmbob}} }

\noindent Let $R(x)=\sum_{1\le n\le x}r(n)$, $\Delta(x)=R(x)-\pi
x$, and $\Delta_0(x)=\Delta(x)$ when $ x\not\in{\Bbb N}$ and
$\Delta_0(x)=\Delta(x)-\frac12r(x)$ when $x\in{\Bbb N}$. Then our
motivation is the formula of Hardy \cite[pg 265]{Hardy15}  which,
for real $x\ge x_0$, we restate as $$\Delta_0(x)=-1+
\sqrt{x}\sum_{n=1}^{\infty}r(n)n^{-1/2} \ J_1\left( 2\pi \sqrt{nx}
\right)$$ where $J_1$ denotes the usual Bessel function. However
the convergence is only conditional and we require a form of this
in which the tail of the infinite series is more readily
accessible.

By Theorem 1 and Lemma 2 of \cite{Hardy25} for any $\delta$ with
$0<\delta<1$ and provided that $x\ge x_0(\delta)$ and
$N>N_0(\delta)$ we have
\begin{eqnarray*} \Delta_0(x) \  =  \ -1  & + & \sqrt{x}\sum_{1\le n\le
N}\frac{r(n)}{n^{\frac12}} \; J_1\left( 2\pi\sqrt{n x} \right) \\
& - & x^{\frac54}\sum_{x(1-\delta) < n < x(1+\delta)}
\frac{r(n)}{\pi n^{\frac54}} \  {\rm{sgn}} \Bigg(
\sqrt{\frac{n}{x}}-1 \Bigg) \int_{2\pi|\sqrt{n} -
\sqrt{x}|\sqrt{N}}^{\infty} \frac{\sin{\alpha}}{\alpha} \, d\alpha  \\
& + &  \  O\left(
(xN)^{-\frac14} + x^{\frac14}N^{-\frac5{12}} \right)\\
\end{eqnarray*}
where we have used $\Delta(x)\ll x^{1/3}$ of \cite{Sie} and
${\rm{sgn}}(u)$ is $-1$, $0$ or $1$ according as $u<0$, $u=0$ or
$u>0$. A standard estimate for $J_1$ \cite[pg 199]{Watson} gives
for $\alpha>\alpha_0$
$$J_1(2\pi\alpha) \ = \ -\frac1{\pi} \; {\alpha}^{-1/2} \cos\left( 2\pi
\alpha +\frac{\pi}4 \right) + O\left( {\alpha}^{-3/2} \right) \
.$$ For convenience we put
$$S(x) \ = \ x^{1/4}\sum_{1\le n\le N}r(n)n^{-3/4}\cos\left(
2\pi\sqrt{nx}+\frac{\pi}4 \right)$$ and
$$E(x)\ = \ x^{\frac54}\sum_{x(1-\delta) < n < x(1+\delta)}
\frac{r(n)}{n^{\frac54}} \ \  {\rm{sgn}}\Bigg(
\sqrt{\frac{n}{x}}-1 \Bigg) \int_{2\pi|\sqrt{n} -
\sqrt{x}|\sqrt{N}}^{\infty} \frac{\sin{\alpha}}{\alpha} \, d\alpha
$$ so that
\begin{eqnarray*}
\Delta_0(x)&= & -1 -\pi^{-1}\big(S(x)+ E(x)\big) + O\left(
x^{-\frac14} + x^{\frac14}N^{-\frac5{12}} \right)\\ &=&
-1 -\pi^{-1}\big(S(x)+ E(x)\big) + O\left( (x/N)^{\frac12} \right)\\
\end{eqnarray*}
\noindent since $x^{-\frac14} + x^{\frac14}N^{-\frac5{12}}\ll
(x/N)^{\frac12}$ whenever $N\ll x^{3/2}$. For $x_0\le x\le y$ we
have $$S(y)-S(x)=\sum_{1\le n\le N }\frac{r(n)}{n^{3/4}}
\int_x^y\Re \left( \left( \frac14 u^{-3/4}+\pi i u^{-1/4}n^{1/2}
\right)e\left( \sqrt{nu}+\frac18 \right) \right)du.$$ The
contribution to $S(y)-S(x)$ from the
$$\frac14 u^{-3/4}e\left( \sqrt{nu}+\frac18 \right)$$
part of the integrand is $\ll x^{-3/4}N^{1/4}(y-x)$.  Here we have
used partial summation and the fact that $r(n)$ is on average
$\pi$.  We shall do this several times hereafter without comment.

\medskip

To prove the theorem we may suppose that $Q>Q_0$.  In particular
$Q_0$ can be chosen so that $q-\psi(q)>2$ whenever $q>Q$.  Thus,
when $Q<q\le 2Q$,
\begin{eqnarray}
\!\!\!\!\!\!\!\!\!\!\!\!\!\!\!\!\!\!\!\!\!\!\!\! \Delta_0
\big((q+\psi(q))^2\big)&-&\Delta_0\big((q-\psi(q))^2\big) \ \ = \
\ - \frac1{\pi} \; T(q,N) \nonumber
\\  & & \!\!\!\!\!\!\!\!\!\!\!\!\!\!\!\!\!\!\!\!\!\!\!\!\!
 \ + \ O\left( N^{\frac14}q^{-\frac12}\psi(q) +
E_+(q,N) + E_-(q,N)+ qN^{-\frac12} \right) \label{tag1.2}
\end{eqnarray}
where $$T(q,N) \ = \ \sum_{1\le n\le N }\frac{r(n)}{n^{1/4}}
\int_{(q-\psi(q))^2}^{(q+\psi(q))^2}\Re \left( \pi i
u^{-1/4}e\left( \sqrt{nu}+\frac18 \right) \right)du$$ and $$
E_{\pm}(q,N)\ = \ \sum_{\frac12Q^2<n<8Q^2}r(n) \min\left( 1,
\frac1{|\sqrt{n}-(q\pm \psi(q)|\sqrt{N}} \right)\ .$$ In the
integral in $T(q,N)$ we make the change of variables, $u=(q+t)^2$,
so that $$T(q,N)=\sum_{1\le n\le N }\frac{r(n)}{n^{1/4}}
\int_{-\psi(q)}^{\psi(q)}\Re \left( 2\pi i (q+t)^{1/2}e\left(
(q+t)\sqrt{n}+\frac18 \right) \right)dt \ .$$ The factor
$(q+t)^{1/2}$ in the integrand is $q^{1/2}+O(|t|q^{-1/2})$ and so
\begin{equation}
T(q,N)=U(q,N)+O\left( q^{-\frac12}\psi(q)^2N^{3/4} \right) \ ,
\label{tag1.3}
\end{equation}
where
$$U(q,N)=\sum_{1\le n\le N }\frac{r(n)}{n^{1/4}}
\int_{-\psi(q)}^{\psi(q)}\Re \left( 2\pi i q^{1/2}e\left(
(q+t)\sqrt{n}+\frac18 \right) \right)dt \ .$$ Collecting together
the estimates (\ref{tag1.2}) and (\ref{tag1.3}) we have
\begin{eqnarray}
\Delta_0 \!\!\!\!\!\!\!\!\!\! & & \big((q+\psi(q))^2\big)
 \ - \ \Delta_0\big((q-\psi(q))^2\big) \ = \ - \frac1{\pi}\ U(q,N)
\nonumber
\\ & &  \ \ \ \ \ \ + \ O \left( N^{\frac34}q^{-\frac12}\psi(q)^2 +
N^{\frac14}q^{-\frac12}\psi(q) + E_-(q,N) + E_+(q,N) +
qN^{-\frac12} \right)  \ . \label{tag1.4}
\end{eqnarray}

\noindent Let $Q^*$ denote the smallest integer $q$ with $q>Q$.
Then
\begin{eqnarray*}
\sum_{Q<q\le 2Q} \!\!\!\!  U(q,n) = \!\!\!\! \sum_{1\le n\le N}
\frac{r(n)}{n^{1/4}} \int_{-\psi(Q^*)}^{\psi(Q^*)} \Re \left( 2\pi
i e\left( t\sqrt{n}+\frac{1}{8} \right) \sum_{ \substack{ Q<q \le
2Q :  \psi(q)\ge|t|} } q^{\frac{1}{2}} e(q\sqrt{n}) \right)dt \  .
\end{eqnarray*}

\noindent We also have
\begin{eqnarray*}
\sum_{Q<q\le 2Q}  \!\!\!\!\!\!\!\!\!\! &&\sum_{\frac12 Q^2<n<8Q^2}
r(n) \ \min\left( 1, \frac1{|\sqrt{n}-(q\pm \psi(q)|\sqrt{N}}
\right) \\ &\ll&  \ \sum_{Q<q\le 2Q} \sum_{\frac12 Q^2<n<8Q^2}
r(n) \
\min\left( 1, \frac Q{|n-(q\pm \psi(q))^2|\sqrt{N}} \right) \\
&\ll &  \ \sum_{Q<q\le 2Q} \sum_{-q^2<h\le 8Q^2} r(q^2+h) \
\min\left( 1, \frac Q{|h \mp 2q\psi(q) -\psi(q)|\sqrt{N}} \right).
\end{eqnarray*} and
\begin{eqnarray*}
\sum_{Q<q\le 2Q} \!\!\!\!\! &&
\sum_{\substack{4q\psi(q)+2\psi(q)^2<|h|\le 8Q^2 \\ h > -q^2 }}
r(q^2+h) \ \min\left( 1, \frac Q{|h \mp 2q\psi(q)
-\psi(q)|\sqrt{N}} \right)
\\ & &\ll  \ \sum_{Q<q\le 2Q} \ \
\sum_{\substack{4q\psi(q)+2\psi(q)^2<|h|\le 8Q^2 \\ h > -q^2 }}
d(q^2+h) \ \min\left( 1, \frac Q{|h|\sqrt{N}} \right).
\end{eqnarray*}

\noindent Here we observe that $$r(n)\le 4d(n)\le
8\sum_{\substack{l|n  \\ l\le\sqrt{n} }}1 \ . $$

Below we state a bound for the number of solutions of a quadratic
congruence which we use several times over and which is readily
established using elementary facts about such congruences.

\begin{lemma}
 Suppose that $m\in{\Bbb N}$, $h\in{\Bbb Z}$ and define $d_1$ and $d_2$
 uniquely by taking $(m,h)=d_1d_2^2$ where $d_1$ is square free.
 Further let $\rho(m;h)$ denote the number of solutions of the
 congruence $y^2+h\equiv 0\pmod m$
in $y$ modulo $m$.  Then $$\rho(m;h)\le 2d_2d\left(
\frac{m}{d_1d_2^2} \right).$$ \label{bob1}
\end{lemma}

\noindent By Lemma \ref{bob1}
\begin{eqnarray*}
\sum_{Q<q\le 2Q} \!\!\!\!\!\!\!\! & &
\sum_{\substack{4q\psi(q)+2\psi(q)^2<|h|\le 8Q^2 \\ h > -q^2 }}
d(q^2+h) \ \min\left( 1, \frac Q{|h|\sqrt{N}} \right)\\ \\ & \ll &
\sum_{Q<q\le 2Q} \ \sum_{4q\psi(q)+2\psi(q)^2<|h|\le 8Q^2} \ \
\sum_{m|q^2+h, m\ll q} \  \min\left( 1, \frac Q{|h|\sqrt{N}}
\right)\\  \\ &\ll& \sum_{0< d_1d_2^2|j| \ll Q} \ \ \sum_{l \ll
Q/{d_1d_2^2}} \frac{Q^2d(l)}{d_1^2d_2^3l|j|\sqrt{N}} \  \ \ll  \ \
Q^2N^{-\frac12}(\log Q)^3.
\end{eqnarray*}
Hence, by (\ref{tag1.4}),
\begin{eqnarray}
 \sum_{Q<q\le 2Q}
\!\!\!\!\!\!\!\!\!\!\! & & \left(
\Delta_0\big((q+\psi(q))^2\big)-\Delta_0\big((q-\psi(q))^2\big)
\right) \ = \ V(Q,N) \nonumber \\  & + & O\Bigg(
N^{\frac34}\sum_{Q<q\le 2Q}q^{-\frac12}\psi(q)^2 +
N^{\frac14}\sum_{Q<q\le 2Q} q^{-\frac12}\psi(q) \nonumber \\ && \
\ \ \ \ \ \ + \ \ Q^2N^{-\frac12}(\log Q)^3 + F_-(Q,N) + F_+(Q,N)
\Bigg) \label{tag1.5}
\end{eqnarray}
where
$$V(Q,N)=\sum_{1\le n\le N }\frac{r(n)}{n^{1/4}} \Im
\int_{-\psi(Q^*)}^{\psi(Q^*)} 2e\left( t\sqrt{n}+\frac18
\right)\sum_{\substack{Q<q\le 2Q \\ \psi(q)\ge|t|}}q^{1/2}e\left(
q\sqrt{n} \right)dt$$ and
\begin{equation}
F_{\pm}(Q,N)= \sum_{Q<q\le 2Q} \sum_{\substack{|h|\le
4q\psi(q)+2\psi(q)^2 \\ h>-q^2}} r(q^2+h) \min\left( 1, \frac
Q{|h\mp2q\psi(q)-\psi(q)^2|\sqrt{N}} \right). \label{tag1.6}
\end{equation}
When $q>Q$, let
$$W(q)=\sum_{r=Q^*}^qe(r\sqrt{n})$$
and suppose $Q^*\le m\le 2Q$. Then
\begin{eqnarray*}
\sum_{q=Q^*}^m \!\! q^{\frac12}e\left( q\sqrt{n} \right) \!\!\!\!
& = & \!\!\!\! \sum_{q=Q^*}^m \!\! q^{\frac12}(W(q)-W(q-1)) =
-\sum_{q=Q^*}^{m-1}\left( (q+1)^{\frac12}-q^{\frac12}
\right)W(q)+m^{\frac12}W(m)\\
& \ll & Q^{1/2}\min\left( m-Q^*+1,\frac1{\|\sqrt n\|} \right)
 =  Q^{1/2}\min\left( \sum_{q=Q^*}^m 1,\frac1{\|\sqrt n\|} \right).\\
\end{eqnarray*}
We have
$$\int_{-\psi(Q^*)}^{\psi(Q^*)} \min\left( \sum_{\substack{Q<q\le
2Q \\ \psi(q)\ge|t|}} 1,\frac1{\|\sqrt n\|} \right) dt \ \ll \
\min\left( \int_{-\psi(Q^*)}^{\psi(Q^*)}\sum_{\substack{Q<q\le 2Q
\\ \psi(q)\ge|t|}} dt ,\frac{\psi(Q^*)}{\|\sqrt n\|}
\right).$$ Therefore,
$$V(Q,N) \ \ll \  Q^{1/2}\sum_{1\le n\le
N}\frac{r(n)}{n^{1/4}} \ \min\left( \sum_{Q<q\le
2Q}\psi(q),\frac{\psi(Q^*)}{\|\sqrt n\|} \right).$$
\medskip

 Suppose that $1\le m\le \sqrt N+\frac12$, and
consider those $n$ with $\left( m-\frac12 \right)^2<n\le \left(
m+\frac12 \right)^2$. Then $\|\sqrt{n}\|=|\sqrt{n}-m| =
\frac{|n-m^2|}{\sqrt{n}+m} \gg \frac{|n-m^2|}{m}$. Hence, when
$m>1$,
\begin{eqnarray*}
\sum_{(m-1/2)^2<n\le
(m+1/2)^2}\frac{r(n)}{n^{1/4}}\!\!\!\!\!\!\!\!\!\! &&\min\left(
\sum_{Q<q\le 2Q}\psi(q),\frac{\psi(Q^*)}{\|\sqrt n\|} \right)\\
&\ll & \frac{r(m^2)}{m^{\frac12}} \sum_{Q<q\le 2Q}\psi(q) +
m^{\frac12} \sum_{0<|h|\le m} \frac{r(m^2+h)}{|h|} \psi(Q^*) \ .\\
\end{eqnarray*}
The Dirichlet series generating function for $r(m^2)$ is
$$4(1+2^{-s})^{-1}\zeta(s)^2L(s)\zeta(2s)^{-1} \ , $$
where $L(s)$ is the Dirichlet $L$-function formed from the
non-trivial character modulo $4$. Thus
\begin{equation}
\sum_{m\le M}r(m^2) \sim \frac4{\pi}M\log M \label{tag1.7}
\end{equation}
and hence
$$\sum_{m\le M}\frac{r(m^2)}{m^{\frac12}} \ll M^{\frac12} \log
M \ .$$ As in the analysis of $E_{\pm}$ above we have
$$\sum_{2\le m\le M} \sum_{0<|h|\le m} \frac{r(m^2+h)}{|h|}  \ \ll \
M(\log M)^3 \ .$$ Hence
$$V(Q,N) \ \ll  \ N^{\frac14}(\log N)Q^{\frac12} \sum_{Q<q\le 2Q}\psi(q)
+ N^{3/4}(\log N)^3 Q^{\frac12}\psi(Q^*) \ .$$ Hence, assuming
$N\le Q^3$, by (\ref{tag1.5}),
\begin{eqnarray}
\sum_{Q<q\le 2Q}\left( \Delta_0\big((q+\psi(q))^2\big) -
\Delta_0\big((q-\psi(q))^2\big) \right)
\!\!\!\!\!\!\!\!\!\!\!\!\!\!\!\!\!\!\!\!\!\!\!\!\!\!\!\!\!\!
\!\!\!\!\!\!\!\!\!\!\!\!\!\!\!\!\!\!\!\!\!\!\!\!\!\!\!\!\!\!\!\!\!\!\!
& & \nonumber
\\ & \ll &
\!\!\!\!\!\!\!\!\!\!\!\!\!\!\!\!\!\!
 N^{\frac14}(\log N)Q^{\frac12}
\sum_{Q<q\le 2Q}\psi(q) + N^{\frac34}(\log N)^3
Q^{\frac12}\psi(Q^*) \nonumber \\  & \ \ \ \ \ \hspace{15mm} + &
Q^2N^{-\frac12}(\log Q)^3 + F_-(Q,N) + F_+(Q,N) \ .\label{tag1.8}
\end{eqnarray}

\medskip

We now turn our attention to $F_{\pm}$. Were the factor $r(q^2+h)$
not to be present this would be a routine matter. The natural way
to remove it is to consider an application of the Cauchy-Schwarz
inequality.  However one is then dependent on being able to bound
$r(n)^2$, or $d(n)^2$ in terms of the divisors of $n$ of order of
magnitude at most $\sqrt{n}$.  This is readily effected by an
application of a combinatorial lemma.

\begin{lemma}
Let $n\in{\Bbb N}$.  Then there is a divisor $m$ of $n$ such that
$m\le \sqrt{n}$ and $d(n)\le \max(2,d(m)^3)$. \label{bob2}
\end{lemma}

\noindent{\bf Proof. \ }  The conclusion follows at once when $n$
has a prime factor $p$ with $p>\sqrt{n}$. Otherwise choose a
sequence $\{m_j\}$ as follows. Let $m_1$ be the largest divisor of
$n$ not exceeding $\sqrt{n}$.  Then given $m_1,m_2,\ldots,m_j$
with $m_1\ldots m_j|n$ and no $m_k$ exceeding $\sqrt{n}$ choose
$m_{j+1}$ to be the largest divisor of $n/(m_1\ldots m_j)$ not
exceeding $\sqrt{n}$.  It follows that $m_4=1$ since otherwise we
would have $m_1 \, m_2>\sqrt{n}$ and $m_3 \, m_4>\sqrt{n}$. Hence
$n=m_1 \, m_2 \, m_3$ and $d(n)\le d(m_j)^3$ for some $j$. \hfill
$\spadesuit$

\medskip

\noindent By Lemma \ref{bob2},
\begin{eqnarray*}
\sum_{Q<q\le 2Q} \ \ \  \sum_{\substack{0<|h|\le
4q\psi(q)+2\psi(q)^2 \\ h>-q^2}}  \!\!\!\!\!\!\! && d(q^2+h)^2
\\ & &    \ll \ \ \  \sum_{0<|h|\le
4Q\psi(Q^*)+2\psi(Q^*)^2}  \hspace{5mm} \sum_{Q<q\le 2Q}
\hspace{5mm} \sum_{m|q^2+h, \ m\ll Q} \!\!\!\!\! d(m)^6\\
\end{eqnarray*}
and, by Lemma \ref{bob1}, this is
\begin{eqnarray*}
&\ll& \   \sum_{0<d_1d_2^2|j| \le 4Q\psi(Q^*)+2\psi(Q^*)^2} \ \ \
\ \ \sum_{l\ll Q/(d_1d_2^2)} d(ld_1d_2^2)^6 \ \frac{Q}{ld_1d_2} \
d(l)\\ \\  &\ll & \   Q^2 \, \psi(Q^*)\, (\log Q)^{128} \, .
\end{eqnarray*}
We also have
$$\sum_{Q<q\le 2Q} \ \ \  \sum_{|h|\le 4q\psi(q)+2\psi(q)^2} \
\min\left( 1, \frac{Q^2}{|h\mp2q\psi(q)-\psi(q)^2|^2 N} \right)\ll
Q+Q^2N^{-\frac12} \ .$$

\vspace{2mm}

\noindent Hence, by (\ref{tag1.7}) and the Cauchy-Schwarz
inequality, $$F_{\pm}(Q,N) \ \ll \  Q\log Q \ + \
Q^{\frac32}\psi(Q^*)^{\frac12} (\log Q)^{64} \ + \ Q^2
\psi(Q^*)^{\frac12}(\log Q)^{64} N^{-\frac14} \ ,$$ and the
theorem follows from (\ref{tag1.8}). \hfill $\spadesuit$

\vspace{10mm}

RCV: Department of Mathematics, McAllister Building, Pennsylvania
State University, University Park, PA 16802, U.S.A.

\end{appendix}

\end{document}